\documentclass[11pt]{amsart}
\usepackage{amsfonts}
\usepackage{amssymb}
\usepackage[letterpaper, margin=0.8in]{geometry}
\usepackage{fancyhdr}

\begin{document}
\author[]{Ali Enayat}
\title[]{Models of Set Theory: Extensions and Dead-ends}

\begin{abstract}
\noindent This paper is a contribution to the study of extensions of
arbitrary models of $\mathsf{ZF}$ (Zermelo-Fraenkel set theory), with no
regard to countability or well-foundedness of the models involved. Our main
results include the theorems below; in Theorems A and B, $\mathcal{N}$ is
said to be a conservative elementary extension of $\mathcal{M}$ if $\mathcal{%
N}$ elementarily extends $\mathcal{M}$, and the intersection of every $%
\mathcal{N}$-definable set with the universe of $\mathcal{M}$ is $\mathcal{M}
$-definable (parameters allowed). In Theorem B, $\mathsf{ZFC}$\ is the
result of augmenting $\mathsf{ZF}$\ with the axiom of choice.\medskip

\noindent \textbf{Theorem A}. \textit{Every model} $\mathcal{M}$ \textit{of} 
$\mathsf{ZF}+\exists p\left( \mathrm{V}=\mathrm{HOD}(p)\right) $ \textit{has
a conservative elementary extension }$\mathcal{N}$ \textit{that contains an
ordinal above all of the ordinals of} $\mathcal{M}$.\medskip

\noindent \textbf{Theorem B}. \textit{If }$\mathcal{N}$\textit{\ is a
conservative elementary extension of a model} $\mathcal{M}$ \textit{of} $%
\mathsf{ZFC}$, \textit{and} $\mathcal{N}$ \textit{has the same natural
numbers as} $\mathcal{M}$, \textit{then} $\mathcal{M}$ \textit{is cofinal in 
}$\mathcal{N}$.\medskip

\noindent \textbf{Theorem C}. \textit{Every consistent extension of} $%
\mathsf{ZF}$ \textit{has a model} $\mathcal{M}$ \textit{of power} $\aleph
_{1}$ \textit{such that} $\mathcal{M}$ \textit{has no proper end extension
to a model of} $\mathsf{ZF}$.
\end{abstract}

\maketitle

\begin{center}
\textbf{1.~INTRODUCTION\bigskip }
\end{center}

This paper concerns the model theory of Zermelo-Fraenkel ($\mathsf{ZF}$) set
theory. More specifically, it is about various types of extensions of 
\textit{arbitrary} models of $\mathsf{ZF}$, as opposed to countable or
well-founded ones. This is a topic that I have worked on intermittently for
over four decades. Indeed, Theorem C of the abstract answers the following
question that was initially posed in my 1984 doctoral dissertation \cite[%
Question 1.1.12]{Ali-thesis}:\medskip

\noindent \textbf{Question} $\heartsuit $.~\textit{Can every model of }$%
\mathsf{ZF}$\textit{\ be properly \textbf{end extended} to a model of} $%
\mathsf{ZF}$?\footnote{%
The negative answer to the version of this question in which
\textquotedblleft end-extended\textquotedblright\ is strengthened to
\textquotedblleft rank extended\textquotedblright\ was established in the
mid 1980s for models of $\mathsf{ZFC}$ \cite{Ali-TAMS}. Theorems 5.18 of
this paper provide a negative answer to Question $\heartsuit $. This
question was re-posed in \cite{Ali-Bogota}. A variant of the same question
was posed recently by Noah Schweber on MathOverflow \cite{SchweberMO}.}%
\medskip

\noindent Note that if there is a proper class of strongly inaccessible
cardinals, then every \textit{well-founded} model of $\mathsf{ZF}$ can be
properly end extended to a model of $\mathsf{ZF}$. Furthermore, by a
well-known result due to Keisler and Morley \cite{Keisler-Morley}, every
model of $\mathsf{ZF}$ whose class of ordinals has \textit{countable
cofinality} has a proper elementary end extension. The motivation for the
above question emerges from the comparative study of the model theory of $%
\mathsf{ZF}$ and $\mathsf{PA}$ (Peano Arithmetic)$\mathrm{.}$ Prima facie, $%
\mathsf{ZF}$ and $\mathsf{PA}$ are unrelated; after all $\mathsf{PA}$
axiomatizes basic intuitions concerning the familiar arithmetical structure
of natural numbers, whereas $\mathsf{ZF}$ codifies the laws governing
Cantor's mysterious universe of sets, a vastly more complex structure.
However, thanks to a clever coding idea, introduced first by Ackermann%
\textbf{\ }\cite{Ackermann}\textbf{,} $\mathsf{ZF}$ and $\mathsf{PA}$ turn
out to be intimately connected at a formal level: $\mathsf{PA}$ is \textit{%
bi-interpretable} with the theory $\mathsf{ZF}_{\mathrm{fin}}+\mathsf{TC}$,
where $\mathsf{ZF}_{\mathrm{fin}}$ is obtained from $\mathsf{ZF}$ by
replacing the axiom of infinity by its negation, and $\mathsf{TC}$ expresses
\textquotedblleft the transitive closure of every set
exists\textquotedblright .\footnote{%
The role of $\mathsf{TC}$ was elucidated in \cite{ESV} by showing that $%
\mathsf{TC}$ cannot be dropped in this bi-interpretability result. The
details of this bi-interpretability result are worked out by Kaye and Wong
in \cite{Richard-Tin Lok}. Their work shows that indeed the two theories are
definitionally equivalent (aka synonymous), i.e., they have a common
definitional extension.} In this light it is natural to compare and contrast
the model theoretic behavior of models of $\mathsf{ZF}$ and $\mathsf{PA}$ in
order to elucidate the role of the axiom of infinity in Zermelo-Fraenkel set
theory. A central result in the model theory of $\mathsf{PA}$\ is the
MacDowell-Specker Theorem \cite{McDowell-Specker}, and its refinement
(independently) by Phillips \cite{Phillips} and Gaifman \cite{Gaifman-PA},
that states that every model of PA (of any cardinality)\ has a proper
conservative elementary end extension.\ Using the aforementioned
bi-interpretation between $\mathsf{PA}$ and $\mathsf{ZF}$, the refinement of
the MacDowell-Specker Theorem is equivalent to:\medskip

\noindent \textbf{Theorem }$\triangledown $.~\textit{Every model of} $%
\mathsf{ZF}_{\mathrm{fin}}+\mathsf{TC}$ (\textit{of any cardinality}) 
\textit{has a proper conservative elementary end extension}.\medskip

\noindent In the above, $\mathcal{N}$ is said to be a conservative
elementary extension of $\mathcal{M}$ if $\mathcal{N}$ elementarily extends $%
\mathcal{M}$, and the intersection of every $\mathcal{N}$-definable set with 
$M$ is $\mathcal{M}$-definable (parameters allowed). It is known that
Theorem $\triangledown $ becomes false if $\mathsf{ZF}_{\mathrm{fin}}+%
\mathsf{TC}$ is replaced by $\mathsf{ZF}$, but the answer to\ Question $%
\heartsuit $ (Theorem 5.18) remained open until now.\footnote{%
As shown by Keisler and Silver \cite{Keisler-Silver}, if $\kappa $ is the
first inaccessible cardinal, then $\left( \mathrm{V}_{\kappa },\in \right) $
does not have a proper elementary end extension. More generally, the results
of Kaufmann \cite{Kanamori} and the author \cite{Ali-TAMS} showed every
consistent extension of $\mathsf{ZFC}$ has a model of power $\aleph _{1}$
that has no proper rank extension to another model of $\mathsf{ZFC}$. For an
exposition of the general framework for the analogues of the
MacDowell-Specker Theorem for set theory, see \cite{Ali-Bogota}.}\medskip

\noindent The foremost motivation in writing this paper was to present the
rather intricate details of the solution to the above 40-year old question
to my younger self. Another motivation was to publicize what can be
described as a long overdue `Additions and Corrections'\ to my 1987 paper 
\cite{AliJSL1987}, namely Theorems A and B of the abstract (Theorem 3.1 and
Corollary 4.2 of the paper, respectively). \medskip

Section 2 is devoted to preliminaries needed for Sections 3, 4, and 5 in
which the main results are presented. Section 3 contains the proof of
Theorem A of the abstract (and closely related material). Section 4 contains
results on faithful extensions (a generalization of conservative
extensions), including Theorem 4.1 that yields Theorem B of the abstract (as
Corollary 4.2). A key result in Section 4 is Theorem 4.4 which is used in
Section 5 in conjunction with other machinery to establish Theorem C of the
abstract. Section 6 is an \textit{Addendum and Corrigendum} to my
aforementioned 1987 paper \cite{AliJSL1987}. The Appendix presents the proof
of a theorem due to the collective effort of Rubin, Shelah, and Schmerl that
plays a key role in the proof of Theorem C of the abstract.\medskip

The main results of the paper were presented online in February 2025 at the
CUNY set theory seminar. I am grateful to Vika Gitman and Gunter Fuchs and
other participants of the seminar for helpful feedback. Thanks also to Corey
Switzer for illuminating discussions about the forcing argument used in the
proof of Theorem A.1 (in the Appendix), and to Junhong Chen for catching
notational inconsistencies of a previous draft. I am especially indebted to
the anonymous referee for meticulous and extensive feedback that was
instrumental in shaping the paper in its current form.\textbf{\bigskip }

\begin{center}
\textbf{2.~PRELIMINARIES\bigskip }
\end{center}

Here we collect the basic definitions, notations, conventions, and results
that are relevant for the subsequent sections. \medskip

\noindent \textbf{2.1.}~\textbf{Definitions and Basic Facts (}Languages and
theories of sets).\medskip

\begin{enumerate}
\item[\textbf{(a)}] $\mathcal{L}_{\mathrm{set}}=\{=,\in \}$ is the language
of $\mathsf{ZF}$-set theory. We treat $\mathsf{ZF}$ as being axiomatized as
usual, except that instead of including the scheme of replacement among the
axioms of $\mathsf{ZF}$, we include the schemes of separation and
collection, as in \cite[Appendix A]{Chang-Keisler}. Thus, in our setup the
axioms of Zermelo set theory $\mathsf{Z}$ are obtained by removing the
scheme of collection from the axioms of $\mathsf{ZF}$. More generally, for $%
\mathcal{L}\supseteq \mathcal{L}_{\mathrm{set}}$ we construe $\mathsf{ZF}(%
\mathcal{L})$ to be the natural extension of $\mathsf{ZF}$ in which the
schemes of separation and collection are extended to $\mathcal{L}$-formulae;
similarly we will use $\mathsf{Z}(\mathcal{L})$ to denote Zermelo set theory
over $\mathcal{L}$, i.e., as the result of extending $\mathsf{Z}$ with the $%
\mathcal{L}$-separation scheme $\mathsf{Sep}(\mathcal{L})$. Officially
speaking, $\mathsf{Sep}(\mathcal{L})$ consists of the universal closures of $%
\mathcal{L}$-formulae of the form%
\begin{equation*}
\forall v\exists w\forall x(x\in w\longleftrightarrow x\in v\wedge \varphi
(x,\vec{y})).
\end{equation*}%
\medskip Thus $\mathsf{ZF}(\mathcal{L})$ is the result of augmenting $%
\mathsf{Z}(\mathcal{L})$ with the $\mathcal{L}$-collection scheme $\mathsf{%
Coll}(\mathcal{L})$, which consists of the universal closures of $\mathcal{L}
$-formulae of the form:

\begin{equation*}
\left( \forall x\in v\text{ }\exists y\text{\ }\varphi (x,y,\vec{z})\right)
\rightarrow \left( \exists w\text{ }\forall x\in v\text{ }\exists y\in w%
\text{ }\varphi (x,y,\vec{z})\right) .
\end{equation*}%
When \textrm{X}\ is a new predicate, we write $\mathcal{L}_{\mathrm{set}}(%
\mathrm{X})$ instead of $\mathcal{L}_{\mathrm{set}}\cup \{\mathrm{X\}}$, and
similarly we write $\mathsf{ZF}(\mathrm{X})$, $\mathsf{Sep}(\mathrm{X})$, $%
\mathsf{Coll}(\mathrm{X})$, etc. instead of $\mathsf{ZF}(\mathcal{L})$, $%
\mathsf{Sep}(\mathcal{L})$, $\mathsf{Coll}(\mathcal{L})$, etc.
(respectively) when $\mathcal{L}=\mathcal{L}_{\mathrm{set}}(\mathrm{X})$.

\item[\textbf{(b)}] For $n\in \omega $, we employ the common notation ($%
\Sigma _{n},$ $\Pi _{n},$ $\Delta _{n}$) for the Levy hierarchy of $\mathcal{%
L}_{\mathrm{set}}$-formulae, as in the standard references in modern set
theory by Kunen \cite{Kunen-book}, Jech \cite{Jechbook-2003} and Kanamori 
\cite{Kanamori}. In particular, $\Delta _{0}=\Sigma _{0}=\Pi _{0}$
corresponds to the collections of $\mathcal{L}_{\mathrm{set}}$-formulae all
of whose quantifiers are bounded. In addition, we will also consider the
Takahashi hierarchy of formulae ($\Sigma _{n}^{\mathcal{P}},$ $\Pi _{n}^{%
\mathcal{P}},$ $\Delta _{n}^{\mathcal{P}}$) (introduced in \cite{Takahashi},
and further studied by Mathias' \cite{Mathias-MacLane}), where $\Delta _{0}^{%
\mathcal{P}}$ is the smallest class of $\mathcal{L}_{\mathrm{set}}$-formulae
that contains all atomic formulae, and is closed both under Boolean
connectives and under quantification in the form $Qx\in y$ and $Qx\subseteq
y $ where $x$ and $y$ are distinct variables, and $Q$ is $\exists $ or $%
\forall $. The classes $\Sigma _{1}^{\mathcal{P}},\Pi _{1}^{\mathcal{P}}$,
etc. are defined inductively from the class $\Delta _{0}^{\mathcal{P}}$ in
the same way that the formula-classes $\Sigma _{1},\Pi _{1}$, etc. are
defined from $\Delta _{0}$-formulae.

\item[\textbf{(c)}] For $n\in \omega $. and $\mathcal{L}\supseteq \mathcal{L}%
_{\mathrm{set}},$ the Levy hierarchy can be naturally extended to $\mathcal{L%
}$-formulae ($\Sigma _{n}(\mathcal{L}),$ $\Pi _{n}(\mathcal{L}),$ $\Delta
_{n}(\mathcal{L})$), where $\Delta _{0}(\mathcal{L})$ is the smallest family
of $\mathcal{L}$-formulae that contains all atomic $\mathcal{L}$-formulae
and is closed under Boolean operations and bounded quantification.

\item[\textbf{(d)}] $\mathsf{KPR}$ (Kripke-Platek with ranks) is the
subtheory of $\mathsf{ZF}$ whose axioms consist of $\mathsf{KP}$ plus an
axiom that asserts that for all ordinals $\alpha ,$ $\{x:\mathrm{rank}%
(x)<\alpha \}$ exists (where $\mathrm{rank}(x)$ is the usual rank function
on sets). Following recent practice (initiated by Mathias \cite%
{Mathias-MacLane}), the foundation scheme of $\mathsf{KP}$ is only limited
to $\Pi _{1}$-formulae; thus in this formulation $\mathsf{KP}$ can prove the
scheme of $\in $-induction for $\Sigma _{1}$-formulae. In contrast to
Barwise's $\mathsf{KP}$ in \cite{Barwise}, which includes the full scheme of
foundation, our version of $\mathsf{KP}$ is finitely axiomatizable. Thus $%
\mathsf{KPR}$ is finitely axiomatizable.

\item[\textbf{(e)}] $\mathsf{ZBQC}$ is the subtheory of $\mathsf{Z}$,
obtained by augmenting $\mathsf{M}_{0}$ with $\mathsf{Foundation}$ (as a
single axiom), $\mathsf{AC}$ (the axiom of choice), and the axiom of
infinity, where $\mathsf{M}_{0}$ is axiomatized by $\mathsf{Extensionality}$%
, $\mathsf{Pairing}$, $\mathsf{Union}$, $\mathsf{Powerset}$, and $\Delta
_{0} $-$\mathsf{Separation}$.\footnote{$\mathsf{ZBQC}$ was championed by Mac
Lane \cite[p.373]{MacLane-book} as a parsimonious foundation for
mathematical practice. Mathias' \cite{Mathias-MacLane} is an excellent
source of information about $\mathsf{ZBQC}.$} Two distinguished
strengthenings of $\mathsf{ZBQC}$ are $\mathsf{Mac}:=$ $\mathsf{ZBQC+TC}$
(recall that $\mathsf{TC}$ asserts that every set has a transitive closure)%
\textsf{,} and $\mathsf{Most}:=\mathsf{ZBQC}+\mathsf{KP}+\Sigma _{1}$-$%
\mathsf{Separation}$. $\mathsf{ZBQC,}$ $\mathsf{Mac}$, and $\mathsf{Most}$
are finitely axiomatizable.\footnote{%
In the presence of the other axioms of $\mathsf{ZBQC}$, $\Delta _{0}$-$%
\mathsf{Separation}$ is well-known to be equivalent to the closure of the
universe under G\"{o}del-operations, see \cite[Theorem 13.4]{Jechbook-2003}.
This makes it clear that $\mathsf{ZBQC}$ and $\mathsf{Mac}$ are finitely
axiomatizable. Another way to see that $\mathsf{Mac}$ is finitely
axiomatizable is to take advantage of the fact that it supports partial
satisfaction classes (see footnote 14). The availability of partial
satisfaction classes within $\mathsf{Mac}$ makes it clear that the extension 
$\mathsf{Most}$ of $\mathsf{Mac}$ is also finitely axiomatizable.}

\item[\textbf{(f)}] For $\mathcal{L}\supseteq \mathcal{L}_{\mathrm{set}}$,
the \textit{dependent choice scheme, }denoted\textit{\ }$\Pi _{\infty }^{1}$%
\textrm{-}$\mathsf{DC}(\mathcal{L})$, consists of the universal closure of $%
\mathcal{L}$-formulae of the following form:%
\begin{equation*}
\forall x\exists y\ \varphi (x,y,\vec{z})\longrightarrow \forall x\exists f\ %
\left[ f:\omega \rightarrow \mathrm{V},\ f(0)=x,\mathrm{and}\ \forall n\in
\omega \ \varphi (f(n),f(n+1),\vec{z})\right] .
\end{equation*}%
In the presence of $\mathsf{ZF}(\mathrm{\mathcal{L}})\mathrm{,}$ thanks to
the Reflection Theorem (Theorem 2.10), $\Pi _{\infty }^{1}$\textrm{-}$%
\mathsf{DC}\mathrm{(\mathcal{L})}$ is equivalent to the single sentence $%
\mathsf{DC}$ below:%
\begin{equation*}
\forall r\ \forall a\ \left[ 
\begin{array}{c}
\left( \forall x\in a\ \exists y\in a\ \left\langle x,y\right\rangle \in
r\right) \longrightarrow \medskip \\ 
\left( \forall x\in a\ \exists f\ \left( f:\omega \rightarrow \mathrm{V},\
f(0)=x,\mathrm{and}\ \forall n\in \omega \ \left\langle
f(n),f(n+1)\right\rangle \in r\right) \right)%
\end{array}%
\right] .
\end{equation*}%
Note that $\mathsf{DC}$ is provable in $\mathsf{ZFC}$.\medskip
\end{enumerate}

\noindent \textbf{2.2.}~\textbf{Definitions and basic facts.~}(Model
theoretic concepts) In what follows we make the blanket assumption that $%
\mathcal{M}$, $\mathcal{N}$, etc. are $\mathcal{L}$-structures, where $%
\mathcal{L}\supseteq \mathcal{L}_{\mathrm{set}}$. \medskip

\begin{enumerate}
\item[\textbf{(a)}] We follow the convention of using $M$, $M^{\ast },$ $%
M_{0}$, etc.~to denote (respectively) the universes of structures $\mathcal{M%
}$, $\mathcal{M}^{\ast },$ $\mathcal{M}_{0},$ etc. We denote the membership
relation of $\mathcal{M}$ by $\in ^{\mathcal{M}}$; thus an $\mathcal{L}_{%
\mathrm{set}}$-structure $\mathcal{M}$ is of the form $(M,\in ^{\mathcal{M}%
}) $.

\item[\textbf{(b)}] $\mathrm{Ord}^{\mathcal{M}}$ is the class of ordinals in
the sense\ of $\mathcal{M}$, i.e., 
\begin{equation*}
\mathrm{Ord}^{\mathcal{M}}:=\left\{ m\in M:\mathcal{M}\models \mathrm{Ord}%
(m)\right\} ,
\end{equation*}%
where $\mathrm{Ord}(x)$ expresses \textquotedblleft $x$ is transitive and is
well-ordered by $\in $\textquotedblright . More generally, for a formula $%
\varphi (\vec{x})$ with parameters from $M$, where $\vec{m}=\left(
x_{1},\cdot \cdot \cdot ,x_{k}\right) $, we write $\varphi ^{\mathcal{M}}$
for: 
\begin{equation*}
\left\{ \vec{m}\in M^{k}:\mathcal{M\models \varphi }\left( m_{1},\cdot \cdot
\cdot ,m_{k}\right) \right\} .
\end{equation*}%
A subset $D$ of $M^{k}$ is $\mathcal{M}$-definable if it is of the form $%
\varphi ^{\mathcal{M}}$ for some choice of $\varphi .$ We write $\mathbb{%
\omega }^{\mathcal{M}}$ for the set of finite ordinals (i.e., natural
numbers) of $\mathcal{M}$, and $\mathbb{\omega }$ for the set of finite
ordinals in the real world, whose members we refer to as\textit{\
metatheoretic natural numbers. }$\mathcal{M}$ is said to be\textit{\ }$%
\omega $-\textit{standard} if $\left( \mathbb{\omega },\mathbb{\in }\right)
^{\mathcal{M}}\cong \left( \mathbb{\omega },\mathbb{\in }\right) .$ For $%
\alpha \in \mathrm{Ord}^{\mathcal{M}}$ we often use $\mathcal{M}_{\alpha }$
to denote the substructure of $\mathcal{M}$ whose universe is 
\begin{equation*}
M_{\alpha }:=\left\{ m\in M:\mathcal{M}\models m\in \mathrm{V}_{\alpha
}\right\} ,
\end{equation*}%
where $\mathrm{V}_{\alpha }$ is defined as usual as $\{x:\mathrm{rank}%
(x)<\alpha \}$, where $\mathrm{rank}(x)=\sup \{\mathrm{rank}(y)+1:y\in x\}.$%
\medskip
\end{enumerate}

\noindent \textbf{(c) }For $c\in M$, $\mathrm{Ext}_{\mathcal{M}}(c)$ is the $%
\mathcal{M}$-extension of $c$, i.e., $\mathrm{Ext}_{\mathcal{M}}(c):=\{m\in
M:m\in ^{\mathcal{M}}c\}.$\medskip

\noindent \textbf{(d) }We say that $X\subseteq M$ \textit{is coded }(\textit{%
in} $\mathcal{M)}$ if there is some $c\in M$ such that $\mathrm{Ext}_{%
\mathcal{M}}(c)=X.$ \medskip

\noindent \textbf{(e) }Suppose\textbf{\ }$\mathcal{M}\models \mathsf{ZF}$.
We say that $X\subseteq M$ is $\mathcal{M}$-\textit{amenable} if $(\mathcal{M%
},X)\models \mathsf{ZF}(\mathrm{X})$\textsf{;} here $X$ is the
interpretation of $\mathrm{X}$. \medskip

\noindent \textbf{2.3.}~\textbf{Definitions.~}(`geometric shapes'\ of
extensions). Suppose $\mathcal{L}\supseteq \mathcal{L}_{\mathrm{set}}$, $%
\mathcal{M}$ and $\mathcal{N}$ are $\mathcal{L}$-structures, and $\mathcal{M}
$ a \textit{submodel}\footnote{%
The notion of a submodel here is the usual one in model theory; i.e., $%
\mathcal{M}\subseteq \mathcal{N}$ means $M\subseteq N$ and the
interpretation of each nonlogical symbol of $\mathcal{L}$ in $\mathcal{M}$
is the restriction to $M$ of the corresponding of interpretation in $%
\mathcal{N}$ (in particular, $\mathrm{Ext}_{\mathcal{M}}(a)\subseteq \mathrm{%
Ext}_{\mathcal{N}}(a)$ for all $a\in M).$} of $\mathcal{N}$ (written $%
\mathcal{M}\subseteq \mathcal{N)}.$ $\medskip $

\begin{enumerate}
\item[\textbf{(a)}] $\mathcal{M}^{\ast }$ \textit{is the} \textit{convex
hull of} $\mathcal{M}$\textit{\ in} $\mathcal{N}$ if $M^{\ast
}=\bigcup\limits_{a\in M}\mathrm{Ext}_{\mathcal{N}}(a).$

\item[\textbf{(b)}] Suppose $a\in M.$ $\mathcal{N}$ \textit{fixes} $a$ if $%
\mathrm{Ext}_{\mathcal{M}}(a)=\mathrm{Ext}_{\mathcal{N}}(a)$, and $\mathcal{N%
}$ \textit{enlarges} $a$ if $\mathrm{Ext}_{\mathcal{M}}(a)\subsetneq \mathrm{%
Ext}_{\mathcal{N}}(a).\smallskip $

\item[\textbf{(c)}] $\mathcal{N}$ \textit{end extends} $\mathcal{M}$
(written $\mathcal{M}\subseteq _{\mathrm{end}}\mathcal{N}$)$,$ if $\mathcal{N%
}$ fixes every $a\in M$. \ End extensions are also referred to in the
literature as \textit{transitive} extensions, and in the old days as \textit{%
outer} extensions.$\smallskip $

\item[\textbf{(d)}] $\mathcal{N}$ is a \textit{powerset-preserving end
extension} of $\mathcal{M}$ (written $\mathcal{M}\subseteq _{\mathrm{end}}^{%
\mathcal{P}}\mathcal{N}$) if $(i)$ $\mathcal{M}\subseteq _{\mathrm{end}}%
\mathcal{N}$, and $(ii)$ for if $a\in M$, $b\in N$, and $\mathcal{N}\models
(b\subseteq a)$, then $b\in M$.

\item[\textbf{(e)}] $\mathcal{N}$ is a \textit{rank extension of }$\mathcal{N%
}$ (written $\mathcal{M}\subseteq _{\mathrm{rank}}\mathcal{N}$) if $\mathcal{%
N}$ is an end extension of $\mathcal{M}$, and for all $a\in M$ and all $b\in
N\backslash M$, $\mathcal{N}\models \mathrm{rank}(a)\in \mathrm{rank}(b)$.
Here we assume that $\mathcal{M}$ and $\mathcal{N}$ are models of a
sufficient fragment of $\mathsf{ZF}$ (such as $\mathsf{KP}$) in which the
rank function is well-defined. $\mathcal{N}$ is a \textit{topped rank
extension} of $\mathcal{M}$ if $\mathrm{Ord}^{\mathcal{N}}\backslash \mathrm{%
Ord}^{\mathcal{M}}$ has a least element.\footnote{%
Recently rank extensions are also referred to as \textit{top extensions}, we
will not use this terminology as it leads to the expression
\textquotedblleft a topped top extension\textquotedblright .}$\smallskip $

\item[\textbf{(f)}] $\mathcal{N}$ is a \textit{cofinal extension of }$%
\mathcal{M}$ (written $\mathcal{M}\subseteq _{\mathrm{cof}}\mathcal{N)}$ if
for every $b\in N$ there is some $a\in M$ such that $b\in \mathrm{Ext}_{%
\mathcal{N}}(a)$, i.e., $\mathcal{N}$ is the convex hull of $\mathcal{M}$ in 
$\mathcal{N}$.$\smallskip $

\item[\textbf{(g)}] $\mathcal{N}$ is \textit{taller} than $\mathcal{M}$
(written $\mathcal{M}\subseteq _{\mathrm{taller}}\mathcal{N)}$ if there is
some $b\in M$ such that $M\subseteq \mathrm{Ext}_{\mathcal{N}}(b).$ \medskip
\end{enumerate}

\noindent \textbf{2.4.}~\textbf{Definitions.~(}`logical shapes'\ of
extensions). Suppose $\mathcal{L}\supseteq \mathcal{L}_{\mathrm{set}}$, and $%
\mathcal{M}$ and $\mathcal{N}$ are $\mathcal{L}$-structures such that $%
\mathcal{M}$ is a submodel of $\mathcal{N}$.\medskip

\begin{enumerate}
\item[\textbf{(a)}] For a subset $\Gamma $ of $\mathcal{L}$-formulae, $%
\mathcal{M}$ is a $\Gamma $-\textit{elementary submodel of }$\mathcal{N}$
(written $\mathcal{M}\preceq _{\Gamma }\mathcal{N}$) if for all $n$-ary
formulae $\gamma \in \Gamma $ and for all $a_{1},\cdot \cdot \cdot ,a_{n}$
in $M$, $\mathcal{M}\models \gamma (a_{1},\cdot \cdot \cdot ,a_{n})$ iff $%
\mathcal{N}\models \gamma (a_{1},\cdot \cdot \cdot ,a_{n}).$ $\mathcal{M}$
is an \textit{elementary submodel} of $\mathcal{N}$ (written $\mathcal{M}%
\preceq \mathcal{N}$) if $\mathcal{M}\preceq _{\Gamma }\mathcal{N}$ for $%
\Gamma =$ $\mathcal{L}$-formulae. We say that $\mathcal{M}$\ is an \textit{%
elementary submodel} of $\mathcal{N}$ if $\mathcal{M}\preceq _{\Gamma }%
\mathcal{N}$ for the set $\Gamma $ of all $\mathcal{L}$-formulae (written $%
\mathcal{M}\preceq \mathcal{N}$). For a given $\gamma \in \Gamma $ we say
that $\gamma $ is \textit{absolute between} $\mathcal{M}$ \textit{and} $%
\mathcal{N}$ if $\mathcal{M}\preceq _{\{\gamma \}}\mathcal{N}$. As usual, we
write $\mathcal{M}\prec _{\Gamma }\mathcal{N}$ if $\mathcal{M}\preceq
_{\Gamma }\mathcal{N}$ and $\mathcal{M}\subsetneq \mathcal{N}$ (similarly
for $\mathcal{M}\prec \mathcal{N}$)$.\smallskip $

\item[\textbf{(b)}] $\mathcal{N}$ is a \textit{self-extension} of $\mathcal{M%
}$ (written $\mathcal{M}\subseteq _{\mathrm{self}}\mathcal{N}$) if $\mathcal{%
N}$ is \textit{interpretable} in $\mathcal{M}$ (i.e., $\mathcal{N}$ is
isomorphic to a model $\mathcal{N}^{\ast }$ whose universe $N^{\ast }$ of
and the interpretation $R^{\mathcal{N}^{\ast }}$ of each $R\in \mathcal{L}$,
are all $\mathcal{M}$-definable), and additionally, there is an $\mathcal{M}$%
-definable embedding $j$ of $\mathcal{M}$ into $\mathcal{N}^{\ast }$. In
this context, $\mathcal{N}$ is an \textit{elementary self-extension} of $%
\mathcal{M}$ if the map $j$ is an elementary embedding.$\smallskip $

\item[\textbf{(c)}] $\mathcal{N}$ is a \textit{conservative} extension of $%
\mathcal{M}$ (written $\mathcal{M}\subseteq _{\mathrm{cons}}\mathcal{N}$) if
for every $\mathcal{N}$-definable $D$, $M\cap D$ is $\mathcal{M}$-definable.$%
\smallskip $

\item[\textbf{(d)}] For\textbf{\ }$\mathcal{M}\models \mathsf{ZF}$, $%
\mathcal{N}$ is a \textit{faithful} extension of $\mathcal{M}$ (written $%
\mathcal{M}\subseteq _{\mathrm{faith}}\mathcal{N)}$, if for every $\mathcal{N%
}$-definable $D$, $M\cap D$ is $\mathcal{M}$-amenable (as in Definition
2.2(e)).\footnote{%
The notion of a faithful extension was first introduced and studied in the
context of \textit{elementary} extensions in \cite{AliJSL1987}.}\medskip
\end{enumerate}

\noindent \textbf{2.5.}~\textbf{Remark.}~The following are readily
verifiable:\medskip

\begin{enumerate}
\item[\textbf{(a)}] If $\left\langle \mathcal{N}_{k}:k\in \omega
\right\rangle $ is a chain of models extending $\mathcal{M}$ whose union is $%
\mathcal{N}$, and $\mathcal{M}\preceq _{\Pi _{k,\mathrm{cons}}}\mathcal{N}%
_{k}\preceq _{\Pi _{k}}\mathcal{N}_{k+1}$ for each $k\in \omega ,$ then $%
\mathcal{M}\preceq _{\mathrm{cons}}\mathcal{N}$.\footnote{%
This fact readily follows from Theorem 2.11. It comes handy in the first
proof of Theorem 3.1.} Here $\mathcal{M}\preceq _{\Pi _{k,\mathrm{cons}}}%
\mathcal{N}_{k}$ is shorthand for: $\mathcal{M}\preceq _{\Pi _{k}}\mathcal{N}%
_{k}$ and $\mathcal{M}\preceq _{_{\mathrm{cons}}}\mathcal{N}_{k}$.$%
\smallskip $

\item[\textbf{(b)}] $\left( \mathcal{M}\subseteq _{\mathrm{self}}\mathcal{N}%
\right) \Rightarrow \left( \mathcal{M}\subseteq _{\mathrm{cons}}\mathcal{N}%
\right) .\smallskip $

\item[\textbf{(c)}] Assuming $\mathcal{M}\models \mathsf{ZF}$, $\left( 
\mathcal{M}\subseteq _{\mathrm{cons}}\mathcal{N}\right) \Rightarrow \left( 
\mathcal{M}\subseteq _{\mathrm{faith}}\mathcal{N}\right) .$\medskip
\end{enumerate}

\noindent The implications in (b) and (c) are not reversible. More
specifically, (b) is not reversible since, e.g., if $\mathcal{N}$ is a
self-extension of $\mathcal{M},$ then $\left\vert M\right\vert =\left\vert
N\right\vert ,$ but it is possible to build a conservative elementary
extension $\mathcal{N}$ of a model $\mathcal{M}$ with $\left\vert
M\right\vert <\left\vert N\right\vert ,$ e.g., if $\mathcal{U}$ is an $%
\mathcal{M}$-ultrafilter, then the ultrapower formation modulo $\mathcal{U}$
can be iterated along any linear order (of arbitrary cardinality) and
results in a conservative elementary extension of $\mathcal{M}$. In
particular, it is possible to arrange $\mathcal{M}\subseteq _{\mathrm{cons}}%
\mathcal{N}$ such that $M$ is countable but $N$ is uncountable. Another way
to see that $(i)$ is not reversible is to note that by a theorem of Gaifman 
\cite{Gaifmanembeddings} every elementary self-extension of a model $%
\mathcal{M}$ of the $\mathsf{ZF}$ is a cofinal extension.\footnote{%
Indeed Gaifman's theorem is stated for models $\mathcal{M}$ of the fragment $%
\mathsf{ZR}$ of $\mathsf{ZF}$, where $\mathsf{ZR}$ is the result of
extending $\mathsf{Z}$ (Zermelo set theory) with an axiom that every element
is a member of some set of the form $\mathrm{V}_{\alpha }.$} However, by
Theorem 3.1 there are models of $\mathsf{ZF}$ that have taller conservative
elementary extensions. To see that (c) is not reversible, note that if $%
\kappa $ and $\lambda $ are strongly inaccessible cardinals with $\kappa
<\lambda $, then $\left( \mathrm{V}_{\kappa },\in \right) $ is faithfully
extended by $\left( \mathrm{V}_{\lambda },\in \right) $, but since the truth
predicate for $\left( \mathrm{V}_{\kappa },\in \right) $ is parametrically
definable in $\left( \mathrm{V}_{\lambda },\in \right) $, by Tarski's
undefinability of truth $\left( \mathrm{V}_{\lambda },\in \right) $ is not a
conservative extension of $\left( \mathrm{V}_{\kappa },\in \right) $.\medskip

It is also noteworthy that the most common examples of definable extensions
in the context of models of set theory are \textit{elementary} extensions.
However, the Boolean-valued approach to forcing provides a wealth of
definable extensions that are not elementary. More specifically, if $%
\mathcal{M}$ is a model of $\mathsf{ZFC}$, $\mathbb{B}$ is an element of $%
\mathcal{M}$ that is a complete Boolean algebra from the point of view of $%
\mathcal{M}$, and $\mathcal{U}$ is an ultrafilter on $\mathbb{B}$ that is in 
$\mathcal{M}$, then there is an $\mathcal{M}$-definable embedding $j$ of $%
\mathcal{M}$ into the model $\mathcal{M}^{\mathbb{B}}/\mathcal{U}$ (where $%
\mathcal{M}^{\mathbb{B}}$ is the $\mathbb{B}$-valued forcing extension of
the universe, as calculated in $\mathcal{M}$). However, in typical cases $%
\mathcal{M}$ and $\mathcal{M}^{\mathbb{B}}/\mathcal{U}$ are not elementarily
equivalent, e.g., we can always arrange $\mathbb{B}$ so that the truth-value
of the continuum hypothesis changes in the transition between the $\mathcal{M%
}$ and $\mathcal{M}^{\mathbb{B}}/\mathcal{U}$. Furthermore, even though $%
\mathcal{M}$ and $\mathcal{M}^{\mathbb{B}}/\mathcal{U}$ might be
elementarily equivalent in some cases, in general the embedding $j$ is not
elementary if $\mathbb{B}$ is a nontrivial notion of forcing, since the
statement \textquotedblleft there is a filter over $\mathbb{B}$ that meets
all the dense subsets of $\mathbb{B}$\textquotedblright\ holds in $\mathcal{M%
}^{\mathbb{B}}/\mathcal{U}$, but not in $\mathcal{M}$ (assuming $\mathbb{B}$
is a nontrivial notion of forcing).\medskip

\noindent \textbf{2.6.}~\textbf{Remark.}~Suppose $\mathcal{M}$ and $\mathcal{%
N}$ are $\mathcal{L}_{\mathrm{set}}$-structures with $\mathcal{M}\subseteq 
\mathcal{N}\mathrm{.}$ The following are readily verifiable from the
definitions involved:\medskip

\begin{enumerate}
\item[\textbf{(a)}] For each $n\in \omega ,$ $\mathcal{M}\preceq _{\Sigma
_{n}}\mathcal{N}$ iff $\mathcal{M}\preceq _{\Pi _{n}}\mathcal{N}.\smallskip $

\item[\textbf{(b)}] $\mathcal{M}\subseteq _{\mathrm{end}}\mathcal{N}$,
implies $\mathcal{M}\preceq _{\Delta _{0}}\mathcal{N}$; and $\mathcal{M}%
\preceq _{\Delta _{0}}\mathcal{N}$ implies $\mathcal{M}\preceq _{\Delta _{1}}%
\mathcal{N}$.$\smallskip $

\item[\textbf{(c)}] Suppose $\mathcal{M}\preceq _{\Delta _{0}}\mathcal{N}$, $%
\mathcal{N}$ fixes $a\in M$, and there is bijection $f$ in $\mathcal{M}$
between $a$ and some $b\in M$. Then $\mathcal{N}$ fixes $b$ as well (by $%
\Delta _{0}$-elementarity, $f$ remains a bijection between $a$ and $b$ as
viewed from $\mathcal{N}$)$.\smallskip $

\item[\textbf{(d)}] Suppose $\mathcal{M}$ and $\mathcal{N}$ are models of $%
\mathsf{KPR}$, $\mathcal{M}\preceq _{\Delta _{0}^{\mathcal{P}}}\mathcal{N}$,
and $\alpha \in \mathrm{Ord}^{M}$. Then $\mathrm{V}_{\alpha }^{\mathcal{M}}=%
\mathrm{V}_{\alpha }^{\mathcal{N}}$; this follows from the fact that the
formula expressing that $x$ is an ordinal and $y=\mathrm{V}_{x}$ as a $%
\Sigma _{1}^{\mathcal{P}}$-formula in $\mathsf{KPR}$ (using the usual
recursive definition of the $\mathrm{V}_{\alpha }$-hierarchy).$\smallskip $

\item[\textbf{(e)}] If $\mathcal{M}\subseteq _{\mathrm{end}}^{\mathcal{P}}%
\mathcal{N}$, then $\mathcal{M}\preceq _{\Delta _{0}^{\mathcal{P}}}\mathcal{N%
}$. On the other hand, since the formula expressing \textquotedblleft $x$ is
an ordinal and $y=\mathrm{V}_{x}$\textquotedblright\ can be written as a $%
\Sigma _{1}^{\mathcal{P}}$-formula, it is absolute between models $\mathcal{M%
}$ and $\mathcal{N}$ such that $\mathcal{M}\subseteq _{\mathrm{end}}^{%
\mathcal{P}}\mathcal{N}$. Thus for models $\mathcal{M}$ and $\mathcal{N}$ of 
$\mathsf{KPR}$:%
\begin{equation*}
\mathcal{M}\subseteq _{\mathrm{end}}^{\mathcal{P}}\mathcal{N}\Rightarrow 
\mathcal{M}\subseteq _{\mathrm{rank}}\mathcal{N}.
\end{equation*}%
As pointed out by the referee, by \cite[Proposition 6.17]{Mathias-MacLane},
the converse of the above implication also holds provided $\mathcal{M}%
\models \mathsf{Most}$ and $\mathcal{N}\models \mathsf{Mac}$.$\smallskip $

\item[\textbf{(f)}] If $\mathcal{M}\preceq _{\Delta _{0}}\mathcal{N}$, and $%
\mathcal{M}^{\ast }$ is the convex hull of $\mathcal{M}$\textit{\ }in $%
\mathcal{N}$, then $\mathcal{M}\subseteq _{\mathrm{cof}}\mathcal{M}^{\ast
}\subseteq _{\mathrm{end}}\mathcal{N}$. If $\mathcal{M}$ and $\mathcal{N}$
are furthermore assumed to be models of $\mathsf{KPR}$, and $\mathcal{M}%
\preceq _{\Delta _{0}^{\mathcal{P}}}\mathcal{N}$, then (e) above implies
that $\mathcal{M}\subseteq _{\mathrm{cof}}\mathcal{M}^{\ast }\subseteq _{%
\mathrm{rank}}\mathcal{N}$.$\smallskip $

\item[\textbf{(g)}] If both $\mathcal{M}$ and $\mathcal{N}$ are models of $%
\mathsf{KPR}$, and $\mathcal{M}\preceq _{\Delta _{0}}\mathcal{N}$, then the
statement \textquotedblleft $\mathcal{N}$ is taller than $\mathcal{M}$%
\textquotedblright\ is equivalent to \textquotedblleft there is some $\gamma
\in \mathrm{Ord}^{\mathcal{N}}$ that exceeds each $\alpha \in \mathrm{Ord}^{%
\mathcal{M}}$\textquotedblright .$\smallskip $

\item[\textbf{(h)}] Generic extensions of models of $\mathsf{ZF}$ are end
extensions but not rank extensions since they have the same ordinals.
However, by (f) above, powerset preserving end extensions of models of $%
\mathsf{KPR}$ are rank extensions, thus for models $\mathcal{M}$ and $%
\mathcal{N}$ of $\mathsf{KPR}$ we have: 
\begin{equation*}
\mathcal{M}\subseteq _{\mathrm{end}}^{\mathcal{P}}\mathcal{N\Longrightarrow M%
}\subseteq _{\mathrm{rank}}\mathcal{N}.
\end{equation*}%
On the other hand, since $\mathcal{P}(x)=y$ is $\Pi _{1}$-statement, we have:%
\begin{equation*}
\mathcal{M}\preceq _{\Sigma _{1},\mathrm{end}}\mathcal{N}\Longrightarrow 
\mathcal{M}\subseteq _{\mathrm{end}}^{\mathcal{P}}\mathcal{N}.
\end{equation*}%
The converse of the above implication need not be true for arbitrary models $%
\mathcal{M}$ and $\mathcal{N}$ of $\mathsf{KPR}$. The referee has offered
the following counterexample: let $\mathcal{N}$ be a model of $\mathsf{ZFC}$
that is $\omega $-standard but whose $\omega _{1}$ is nonstandard, and let $%
\mathcal{M}$ be the well-founded part of $\mathcal{N}$. Then $\mathcal{M}%
\models \mathsf{KPR}$, and the property of being equinumerous to a von
Neumann ordinal is a $\Sigma _{1}$-property of $\mathcal{P}^{\mathcal{N}%
}(\omega )$ that holds in $\mathcal{N}$ but fails in $\mathcal{M}$. However,
if $\mathcal{M}\models \mathsf{Most}$ and $\mathcal{N}\models \mathsf{Mac}$,
then by \cite[Proposition 6.17]{Mathias-MacLane} the converse of the above
implication holds.$\smallskip $

\item[\textbf{(i)}] As we shall see in Lemma 4.3, faithful end extensions of
models of $\mathsf{ZF}$ are rank extensions.$\smallskip $

\item[\textbf{(g)}] Thanks to the provability of the induction scheme (over
natural numbers) in $\mathsf{ZF}$, if $\mathcal{N}$ is a faithful extension
of a model $\mathcal{N}$ (where $\mathcal{M}$ is a model of $\mathsf{ZF}$),
then $\mathcal{N}$ fixes each element of $\omega ^{\mathcal{M}}$ (but $%
\mathcal{N}$ need not fix $\omega ^{\mathcal{M}}$, e.g., consider the case
when $\mathcal{N}$ is an internal ultrapower of $\mathcal{M}$ modulo a
nonprincipal ultrafilter over $\omega ^{\mathcal{M}}$).\medskip
\end{enumerate}

\noindent \textbf{2.7.~Definition.~}Reasoning within $\mathsf{KP}$, for each
object $a$ in the universe of sets, let $\dot{a}$ be a constant symbol
denoting $a$ (where the map $a\mapsto \dot{a}$ is $\Delta _{1}$, e.g., $\dot{%
a}=\left\langle 3,a\right\rangle $ as in Devlin's monograph \cite%
{Devlin-book-on-L})$.$ Let $\mathrm{Sent}_{\mathcal{L}_{\mathrm{set}}}(x)$
be the $\mathcal{L}_{\mathrm{set}}$-formula that defines the proper class $%
\mathrm{Sent}_{\mathcal{L}_{\mathrm{set}}}$ of sentences in the
language\medskip\ 

\begin{center}
$\mathcal{L}^{+}=$ $\mathcal{L}\cup \{\dot{a}:a\in \mathrm{V}\}$, \medskip
\end{center}

\noindent and let $\mathrm{Sent}_{\mathcal{L}_{\mathrm{set}}}(i,x)$ be the $%
\mathcal{L}_{\mathrm{set}}$-formula that expresses \textquotedblleft $i\in
\omega ,$ $x\in \mathrm{Sent}_{\mathcal{L}_{\mathrm{set}}},$ and $x$ is a $%
\Sigma _{i}$-sentence\textquotedblright . In (a) and (b) below, $\mathcal{M}%
\models \mathsf{KP}$, $S\subseteq M$, and $k\in \omega ^{\mathcal{M}}$.%
\footnote{$\mathsf{KP}$ is chosen here only for convenience; much weaker
theories of sets do the job here.}\medskip

\begin{enumerate}
\item[\textbf{(a)}] Given $k\in \omega $, $S$ is a $\Sigma _{k}$-\textit{%
satisfaction class for} $\mathcal{M}$ if $\left( \mathcal{M},S\right)
\models \mathrm{Sat}(k,\mathrm{S})$, where $\mathrm{Sat}(k,\mathrm{S})$ is
the universal generalization of the conjunction of the axioms $(i)$ through $%
(iv)$ below. We assume that first order logic is formulated using only the
logical constants $\left\{ \lnot ,\vee ,\exists \right\} .$\medskip

\begin{enumerate}
\item[$(i)$] $\left[ \left( \mathrm{S}\left( \dot{x}=\dot{y}\right)
\leftrightarrow x=y\right) \wedge \left( \mathrm{S}\left( \dot{x}\in \dot{y}%
\right) \leftrightarrow x\in y\right) \right] .$\medskip

\item[$(ii)$] $\left[ \mathrm{Sent}_{\mathcal{L}_{\mathrm{set}}}(k,\varphi
)\wedge \left( \varphi =\lnot \psi \right) \right] \rightarrow \left[ 
\mathrm{S}(\varphi )\leftrightarrow \lnot \mathrm{S}\mathsf{(}\psi \mathsf{)}%
\right] \mathsf{.}$\medskip

\item[$(iii)$] $\left[ \mathrm{Sent}_{\mathcal{L}_{\mathrm{set}}}(k,\varphi
)\wedge \left( \varphi =\psi _{1}\vee \psi _{2}\right) \right] \rightarrow %
\left[ \mathrm{S}(\varphi )\leftrightarrow \left( \mathrm{S}(\psi _{1})\vee 
\mathrm{S}(\psi _{2})\right) \right] \mathsf{.}$\medskip

\item[$(iv)$] $\left[ \mathrm{Sent}_{\mathcal{L}_{\mathrm{set}}}(k,\varphi
)\wedge \left( \varphi =\exists v\ \psi (v)\right) \right] \rightarrow \left[
\mathrm{S}(\varphi )\leftrightarrow \exists x\ \mathrm{S}\mathsf{(}\psi 
\mathsf{(}\dot{x}\mathsf{))}\right] .$\medskip
\end{enumerate}

\item[\textbf{(b)}] $S$ is a \textit{full} \textit{satisfaction class for} $%
\mathcal{M}$ if $\left( \mathcal{M},S\right) \models \mathrm{Sat}(k,\mathrm{S%
})$ for every $k\in \omega $.\medskip
\end{enumerate}

\noindent \textbf{2.8.~Tarski's Definability and Undefinability of Truth
Theorems}.~\textit{Suppose} $\mathcal{M}$ \textit{is an} $\mathcal{L}_{%
\mathrm{set}}$-\textit{structure.}\medskip

\begin{enumerate}
\item[\textbf{(a)}] (Definability/Codability of Truth) \textit{If} $\mathcal{%
M}$ \textit{is a structure coded as an element }$m$ \textit{of a model }$%
\mathcal{N}$ \textit{of} \textit{a sufficiently strong}\footnote{%
There are two canonical fragments of $\mathsf{ZF}$ that are `sufficiently
strong' for this purpose, namely:
\par
\begin{enumerate}
\item[(1)] $\mathsf{KP}$\ (see Definition 2.1(d)), as shown by Friedman, Lu,
and Wong in \cite[Lemma 4.1]{Friedman at al}.
\par
\item[(2)] The fragment $\mathsf{M}_{0}$ of $\mathsf{Z}$ (see Definition
2.1(e)), as shown by Mathias \cite[Proposition 3.10]{Mathias-MacLane}.
\end{enumerate}
\par
\noindent Note, however, that much weaker systems suffice if one only wishes
to have a set theory within which the Tarskian satisfaction relation of
every internal set structure is \textit{definable}, as opposed to: \textit{%
coded as a set}. One such weak system is $\mathsf{DS}$ (for `Devlin
strengthened'), which is shown by Mathias \cite[Proposition 10.37]%
{Mathias-fixingDevlin} to be capable of defining the Tarskian satisfaction
predicate for set structures.}\textit{\ fragment of }$\mathsf{ZF}$, \textit{%
then there is }$s\in N$ \textit{such that }$\mathcal{N}$\textit{\ views }$s$ 
\textit{as the Tarskian satisfaction relation on }$m$\textit{, written }$%
\mathcal{N}\models \mathrm{sat}(s,m).$\footnote{%
In other words, within $\mathcal{N}$, $s$ is the set consisting of ordered
pairs $(\varphi (\vec{x}),\vec{a})$ such that\ $m\models \varphi \left( \vec{%
x}/\vec{a}\right) .$} \textit{Moreover}, $\mathrm{sat}(x,y)$\textit{\ is} $%
\Delta _{1}$ \textit{in} $\mathcal{N}.$\medskip

\item[\textbf{(b)}] (Undefinability of Truth) \textit{If} $S$ \textit{is a
full satisfaction class for} $\mathcal{M}$, \textit{then} $S$ \textit{is not}
$\mathcal{M}$-\textit{definable.} \medskip
\end{enumerate}

It is a well-known result of Levy that if $\mathcal{M}\models \mathsf{ZF}%
\mathrm{,}$ then there is a $\Delta _{0}$-satisfaction class for\textit{\ }$%
\mathcal{M}$\textit{\ }that is\textit{\ }definable in $\mathcal{M}$ both by
a $\Sigma _{1}$-formula and a $\Pi _{1}$-formula (see \cite[p.~186]%
{Jechbook-2003} for a proof). This makes it clear that for each $n\geq 1,$
there is a $\Sigma _{n}$-satisfaction class for\textit{\ }$\mathcal{M}$%
\textit{\ }that is\textit{\ }definable in $\mathcal{M}$ by a $\Sigma _{n}$%
-formula. Levy's result extends to models of $\mathsf{ZF}\mathrm{(}\mathcal{L%
}\mathrm{)}$ if $\mathcal{L}$\ is finite as follows: \medskip

\noindent \textbf{2.9.~Levy's Partial Definability of Truth Theorem.}~%
\textit{Let} $\mathcal{L}$ \textit{be a finite extension of} $\mathcal{L}_{%
\mathrm{set}}$. \textit{For each} $n\in \omega $ \textit{there is an} $%
\mathcal{L}$-\textit{formula} $\mathrm{Sat}_{\Sigma _{n}(\mathcal{L})}$ 
\textit{such that for all models} $\mathcal{M}$ \textit{of} \textit{a
sufficiently strong}\footnote{%
There are two canonical fragments of $\mathsf{ZF}$ that are `sufficiently
strong' for this purpose: namely: $\mathsf{KP}$, and $\mathsf{M}_{0}+\mathsf{%
TC},$ see, e.g., Definitions 2.9 and 2.10 of McKenzie's \cite{Zach
collection} for the case of $\mathsf{KP}$, a similar construction works for $%
\mathsf{M}_{0}+\mathsf{TC}$. What is needed in both cases is $\mathsf{TC}$,
plus the ability of the theory to define the Tarskian satisfaction predicate
for set structures.} $\mathsf{ZF}\mathrm{(}\mathcal{L}\mathrm{)}$, \textrm{%
Sat}$_{\Sigma _{n}(\mathcal{L})}^{\mathcal{M}}$ \textit{is a} $\Sigma _{n}(%
\mathcal{L})$-\textit{satisfaction class for} $\mathcal{M}$. \textit{%
Furthermore, for} $n\geq 1$, $\mathrm{Sat}_{\Sigma _{n}(\mathcal{L})}$ 
\textit{is equivalent to a} $\Sigma _{n}(\mathcal{L})$-\textit{formula} (%
\textit{provably in} $\mathsf{ZF}(\mathcal{L})$\textrm{)}$\mathrm{.}$\medskip

The Reflection Theorem is often formulated as a theorem scheme of $\mathsf{ZF%
}$ (e.g., as in \cite{Jechbook-2003}), the proof strategy of the Reflection
Theorem applies equally well to $\mathsf{ZF}(\mathcal{L})$ for finite
extensions $\mathcal{L}$ of $\mathcal{L}_{\mathrm{set}}$. \medskip

\noindent \textbf{2.10.~Montague-Vaught-Levy Reflection Theorem.}~\textit{Let%
} $\mathcal{L}$ \textit{be a finite extension of} $\mathcal{L}_{\mathrm{set}%
}.$ \textit{For each} $n\in \omega $, $\mathsf{ZF}(\mathcal{L})$\textit{\
proves that there are arbitrarily large ordinals} $\gamma $ \textit{such that%
} \textit{the submodel of the universe determined by }\textrm{$V$}$_{\gamma
} $ \textit{is a} $\Sigma _{n}(\mathcal{L})$-\textit{elementary submodel of
the universe}.\footnote{%
This is expressible with the help of $\mathrm{Sat}_{\Sigma _{n}(\mathcal{L}%
)}.$}\medskip

\noindent \textbf{2.11.~}$\mathbf{\Sigma }_{\mathbf{n}}$-\textbf{Elementary
Chains Theorem.}~\textit{Suppose} $n\in \omega ,$ $\mathcal{L}\supseteq 
\mathcal{L}_{\mathrm{set}}$, $(I,<_{I})$ \textit{is a linear order, and} $%
\left\{ \mathcal{M}_{i}:i\in I\right\} $ \textit{is a collection of} $%
\mathcal{L}$-\textit{structures such that }$\mathcal{M}_{i}\ \mathcal{%
\preceq }_{\Sigma _{n}(\mathcal{L})}\ \mathcal{M}_{i^{\prime }}$ \textit{%
whenever }$i<_{I}i^{\prime }.$ \textit{Then we have}:\medskip

\begin{center}
$\mathcal{M}_{i}\ \mathcal{\preceq }_{\Sigma _{n}(\mathcal{L})}\
\bigcup\limits_{i\in I}\mathcal{M}_{i}.$\medskip
\end{center}

The following is a special case of a result of Gaifman \cite[Theorem 1, p.54]%
{Gaifmanembeddings}. In the statement and the proof, $\mathcal{M}^{\ast
}\prec _{_{\Delta _{0}(\mathcal{L}),\mathrm{end}}}\mathcal{N}$ is shorthand
for the conjunction of $\mathcal{M}^{\ast }\prec _{\Delta _{0}(\mathcal{L)}}%
\mathcal{N}$ and $\mathcal{M}^{\ast }\prec _{\mathrm{end}}\mathcal{N}$%
.\medskip

\noindent \textbf{2.12.}~\textbf{Gaifman Splitting Theorem.}~\textit{Let} $%
\mathcal{L}\supseteq \mathcal{L}_{\mathrm{set}}.$ \textit{Suppose }$\mathcal{%
M}\models \mathsf{ZF}(\mathcal{L})$, $\mathcal{N}$ \textit{is an} $\mathcal{L%
}$-\textit{structure with }$\mathcal{M}\prec _{\Delta _{0}(\mathcal{L})}%
\mathcal{N}$, \textit{and} $\mathcal{M}^{\ast }$ \textit{is the convex hull
of} $\mathcal{M}$ \textit{in} $\mathcal{N}$. \textit{Then the following hold}%
:\medskip

\begin{enumerate}
\item[\textbf{(a)}] $\mathcal{M}\preceq _{\mathrm{cof}}\mathcal{M}^{\ast
}\preceq _{\Delta _{0}(\mathcal{L}),\mathrm{end}}\mathcal{N}$.\medskip

\item[\textbf{(b)}] $\mathcal{M\prec N}\Rightarrow \mathcal{M}\preceq _{%
\mathrm{cof}}\mathcal{M}^{\ast }\preceq _{\mathrm{end}}\mathcal{N}$.\medskip
\end{enumerate}

\noindent \textbf{Proof.} To simplify the notation, we present the proof for 
$\mathcal{L}=\mathcal{L}_{\mathrm{set}}$ since the same reasoning handles
the general case. In our proof we use the notation $\varphi ^{v}$, where $%
\varphi $ is an $\mathcal{L}_{\mathrm{set}}$-formula and $v$ is a parameter,
to refer to the $\Delta _{0}$-formula obtained by restricting all the
quantifiers of $\varphi $ to the elements of $v$. A straightforward
induction on the complexity of formulae shows that:\medskip

\noindent (1) If $\mathcal{K}$ is an $\mathcal{L}_{\mathrm{set}}$-structure,
and $m\in K$, then $\mathcal{K}\models \varphi ^{m}$ iff $(m,\in )^{\mathcal{%
K}}\models \varphi $. \medskip

\noindent Note that in the above, $\varphi $ is allowed to contain
parameters in $\mathrm{Ext}_{\mathcal{K}}(m).$ Suppose $\mathcal{M}\prec
_{\Delta _{0}}\mathcal{N},$ where $\mathcal{M}\models \mathsf{ZF.}$ Note
that we are not assuming that $\mathcal{N}\models \mathsf{ZF}$.\textsf{\ }%
Putting (1) together with $\mathcal{M}\prec _{\Delta _{0}}\mathcal{N}$
implies that $(m,\in )^{\mathcal{M}}\preceq (m,\in )^{\mathcal{M}^{\ast }}$
for all $m\in M.$ On the other hand, the definition of $\mathcal{M}^{\ast }$
makes it clear that $(m,\in )^{\mathcal{N}}=(m,\in )^{\mathcal{M}^{\ast }}$
if $m\in M.$ Putting all this together, we have: \medskip

\noindent (2) \ \ For all $m\in M$, $(m,\in )^{\mathcal{M}}\preceq (m,\in )^{%
\mathcal{M}^{\ast }}=(m,\in )^{\mathcal{N}}$.\medskip

\noindent It is easy to see, using the definition of $\mathcal{M}^{\ast }$,
that $\mathcal{M}\subseteq _{\mathrm{cof}}\mathcal{M}^{\ast }\subseteq _{%
\mathrm{end}}\mathcal{N}$, and in particular $\mathcal{M}^{\ast }\preceq
_{\Delta _{0}}\mathcal{N}$ by Remark 2.6(b). Thus the proof of (a) is
complete once we show that $\mathcal{M}\preceq _{\Sigma _{n}}\mathcal{M}%
^{\ast }$ for each $n\in \omega .$ Towards this aim, given $n\in \omega $,
we can invoke the Reflection Theorem in $\mathcal{M}$ to get hold of some $%
U\subseteq \mathrm{Ord}^{\mathcal{M}}$ such that:\medskip

\noindent (3) \ \ $(\mathcal{M},U)$ satisfies \textquotedblleft $U$ is
unbounded in $\mathrm{Ord},$ and $\left[ \left( \mathrm{V}_{\alpha },\in
\right) \prec _{\Sigma _{n}}\left( \mathrm{V},\in \right) \right] $ for
every $\alpha \in U$\textquotedblright $.$\medskip

\noindent For each $\alpha \in U$ let $m_{\alpha }\in M$ such that $%
m_{\alpha }=\mathrm{V}_{\alpha }^{\mathcal{M}}.$ By (3), we have:\medskip

\noindent (4) \ \ If $\alpha ,\beta \in U$ with $\alpha <\beta $, then $%
(m_{\alpha },\in )^{\mathcal{M}}\prec _{\Sigma _{n}}(m_{\beta },\in )^{%
\mathcal{M}}.$\medskip

\noindent Next, we observe that if $a,b$ are in $M$ with $(a,\in )^{\mathcal{%
M}}\prec _{\Sigma _{n}}(b,\in )^{\mathcal{M}}$, then for each $k$-ary $%
\Sigma _{n}$-formula $\sigma (\vec{x})$, $\mathcal{M}$ satisfies the $\Delta
_{0}$-statement $\delta _{\sigma }(a,b),$ where: \medskip

\begin{center}
$\delta _{\sigma }(a,b):=\forall \vec{x}\in a\ \left[ \sigma ^{a}(\vec{x}%
)\leftrightarrow \sigma ^{b}(\vec{x})\right] $,\medskip
\end{center}

\noindent where the prefix $\forall \vec{x}\in a$ is shorthand for $\forall
x_{0}\in a\cdot \cdot \cdot \forall x_{k-1}\in a.$ Therefore, by putting (4)
and (2) together with the assumption $\mathcal{M}\prec _{\Delta _{0}}%
\mathcal{N}$ we can conclude:\medskip

\noindent (5) \ \ If $\alpha ,\beta \in U$ with $\alpha <\beta $, then $%
(m_{\alpha },\in )^{\mathcal{M}^{\ast }}\prec _{\Sigma _{n}}(m_{\beta },\in
)^{\mathcal{M}^{\ast }}.$ \medskip

\noindent On the other hand, since $U$ is unbounded in $\mathrm{Ord}^{%
\mathcal{M}}$, the definition of $m_{\alpha }$, together with the fact that $%
\mathcal{M}\subseteq _{\mathrm{cof}}\mathcal{M}^{\ast }$, imply: \medskip

\noindent (6) $\ \ \mathcal{M}=\bigcup\limits_{\alpha \in U}(m_{\alpha },\in
)^{\mathcal{M}}$ and $\mathcal{M}^{\ast }=\bigcup\limits_{\alpha \in
U}(m_{\alpha },\in )^{\mathcal{M}^{\ast }}.$\medskip

\noindent Thanks to (2), (5), (6), and Theorem 2.11 ($\Sigma _{n}$%
-Elementary Chains Theorem) we can conclude:\medskip

\noindent (7) \ \ If $\alpha \in U$, then $(m_{\alpha },\in )^{\mathcal{M}%
}\preceq (m_{\alpha },\in )^{\mathcal{M}^{\ast }}\preceq _{\Sigma _{n}}%
\mathcal{M}^{\ast }\mathrm{.}$\medskip

\noindent We can conclude that $\mathcal{M}\preceq _{\Sigma _{n}}\mathcal{M}%
^{\ast }$ by (4), (7), and Theorem 2.11. This concludes the proof of (a).
\medskip

To prove (b), suppose $\mathcal{M}\prec \mathcal{N}.$ With (a) at hand, we
only need to show that $\mathcal{M}^{\ast }\preceq _{\Sigma _{n}}\mathcal{N}$
for each $n\in \omega $. We observe that if $\alpha \in U,$ then by (5) $%
(m_{\alpha },\in )^{\mathcal{M}^{\ast }}\prec _{\Sigma _{n}}\mathcal{M}%
^{\ast },$ and by (2), $(m_{\alpha },\in )^{\mathcal{M}^{\ast }}=(m_{\alpha
},\in )^{\mathcal{N}}.$ On the other hand, the assumption $\mathcal{M}\prec 
\mathcal{N}$ together with (3) makes it clear that if $\alpha \in U,$ $%
(m_{\alpha },\in )^{\mathcal{N}}\prec _{\Sigma _{n}}\mathcal{N}$. Thus $%
\mathcal{M}^{\ast }\preceq _{\Sigma _{n}}\mathcal{N}$ by Theorem 2.11.\hfill 
$\square $\medskip

\noindent \textbf{2.13.~Remark.~}In general, if $\mathcal{M}\subseteq 
\mathcal{N}$, where $\mathcal{M}$ and $\mathcal{N}$ are both models of $%
\mathsf{ZF}$, the condition \textquotedblleft $\mathrm{Ord}^{\mathcal{M}}$
is cofinal in $\mathrm{Ord}^{\mathcal{N}}$\textquotedblright\ does not imply
that $\mathcal{M}$ is cofinal in $\mathcal{N}$ (for example consider $%
\mathcal{M}$ and $\mathcal{N}$ where $\mathcal{N}$\ is a set-generic
extension of $\mathcal{M}$). However, the conjunction of $\Sigma _{1}$%
-elementarity and \textquotedblleft $\mathrm{Ord}^{\mathcal{M}}$ is cofinal
in $\mathrm{Ord}^{\mathcal{N}}$\textquotedblright , where $\mathcal{M}%
\models \mathsf{ZF}$ (or even $\mathsf{KPR}$) and $\mathcal{N}$ is an $%
\mathcal{L}_{\mathrm{set}}$-structure implies that $\mathcal{M}$ is cofinal
in $\mathcal{N}$, since the statement $(y\in \mathrm{Ord}\wedge x=\mathrm{V}%
_{y})$ can be written as a $\Pi _{1}$-statement within $\mathsf{KPR}$, and
therefore it is absolute between $\mathcal{M}$ and $\mathcal{N}$ if $%
\mathcal{M}\prec _{\Sigma _{1}}\mathcal{N}.$\footnote{%
As noted in \cite[Lemma 0.2, p5]{Kanamori}, the rank function is $\Delta
_{1}^{\mathsf{ZF}}$, and this fact can be used to show that $(y\in \mathrm{%
Ord}\wedge x=\mathrm{V}_{y})$ can be written as a $\Pi _{1}$-statement
within $\mathsf{ZF}$; an examination of the proof makes it clear that $%
\mathsf{KPR}$ suffices for this purpose.} This is how elementary embeddings
between inner models are usually formalized within $\mathsf{ZF}$, as in
Kanamori's monograph \cite[Proposition 5.1]{Kanamori}.\textbf{\medskip }

\noindent Recall that for models of $\mathsf{ZF}$, the $\mathcal{L}_{\mathrm{%
set}}$-sentence $\exists p\left( \mathrm{V}=\mathrm{HOD}(p)\right) $
expresses: \textquotedblleft there is some $p$ such that every set is first
order definable in some structure of the form $\left( \mathrm{V}_{\alpha
},\in ,p,\beta \right) _{\beta <\alpha }$ with $p\in \mathrm{V}_{\alpha }$%
\textquotedblright . The following result is classical.\medskip

\noindent \textbf{2.14.~Myhill-Scott Theorem \cite{Myhill-Scott}.~}\textit{%
The following statements are equivalent for} $\mathcal{M}\models \mathsf{ZF}$%
\textrm{:}\textbf{\medskip }

\begin{enumerate}
\item[\textbf{(a)}] $\mathcal{M}\models \exists p\left( \mathrm{V}=\mathrm{%
HOD}(p)\right) .$\textbf{\medskip }

\item[\textbf{(b)}] $\mathcal{M}$ \textit{carries a} \textit{definable
global well-ordering, i.e.,} \textit{for some} $p\in M$ \textit{and some
set-theoretic formula} $\varphi (x,y,z)$, $\mathcal{M}$ \textit{satisfies }%
\textquotedblleft $\varphi (x,y,p)$ \textit{well-orders the universe}%
\textquotedblright .\textbf{\medskip }

\item[\textbf{(c)}] $\mathcal{M}$ \textit{carries a} \textit{definable
global choice function, i.e., for some} $p\in M$ \textit{and some
set-theoretic formula} $\psi (x,y,z)$, $\mathcal{M}$ \textit{satisfies }%
\textquotedblleft $\psi (x,y,p)$ \textit{is the graph of a global choice
function}\textquotedblright .
\end{enumerate}

\textbf{\bigskip }

\begin{center}
\textbf{3.~TALLER CONSERVATIVE ELEMENTARY EXTENSIONS}\bigskip
\end{center}

\noindent \textbf{3.1.~Theorem.}~\textit{Every model }$\mathcal{M}$\textit{\
of }$\mathsf{ZF}+\exists p\left( \mathrm{V}=\mathrm{HOD}(p)\right) $ (%
\textit{of any cardinality}) \textit{has a conservative elementary extension 
}$\mathcal{N}$\textit{\ such that }$\mathcal{N}$\textit{\ is taller than }$%
\mathcal{M}$\textit{.}\medskip

We shall present two proofs of Theorem 3.1, the first one is based on a
class-sized syntactic construction taking place within a model of set
theory; it was inspired by Kaufmann's proof of the MacDowell-Specker theorem
using the Arithmetized Completeness Theorem, as described by Schmerl \cite%
{Schmerl Matt proof}.\footnote{%
As noted by Phillips \cite{Phillips} and Gaifman \cite{Gaifman-PA}, the
McDowell-Specker proof lends itself to a fine-tuning that ensures that the
conservative elementary extension $\mathcal{N}$ of a given model $\mathcal{M}
$ has the minimality property, i.e., there is no $\mathcal{K}$ such that $%
\mathcal{M}\prec \mathcal{K}\prec \mathcal{N}$. However, Kaufmann's slick
proof does not seem to lend itself to this embellishment.} The second one is
based on a model-theoretic construction that is reminiscent of the original
ultrapower proof\footnote{%
The original proof by MacDowell and Specker of their result uses the concept
of a \textquotedblleft finitely additive 2-valued measure\textquotedblright
, which is an alternative formulation of the notion of an ultrafilter.} of
the MacDowell-Specker theorem \cite{McDowell-Specker}. The first proof is
short and devilish; the second proof is a bit longer but is more transparent
due to its combinatorial flavor; it also lends itself to a refinement, as
indicated in Remark 3.3 and Theorem 3.4.\medskip

\noindent \textbf{First proof of Theorem 3.1}.~We need the following two
Facts 1 and 2 below.\medskip

\noindent FACT 1: \textit{Suppose} $\mathcal{M}$ \textit{is a model of }$%
\mathsf{ZF}$ \textit{that carries an }$\mathcal{M}$\textit{-definable global
well-ordering}, \textit{and} $T$ \textit{is an }$\mathcal{M}$-\textit{%
definable class of first order sentences such that }$\mathcal{M}$ \textit{%
satisfies} \textquotedblleft $T$\ \textit{is a consistent first order theory}%
\textquotedblright . \textit{Then} \textit{there is a model }$\mathcal{N}%
\models T^{\mathrm{st}}$ \textit{such that the elementary diagram of} $%
\mathcal{N}$ \textit{is }$\mathcal{M}$\textit{-definable, where} $T^{\mathrm{%
st}}$ \textit{is the collection of sentences in} $T$ \textit{with standard
shape\footnote{%
More specifically, $\varphi \in T^{\mathrm{st}}$ iff $\varphi \in T$ and
there is some standard formula $\varphi ^{\ast }$ such that $\mathcal{M}$
believes that $\varphi $ is the result of substituting constants for the
free variables of $\varphi ^{\ast }.$}}. \medskip

\noindent Proof of Fact 1. Since $\mathcal{M}$ has a definable global
well-ordering, the Henkin proof of the completeness theorem of first order
logic can be applied within $\mathcal{M}$ to construct a Henkinized complete
extension $T^{\mathrm{Henkin}}$ of $T$ (in a language extending the language
of $T$\ by class-many new constant symbols) such that $T$\ is definable in $%
\mathcal{M}$. This in turn allows $\mathcal{M}$ to define $\mathcal{N}$ by
reading it off $T^{\mathrm{Henkin}},$ as in the usual Henkin proof of the
completeness theorem.\medskip

\noindent FACT 2: \textit{If}$\mathcal{M}\models \mathsf{ZF}$, \textit{then} 
\textit{for each} $n\in \omega $ $\mathcal{M}\models \mathrm{Con}(\mathrm{Th}%
_{\Pi _{n}}($\textrm{$V$}$,\in ,\dot{a})_{a\in \mathrm{V}}$); \textit{here} $%
\mathrm{Con}(\mathrm{X})$ \textit{expresses the formal consistency of} $%
\mathrm{X}$, \textit{and} $\mathrm{Th}_{\Pi _{n}}($\textrm{$V$}$,\in ,\dot{a}%
)_{a\in \mathrm{V}}$ \textit{is the} $\Pi _{n}$-\textit{fragment of the
elementary} \textit{diagram of the universe }(\textit{which is a definable
class}, \textit{using }$\lnot \mathrm{Sat}_{\Sigma _{n}}$, \textit{where} $%
\mathrm{Sat}_{\Sigma _{n}}$ \textit{is as in Theorem} 2.9). \medskip

\noindent Proof of Fact 2. This is an immediate consequence of the
Reflection Theorem 2.10.\medskip

\noindent Starting with a model $\mathcal{M}$ of $\mathsf{ZF}+\exists
p\left( \mathrm{V}=\mathrm{HOD}(p)\right) ,$ by Theorem 2.14, $\mathcal{M}$
carries a global definable well-ordering. We will construct an increasing
sequence of $\mathcal{L}_{\mathrm{set}}$-structures $\left\langle \mathcal{N}%
_{k}:k\in \omega \right\rangle $ that satisfies the following properties for
each $k\in \omega $:\medskip

\begin{enumerate}
\item[(1)] $\mathcal{N}_{0}=\mathcal{M},$\medskip

\item[(2)] $\mathcal{M}\preceq _{\Pi _{k+2}}\mathcal{N}_{k}\preceq _{\Pi
_{k+1}}\mathcal{N}_{k+1}$.\medskip

\item[(3)] There is some $\alpha \in \mathrm{Ord}^{\mathcal{N}_{1}}$ that is
above each $\beta \in \mathrm{Ord}^{\mathcal{M}}$, thus $\mathcal{M}\preceq
_{\mathrm{taller}}\mathcal{N}_{k}$ for all $k\geq 1$.\medskip

\item[(4)] $\mathcal{N}_{k}$ is a self-extension of $\mathcal{M}$, and in
particular $\mathcal{M}\preceq _{\mathrm{cons}}\mathcal{N}_{k}$.\footnote{%
See Definition 2.4(b) and Remark 2.5(a).}\medskip
\end{enumerate}

\noindent In other words: \medskip

\begin{center}
$\mathcal{M=N}_{0}\preceq _{\Pi _{3},\mathrm{taller}}\mathcal{N}_{1}\preceq
_{\Pi _{2}}\mathcal{N}_{2}\preceq _{\Pi _{3}}\mathcal{N}_{3}\ ...,$\medskip

and for each $k\in \omega $, $\mathcal{M}\preceq _{\Pi _{k+2},\mathrm{%
cons,taller}}\mathcal{N}_{k}.$ \medskip
\end{center}

\noindent In light of part (a)\ of Remark 2.5, this shows that $\mathcal{M}%
\prec _{\mathrm{cons,taller}}\mathcal{N}$, where $\mathcal{N}%
:=\bigcup\limits_{k\in \omega }\mathcal{N}_{k}$. So the proof will be
complete once we explain how to recursively build the desired chain $%
\left\langle \mathcal{N}_{k}:k\in \omega \right\rangle .$ \medskip

\noindent Using the notation of Definition 2.7, let $\mathcal{L}_{\mathrm{set%
}}^{+}=$ $\mathcal{L}_{\mathrm{set}}\cup \{\dot{a}:a\in \mathrm{V}\}$ be the
language $\mathcal{L}_{\mathrm{set}}$ of set theory with a constant $\dot{a}$
for each $a\in \mathrm{V}$. To build $\mathcal{N}_{1}$ we argue within $%
\mathcal{M}$: add a new constant $c$ to $\mathcal{L}_{\mathrm{set}}^{+}$ and
consider the theory $T_{1}$ defined as follows:\medskip

\begin{center}
$T_{1}=\left\{ \mathrm{Ord}(c)\right\} \cup \{\dot{\alpha }\in c:\mathrm{Ord}%
(\alpha )\}\cup \mathrm{Th}_{\Pi _{3}}($\textrm{$V$}$,\in ,\dot{a})_{a\in 
\mathrm{V}}.$\medskip
\end{center}

\noindent Thus $T_{1}$ is a proper class within $\mathcal{M}$. Arguing
within $\mathcal{M}$, and using Fact 2, it is easy to see that $\mathrm{Con}%
(T_{1})$\ holds, and therefore by Fact 1 we can get hold of a proper class
model $\mathcal{N}_{1}$ of $T_{1}$ such that $\mathcal{N}_{1}$ is a
self-extension of $\mathcal{M}$. Thus $\mathcal{M}\prec _{_{\Pi _{3,\mathrm{%
cons}}}}\mathcal{N}_{1}$, and there is an ordinal in $\mathcal{N}_{1}$ that
is above all of the ordinals of $\mathcal{M}$.\medskip

\noindent Next suppose we have built $\left\langle \mathcal{N}_{i}:1\leq
i\leq k\right\rangle $ for some $k\in \omega $ while complying with (1)
through (4). Consider the theory $T_{k}$ defined in $\mathcal{M}$ as
follows:\medskip

\begin{center}
$T_{k+1}=\mathrm{Th}_{\Pi _{k+1}}(\mathcal{N}_{k},b)_{b\in N_{k}}\cup 
\mathrm{Th}_{\Pi _{k+2}}($\textrm{$V$}$,\in ,a)_{a\in \mathrm{V}}.$\medskip
\end{center}

\noindent Fact 2 together with the inductive assumption that $\mathcal{M}%
\prec _{\Pi _{k+1}}\mathcal{N}_{k}$ makes it evident that $\mathrm{Con}%
(T_{k+1})$ holds in $\mathcal{M}$. Hence by Fact 1 we can get hold of a
model $\mathcal{N}_{k+1}$ of $T$ such that $\mathcal{N}_{k+1}$ is a
self-extension of $\mathcal{M}$. This makes it clear that $\left\langle 
\mathcal{N}_{i}:1\leq i\leq k+1\right\rangle $ meets (1) through (4).\hfill $%
\square $\medskip

\noindent \textbf{Second proof of Theorem 3.1}.~Given $\mathcal{M}\models 
\mathsf{ZF}+\exists p\left( \mathrm{V}=\mathrm{HOD}(p)\right) ,$ let $%
\mathbb{B}$ be the Boolean algebra consisting of $\mathcal{M}$-definable
subsets of $\mathrm{Ord}^{\mathcal{M}}$ (ordered by $\subseteq $). In Stage
1 of the proof we will build an appropriate ultrafilter $\mathcal{U}$ on $%
\mathbb{B},$ and in Stage 2 we will verify that the definable ultrapower of $%
\mathcal{M}$ modulo $\mathcal{U}$, denoted $\mathcal{M}^{\mathrm{Ord}^{%
\mathcal{M}}}/\ \mathcal{U}$, is a conservative elementary extension of $%
\mathcal{M}$ that is taller than $\mathcal{M}$. \medskip

\noindent STAGE 1: The following facts come handy in this stage:\medskip

\noindent Fact 1. By Theorem 2.9, for each $n\in \omega $ there is a
definable satisfaction predicate for $\Sigma _{n}$-formulae.\medskip

\noindent Fact 2. By Theorem 2.14, $\mathcal{M}$ has a parametrically
definable global well-ordering.\medskip

\noindent To facilitate the construction of $\mathcal{U}$, let us introduce
some conventions:\medskip

\begin{enumerate}
\item[$(i)$] Given $\mathfrak{X}\subseteq \mathcal{P}(M)$, we say that $%
\mathfrak{X}$ is $\mathcal{M}$-\textit{listable}, if that there is a
parametric formula $\psi (\alpha ,\beta )$ such that for all $X\subseteq M$
we have:\medskip
\end{enumerate}

\begin{center}
$X\in \mathfrak{X}$ iff $\exists \alpha \in \mathrm{Ord}^{\mathcal{M}}$ $%
X=\{m\in M:\mathcal{M}\models \psi (\alpha ,m)\}.$\medskip
\end{center}

\begin{enumerate}
\item[\ ] In the above context, if $\psi $ is a $\Sigma _{n}$-formula, $%
\mathfrak{X}$ is said to be\textit{\ }$\Sigma _{n}$\textit{-listable in }$%
\mathcal{M}$. Moreover, we say that $\left\langle X_{\alpha }:\alpha \in 
\mathrm{Ord}^{\mathcal{M}}\right\rangle $ is an $\mathcal{M}$\textit{-list}
of $\mathfrak{X}$ if $X_{\alpha }=\{m\in M:\mathcal{M}\models \psi (\alpha
,m)\}$ for all $\alpha \in \mathrm{Ord}^{\mathcal{M}}.$

\item[$(ii)$] Given an $\mathcal{M}$-listable $\mathfrak{X}\subseteq \mathbb{%
B}$, we say that $\mathfrak{X}$ is $\mathcal{M}$-\textit{thick} if every $%
\mathcal{M}$-finite intersection of elements of\textit{\ }$\mathfrak{X}$%
\textit{\ }is nonempty. More precisely, $\mathfrak{X}$ is $\mathcal{M}$%
-thick if $\mathcal{M}$ satisfies the following, where $\psi $ is a listing
formula for $\mathfrak{X}$:\medskip
\end{enumerate}

\begin{center}
$\forall k\in \omega $ $\forall \left\langle \alpha _{i}:i<k\right\rangle
\in \mathrm{Ord}^{k}\ \exists \gamma \in \mathrm{Ord}\ \forall i<k\ \psi
(\alpha _{i},\gamma ).$\medskip
\end{center}

\begin{enumerate}
\item[$(iii)$] For each $n\in \omega $, let $\mathbb{B}_{n}$\textrm{\ }be
the collection of $X\in \mathbb{B}$ such that $X$ is $\Sigma _{n}$-definable
in $\mathcal{M}$ (parameters allowed). Note that by Facts 1 and 2, $\mathbb{B%
}_{n}$ is $\mathcal{M}$-listable for each $n\in \omega $.

\item[$(iv)$] $\mathcal{F}$ is the collection of all $X\in \mathbb{B}$ of
the form $\mathrm{Ord}^{\mathcal{M}}\backslash \mathrm{Ext}_{\mathcal{M}%
}(\alpha )$, where $\alpha \in \mathrm{Ord}^{\mathcal{M}}.$ $\mathcal{F}$ is
clearly $\mathcal{M}$-listable and $\mathcal{M}$-thick.\medskip
\end{enumerate}

\noindent We wish to construct an ultrafilter $\mathcal{U}$ $\subseteq $ $%
\mathbb{B}$ that satisfies the following properties:\medskip

\begin{enumerate}
\item[$(1)$] $\mathcal{U}$ is an ultrafilter on $\mathbb{B}$.$\smallskip $

\item[$(2)$] $\mathcal{F\subseteq U}$.$\smallskip $

\item[$(3)$] Given any $\mathcal{M}$-definable $f:\mathrm{Ord}^{\mathcal{M}%
}\rightarrow \mathcal{M}$, there is some $P\in \mathcal{U}$ such that either 
$f\upharpoonright P$ is one-to-one, or the range of $f\upharpoonright P$ is
a set (as opposed to a class).$\smallskip $

\item[$(4)$] Given any $n\in \mathbb{\omega }$, $\mathcal{U}\cap \mathbb{B}%
_{n}$ is $\mathcal{M}$-listable. This condition is equivalent to: for any $%
\mathcal{M}$-list $\left\langle X_{\alpha }:\alpha \in \mathrm{Ord}^{%
\mathcal{M}}\right\rangle $, $\left\{ \alpha \in \mathrm{Ord}^{\mathcal{M}%
}:X_{\alpha }\in \mathcal{U}\right\} \in \mathbb{B}.$\medskip
\end{enumerate}

\noindent The desired ultrafilter $\mathcal{U}$ will be defined as $%
\bigcup\limits_{n\in \omega }\mathcal{U}_{n}$, where $\mathcal{U}_{n+1}$
will be defined from $\mathcal{U}_{n}$ using an internal recursion within $%
\mathcal{M}$ of length $\mathrm{Ord}^{\mathcal{M}}$. Thus, intuitively
speaking, $\mathcal{U}$ will be constructed in $\omega \times \mathrm{Ord}^{%
\mathcal{M}}$-stages. \medskip

\begin{itemize}
\item In what follows, for $n\in \omega $, we fix an $\mathcal{M}$-list $%
\left\langle f_{n,\alpha }:\alpha \in \mathrm{Ord}^{\mathcal{M}%
}\right\rangle $ of all parametrically definable $\Sigma _{n}$-functions%
\footnote{%
The the notion of $\mathcal{M}$-listability can be naturally extended to
family $\mathcal{Y}$ of subsets of $M^{2}$, since $\mathcal{Y}$ can be coded
as a family $\mathfrak{X}_{\mathcal{Y}}$ of subsets of $M$ with the help of
a canonical pairing function.} of $\mathcal{M}$ with domain $\mathrm{Ord}^{%
\mathcal{M}}$, and an $\mathcal{M}$-list $\left\langle S_{n,\alpha }:\alpha
\in \mathrm{Ord}^{\mathcal{M}}\right\rangle $ of all elements of $\mathbb{B}%
_{n}$. Such $\mathcal{M}$-lists exist by Facts 1 and 2 listed at the
beginning of Stage 1.\medskip
\end{itemize}

\noindent Let $\mathcal{U}_{0}$ $:=\mathcal{F}.$ Thus $\mathcal{U}_{0}$ is $%
\mathcal{M}$-listable and $\mathcal{M}$-thick. We now describe the inductive
construction of $\mathcal{U}_{n+1}$ from $\mathcal{U}_{n}$. So suppose $n\in
\omega $, and we have $\mathcal{U}_{n}\subseteq \mathbb{B}$ that satisfies
the clauses $\mathbb{C}_{1}(n)$ through $\mathbb{C}_{3}(n)$ below. Note that
conditions $\mathbb{C}_{2}(n)$ and $\mathbb{C}_{3}(n)$ only kick in for $%
n\geq 1$ and ensure that $\mathcal{U}_{n}$ abides by certain $\Sigma _{n-1}$%
-obligations for $n\geq 1$.\medskip

\begin{enumerate}
\item[$\mathbb{C}_{1}(n):$] $\mathcal{U}_{n}$ is $\mathcal{M}$-listable and $%
\mathcal{M}$-thick.$\smallskip $

\item[$\mathbb{C}_{2}(n):$] If $n\geq 1$, then $\forall \alpha \in \mathrm{%
Ord}^{\mathcal{M}}\ \exists X\in \mathcal{U}_{n}$ ($f_{n-1,\alpha
}\upharpoonright X$ is one-to-one, or the range of $f_{n-1,\alpha
}\upharpoonright X$ is a set).$\smallskip $

\item[$\mathbb{C}_{3}(n):$] If $n\geq 1,$ then$\ \forall \alpha \in \mathrm{%
Ord}^{\mathcal{M}}$ ($S_{n-1,\alpha }\in \mathcal{U}_{n+1}$ or $\mathrm{Ord}%
^{\mathcal{M}}\backslash S_{n-1,\alpha }\in \mathcal{U}_{n}).$\medskip
\end{enumerate}

\noindent The following lemma is crucial for the construction of $\mathcal{U}%
_{n+1}\supseteq \mathcal{U}_{n}$ such that $\mathcal{U}_{n+1}$ satisfies $%
\mathbb{C}_{1}(n+1)$ through $\mathbb{C}_{3}(n+1).$\medskip

\noindent \textbf{Lemma 3.1.1.}\ \textit{Suppose} $n\in \omega $, $\mathcal{F%
}\subseteq \mathfrak{X}\subseteq \mathbb{B},$ $\mathfrak{X}$ \textit{is }$%
\Sigma _{n}$\textit{-listable in }$\mathcal{M}$\textit{, and} $\mathfrak{X}$ 
\textit{is} $\mathcal{M}$\textit{-thick.} \textit{Then there is a large
enough} $k\in \omega $ \textit{that satisfies the following two conditions}%
:\medskip

\begin{enumerate}
\item[\textbf{(a)}] \textit{Given any} $A\subseteq \mathrm{Ord}^{\mathcal{M}%
} $ \textit{that} \textit{is} $\Sigma _{m}$\textit{-definable in }$\mathcal{M%
}$, \textit{let} $\widehat{A}\in \mathbb{B}$ \textit{be defined by}: $%
\widehat{A}=A$ \textit{if} $\mathfrak{X\cup \{}A\}$\textit{\ is }$\mathcal{M}
$\textit{-thick, and otherwise }$\widehat{A}:=\mathrm{Ord}^{\mathcal{M}%
}\backslash A$. \textit{Then} $\mathfrak{X\cup \{}\widehat{A}\}$ \textit{is} 
$\mathcal{M}$\textit{-thick}, \textit{and} $\widehat{A}$ \textit{is} $\Sigma
_{p}$-\textit{definable in }$\mathcal{M}$ \textit{for }$p=\max \{m,n\}+k$.$%
\smallskip $

\item[\textbf{(b)}] \textit{Given any} $f:\mathrm{Ord}^{\mathcal{M}%
}\rightarrow M$ \textit{that} \textit{is} $\Sigma _{m}$\textit{-definable in 
}$\mathcal{M}$\textit{, there is some }$P\in \mathbb{B}$ \textit{that is} $%
\Sigma _{p}$\textit{-definable in }$\mathcal{M}$\textit{\ for }$p=\max
\{m,n\}+k$ \textit{such that }$\mathfrak{X\cup \{}P\}$ \textit{is} $\mathcal{%
M}$\textit{-thick, and }$f\upharpoonright $ $P$ \textit{is either
one-to-one, or the range of }$f\upharpoonright $ $P$ \textit{is coded as an
element of }$\mathcal{M}$. \medskip
\end{enumerate}

\noindent \textbf{Proof.}~(a) is easy to see, so we leave the proof to the
reader. Our proof of (b) will not explicitly specify the bound $k$, instead
we focus on the combinatorics, which is uniform enough to yield the
existence of such a bound. Given an $\mathcal{M}$-definable $f:\mathrm{Ord}^{%
\mathcal{M}}\rightarrow M$, we distinguish two cases:\medskip

\noindent \textbf{Case 1}: There is some $\gamma \in \mathrm{Ord}^{\mathcal{M%
}}$ such that $\mathfrak{X}\cup \{Y_{\gamma }\}$ is $\mathcal{M}$-thick,
where 
\begin{equation*}
Y_{\gamma }=\{\alpha \in \mathrm{Ord}^{\mathcal{M}}:\mathcal{M}\models
f(\alpha )<\gamma \}.
\end{equation*}

\noindent \textbf{Case 2}: Not case 1.\medskip

\noindent If case 1 holds, then we let $P:=Y_{\gamma _{0}}$, where $\gamma
_{0}$ is the first ordinal in $\mathrm{Ord}^{\mathcal{M}}$ witnessing the
veracity of case 1. Clearly the range of\textit{\ }$f\upharpoonright $ $P$
is coded as an element of\textit{\ }$\mathcal{M}$.\medskip

\noindent If case 2 holds, then let $\mathfrak{X}^{\ast }:=$ the set of all $%
\mathcal{M}$-finite intersections of $\mathfrak{X}$. Note that $\mathfrak{X}%
^{\ast }$ is $\mathcal{M}$-listable, so we can let $\left\langle D_{\eta
}:\eta \in \mathrm{Ord}^{\mathcal{M}}\right\rangle $ be an $\mathcal{M}$%
-list of $\mathfrak{X}^{\ast }$. It is also clear that the $\mathcal{M}$%
-thickness of $\mathfrak{X}$ is inherited by $\mathfrak{X}^{\ast }$. We will
construct an $\mathcal{M}$-definable set $P=\{p_{\eta }:\eta \in \mathrm{Ord}%
^{\mathcal{M}}\}$ by recursion (within $\mathcal{M}$) such that the
restriction of $f$ to $P$ is one-to-one, and such that $\mathfrak{X}\cup
\{P\}$ is $\mathcal{M}$-thick$.$ Suppose we have already built $P$ up to
some $\gamma \in \mathrm{Ord}^{\mathcal{M}}$ as $\{p_{\eta }:\eta <\gamma \}$
such that the following two conditions hold in $\mathcal{M}$:\medskip

\noindent $(\ast )$ $\ p_{\eta }<p_{\eta ^{\prime }}$ and $f(p_{\eta
})<f(p_{\eta ^{\prime }})$ whenever $\eta <\eta ^{\prime }<\gamma $.\medskip

\noindent $(\ast \ast )$ $\ p_{\eta }\in $ $D_{\eta }$ for all $\eta <\gamma
.$\medskip

\noindent Within $\mathcal{M}$, let $\theta =\sup \{f(p_{\eta }):\eta
<\gamma \}.$ Note that $\theta $ is well-defined because $\mathcal{M}$
satisfies the replacement scheme. Since we know that case 2 holds, there
must exist some $D_{\delta }\in \mathfrak{X}^{\ast }$ such that:\medskip

\begin{center}
$D_{\delta }\cap Y_{\theta +1}=\varnothing $, where $Y_{\theta +1}:=\{\alpha
\in \mathrm{Ord}^{\mathcal{M}}:\mathcal{M}\models f(\alpha )<\theta +1\}.$%
\medskip
\end{center}

\noindent Therefore within $\mathcal{M}$ can choose $p_{\gamma }$ to be the
first member of $D_{\delta }\cap D_{\gamma \text{ }}$that is above $%
\{s_{\eta }:\eta <\gamma \}$. Since $\mathfrak{X}$ is $\mathcal{M}$-thick
and $\mathcal{F}\subseteq \mathfrak{X}$, $D_{\delta }\cap D_{\gamma \text{ }%
} $is unbounded in $\mathrm{Ord}^{\mathcal{M}},$ hence $s_{\gamma }$ is
well-defined. Thus we have constructed an $\mathcal{M}$-definable set $%
P=\{p_{\eta }:\eta \in \mathrm{Ord}^{\mathcal{M}}\}$ such that $p_{\eta }\in 
$ $D_{\eta }$ for all $\eta \in \mathrm{Ord}^{\mathcal{M}}$ (so $\mathfrak{X}%
\cup \{P\}$ is $\mathcal{M}$-thick), and $f$ is strictly monotone increasing
on $P$. This concludes the proof of part (b) of Lemma 3.1.1.\hfill $\square $%
\medskip

\noindent Using Lemma 3.1.1, it is now straightforward to use transfinite
recursion within $\mathcal{M}$ of length $\mathrm{Ord}^{\mathcal{M}}$ to
construct a family $\mathcal{U}_{n+1}\subseteq \mathbb{B}$ that extends $%
\mathcal{U}_{n}$ and which satisfies clauses $\mathbb{C}_{1}(n+1)$ through $%
\mathbb{C}_{3}(n+1)$. $\mathcal{U}_{n+1}$ will be of the form: \medskip

\begin{center}
$\mathcal{U}_{n}\ \cup \bigcup\limits_{\alpha \in \mathrm{Ord}^{\mathcal{M}%
}}\{A_{\alpha },P_{\alpha }\}$,
\end{center}

\noindent where $A_{\alpha }$ and $P_{\alpha }$ are defined within $\mathcal{%
M}$ by transfinite recursion on $\alpha $. \ Suppose for some $\alpha \in 
\mathrm{Ord}^{\mathcal{M}}$ we have built subsets $A_{\beta }$ and $P_{\beta
}$ of $\mathbb{B}$ for each $\beta <\alpha ,$ and $\mathcal{U}_{n}\cup
\bigcup\limits_{\beta <\alpha }\{A_{\beta },P_{\beta }\}$ is $\mathcal{M}$%
-thick. Now we can use part (a) of Lemma 3.1.1 to construct $A_{\alpha }\in 
\mathbb{B}$ such that $\mathcal{U}_{n}\cup \bigcup\limits_{\beta <\alpha
}\{A_{\beta },P_{\beta }\}\cup \{A_{\alpha }\}$ is $\mathcal{M}$-thick and $%
A_{\alpha }=S_{n,\alpha }$ or $A_{\alpha }=\mathrm{Ord}^{\mathcal{M}%
}\backslash S_{n,\alpha }.$ In the next step, using part (b) of Lemma 3.1.1,
we can construct $P_{\alpha }\in \mathbb{B}$ which `takes care of'\ $%
f_{n,\alpha }$ (i.e., the restriction of $f_{n,\alpha }$ to $P_{\alpha }$ is
either one-to-one, or its range is a set) and which has the property that $%
\mathcal{U}_{n}\cup \bigcup\limits_{\beta <\alpha }\{A_{\beta },P_{\beta
}\}\cup \{A_{\alpha },P_{\alpha }\}$ is $\mathcal{M}$-thick. This concludes
the recursive construction of $A_{\alpha }$ and $P_{\alpha }$ for $\alpha
\in \mathrm{Ord}^{\mathcal{M}}$, thus we can define:\medskip

\begin{center}
$\mathcal{U}_{n+1}:=\mathcal{U}_{n}\ \cup \bigcup\limits_{\alpha \in \mathrm{%
Ord}^{\mathcal{M}}}\{A_{\alpha },P_{\alpha }\}$.\medskip
\end{center}

\noindent The construction makes it clear that $\mathcal{U}_{n+1}$ satisfies 
$\mathbb{C}_{1}(n+1)$ through $\mathbb{C}_{3}(n+1)$. This, in turn, makes
evident that for \medskip

\begin{center}
$\mathcal{U}:=\bigcup\limits_{n\in \omega }\mathcal{U}_{n}$, \medskip
\end{center}

\noindent $\mathcal{U}$\ has properties (1) through (4) promised at the
beginning of Stage 1. Note that $\mathcal{U}$ satisfies condition (4) since $%
\mathcal{U}_{n+1}$ and $\mathbb{B}_{n}$ are both $\mathcal{M}$-listable, and
by clause $\mathbb{C}_{3}(n+1)$, $\mathcal{U}\cap \mathbb{B}_{n}\subseteq 
\mathcal{U}_{n+1}.$ This concludes Stage 1.\medskip

\noindent STAGE 2: Let $\mathcal{N}:=\mathcal{M}^{\mathrm{Ord}^{\mathcal{M}%
}}/\ \mathcal{U}$ be the definable ultrapower of $\mathcal{M}$ modulo $%
\mathcal{U}$ whose universe consists of $\mathcal{U}$-equivalence classes $%
[f]_{\mathcal{U}}$ of $\mathcal{M}$-definable functions $f:\mathrm{Ord}^{%
\mathcal{M}}\rightarrow M$. Thanks to the availability of an $\mathcal{M}$%
-definable global well-ordering, the following \L o\'{s}-style result can be
readily verified. \medskip

\noindent \textbf{Theorem 3.1.2.}~\textit{For all} $k$\textit{-ary formulae} 
$\varphi (x_{0},\cdot \cdot \cdot ,x_{k-1})$ \textit{of} $\mathcal{L}_{%
\mathrm{set}}$ \textit{and all} $k$\textit{-tuples of} $\mathcal{M}$\textit{%
-definable functions} $\left\langle f_{i}:i<k\right\rangle $ \textit{from} $%
\mathrm{Ord}^{\mathcal{M}}$ \textit{to} $M$, \textit{we have}:\medskip

\begin{center}
$\mathcal{N}\models \varphi ([f_{0}]_{\mathcal{U}},\cdot \cdot \cdot
,[f_{k-1}]_{\mathcal{U}})$ iff $\left\{ \alpha \in \mathrm{Ord}^{\mathcal{M}%
}:\mathcal{M}\models \varphi \left( f_{0}(\alpha ),\cdot \cdot \cdot
,f_{k-1}(\alpha )\right) \right\} \in \mathcal{U}$.\smallskip
\end{center}

\noindent For each $m\in M$ let $\widetilde{m}:$ $\mathrm{Ord}^{\mathcal{M}%
}\rightarrow \{m\}$. Theorem 3.1.2 assures us that the map $j:\mathcal{M}%
\rightarrow \mathcal{N}$ given by $j(m)=[\widetilde{m}]_{\mathcal{U}}$ is an
elementary embedding, thus by identifying $\mathcal{M}$\ with the image of $%
j $, we can construe $\mathcal{N}$ as an elementary extension of $\mathcal{M}
$. Let $\mathrm{id}:\mathrm{Ord}^{\mathcal{M}}\rightarrow \mathrm{Ord}^{%
\mathcal{M}}$ be the identity map. By property (2) of $\mathcal{U}$
specified in Stage 1, every element of $\mathcal{U}$ is unbounded in $%
\mathrm{Ord}^{\mathcal{M}}.$ So by Theorem 3.1.2, $\mathcal{N}$ is taller
than $\mathcal{M}$\ since $[\mathrm{id}]_{\mathcal{U}}$ exceeds all ordinals
of $\mathcal{M}$. \medskip

\noindent It remains to verify that $\mathcal{N}$ is a conservative
extension of $\mathcal{M}$. Given a formula $\psi (x,y_{0},\cdot \cdot \cdot
,y_{k-1})$ and parameters $s_{0},\cdot \cdot \cdot ,s_{k-1}$ in $N$, let:
\medskip

\begin{center}
$Y=\left\{ m\in M:\mathcal{N}\models \varphi ([\widetilde{m}]_{\mathcal{U}%
},s_{0},\cdot \cdot \cdot ,s_{k-1})\right\} .$\medskip
\end{center}

\noindent We wish to show that $Y$ is $\mathcal{M}$-definable$\mathfrak{.}$
Choose $f_{0},\cdot \cdot \cdot ,f_{k-1}$ so that $s_{i}=[f_{i}]_{\mathcal{U}%
}$ for all $i<p$ and consider the formula $\theta (x,\alpha )$ defined
below:\medskip

\begin{center}
$\theta (x,\alpha ):=\alpha \in \mathrm{Ord}\wedge \varphi (x,f_{0}(\alpha
),\cdot \cdot \cdot ,f_{k-1}(\alpha ))$.\medskip
\end{center}

\noindent Pick $n\in \omega $ large enough so that $\theta (x,\alpha )$ is a 
$\Sigma _{n}$-formula. For each $m\in M,$ let:\medskip\ 

\begin{center}
$S_{m}:=\{\alpha \in \mathrm{Ord}^{\mathcal{M}}:\mathcal{M}\models \theta
(m,\alpha )\}$. \medskip
\end{center}

\noindent Note that the choice of $n$, $S_{m}\in \mathbb{B}_{n}$ for each $%
m\in M$, so $S_{m}\in \mathcal{U}$ iff $S_{m}\in \mathcal{U}\cap \mathbb{B}%
_{n}$ for each $m\in M$. On the other hand, by Theorem 3.1.2 we have:\medskip

\begin{center}
$m\in Y$ iff $S_{m}\in \mathcal{U}$ for each $m\in M$. \medskip
\end{center}

\noindent So we can conclude that $m\in Y$ iff $S_{m}\in \mathcal{U}\cap 
\mathbb{B}_{n}$ for each $m\in M.$ Since $\mathcal{U}\cap \mathbb{B}_{n}$ is 
$\mathcal{M}$-listable by property (4) of $\mathcal{U}$ specified in Stage
1, we can write $\mathcal{U}\cap \mathbb{B}_{n}$ as an $\mathcal{M}$-list $%
\left\langle X_{\eta }:\eta \in \mathrm{Ord}^{\mathcal{M}}\right\rangle $
given by some formula $\psi (x,y).$ Thus the following holds for all $m\in M$%
, which shows that $Y$ is $\mathcal{M}$-definable.\medskip

\begin{center}
$m\in Y$ iff $\exists \eta \left( S_{m}=X_{\eta }\right) $ iff $\mathcal{M}%
\models \exists \eta \ \forall \alpha \left( \theta (m,\alpha
)\leftrightarrow \psi (\eta ,m)\right) .$
\end{center}

\hfill $\square $\medskip

\noindent \textbf{3.2.~Corollary.}~\textit{Work in }$\mathsf{ZF}+\exists
p\left( \mathrm{V}=\mathrm{HOD}(p)\right) $ \textit{in the metatheory}, 
\textit{and let} $\mathcal{M}$ \textit{be an arbitrary model of} $\mathsf{ZF}%
+\exists p\left( \mathrm{V}=\mathrm{HOD}(p)\right) .$\medskip

\begin{enumerate}
\item[\textbf{(a)}] \textit{For every regular cardinal }$\kappa $\textit{\
there is a conservative elementary extension} $\mathcal{N}_{\kappa }$ 
\textit{of} $\mathcal{M}$ \textit{that is taller than }$\mathcal{M}$ \textit{%
and} $\mathrm{cf}(\mathrm{Ord}^{\mathcal{N}})=\kappa $.\footnote{%
As indicated by the proof of this part, $\mathsf{ZFC}$ suffices as the
metatheory for this part.}\medskip

\item[\textbf{(b)}] \textit{There is a proper class model }$\mathcal{N}_{%
\mathrm{Ord}}$ \textit{whose elementary diagram is a definable class such
that }$(1)$\textit{\ }$\mathcal{M\prec }_{\mathrm{cons,taller}}\mathcal{N}_{%
\mathrm{Ord}}$, \textit{and }$(2)$\textit{\ }$\mathrm{Ord}^{\mathcal{N}_{%
\mathrm{Ord}}}$ \textit{contains a cofinal class of order-type }$\mathrm{Ord}
$.\medskip
\end{enumerate}

\noindent \textbf{Proof.}~To establish (a), we observe that given a regular
cardinal\textit{\ }$\kappa $, Theorem 3.1 allows us to build an elementary
chain of models $\left\langle \mathcal{M}_{\alpha }:\alpha \in \kappa
\right\rangle $ with the following properties:\medskip

\noindent (1) $\ \ \mathcal{M}_{0}=\mathcal{M}.$\medskip

\noindent (2) \ \ For each $\alpha \in \kappa $ $\mathcal{M}_{\alpha +1}$ is
a conservative elementary extension of $\mathcal{M}$ that is taller than $%
\mathcal{M}.$\medskip

\noindent (3) \ \ For limit $\alpha \in \kappa $, $\mathcal{M}_{\alpha
}=\bigcup\limits_{\beta \in \alpha }\mathcal{M}_{\beta }.$

\noindent The desired model $\mathcal{N}$ is clearly $\mathcal{N}_{\kappa
}=\bigcup\limits_{\alpha \in \kappa }\mathcal{M}_{\alpha }.$ This concludes
the proof of part (a).

\noindent With part (a) at hand, it is evident that in the presence of a
global definable well-ordering (available thanks to the veracity of $\exists
p\left( \mathrm{V}=\mathrm{HOD}(p)\right) $ in the metatheory), we can
construct the desired model $\mathcal{N}_{\mathrm{Ord}}$ as the union of an
elementary chain of models $\left\langle \mathcal{M}_{\alpha }:\alpha \in 
\mathrm{Ord}\right\rangle $ that satisfies (1) and it also satisfies the
result of replacing $\kappa $ with $\mathrm{Ord}$ in (2) and (3).\hfill $%
\square $\medskip

\noindent \textbf{3.3.~Remark.}~We should point out that condition (3) of
the proof of Theorem 3.1 was inserted so as to obtain the following
corollary which is of interest in light of a result of \cite{EnayatMIN}~that
states that every model of $\mathsf{ZFC}$ has a cofinal conservative
extension that possesses a minimal \textit{cofinal} elementary extension%
\footnote{$\mathcal{N}$ is a \textit{minimal} elementary extension of $%
\mathcal{M}$ if $\mathcal{M}\prec \mathcal{N}$ and there is no $\mathcal{K}$
such that $\mathcal{M}\prec \mathcal{K}\prec \mathcal{N}$.}.\medskip

\noindent \textbf{3.4.~Theorem.}~\textit{Every model of }$\mathsf{ZF}%
+\exists p\left( \mathrm{V}=\mathrm{HOD}(p)\right) $\textit{\ has a cofinal
conservative elementary extension that possesses a minimal elementary end
extension.}\medskip

\noindent \textbf{Proof}.~Let $\mathcal{M}$ be a model of $\mathsf{ZF}%
+\exists p\left( \mathrm{V}=\mathrm{HOD}(p)\right) $\textit{, }$\mathcal{N}$%
\textit{\ }be the model constructed in the second proof of Theorem 3.1, and
let\textit{\ }$\mathcal{M}^{\ast }$ be the convex hull of $\mathcal{M}$ in $%
\mathcal{N}$. Note that by Gaifman Splitting Theorem we have:\smallskip

\begin{center}
$\mathcal{M}\preceq _{\mathrm{cof}}\mathcal{M}^{\ast }\prec _{\mathrm{end}}%
\mathcal{N}$.\smallskip
\end{center}

\noindent The above together with the fact that $\mathcal{N}$ is
conservative extension of $\mathcal{M}$ makes it clear that $\mathcal{M}%
^{\ast }$ is a conservative extension of $\mathcal{M}$.\footnote{%
By Theorem 4.1, $\mathcal{M}^{\ast }$ is a proper elementary extension of $%
\mathcal{M}$.} We claim that $\mathcal{N}$ is a minimal elementary end
extension of $\mathcal{M}^{\ast }$. First note that $\mathcal{N}$ is
generated by members of $\mathcal{M}$ and the $\mathcal{U}$-equivalence
class of the identity function $[\mathrm{id}]_{\mathcal{U}}$ via the
definable terms of $\mathcal{M}$, i.e., the following holds:\medskip

\begin{center}
For any $a\in N$, there exists some $\mathcal{M}$-definable function $f$
such that $\mathcal{N}\models f([\mathrm{id}]_{\mathcal{U}})=a.$\medskip
\end{center}

\noindent To verify that $\mathcal{N}$ is a minimal elementary end extension
of $\mathcal{M}^{\ast }$ it is sufficient to show that if $s\in N\backslash
M^{\ast }$ then $[\mathrm{id}]_{\mathcal{U}}$ is definable in $\mathcal{N}$
using parameters from $M\cup \{s\}.$ Choose $h$ such that $s=[h]_{\mathcal{U}%
}$ and recall that by our construction of $\mathcal{U}$ we know that there
exists $X\subseteq \mathrm{Ord}^{\mathcal{M}}$ with $X\in \mathcal{U}$ such
that either $h\upharpoonright X$ is one-to-one, or the image of $X$ under $h$
is a set (and thus of bounded $\mathcal{M}$-rank). However, since $s$ was
chosen to be in $N\backslash M^{\ast }$ we can rule out the latter
possibility. Note that since $X\in \mathcal{U}$ the extension of $X$ in $%
\mathcal{N}$ will contain $[\mathrm{id}]_{\mathcal{U}}$, i.e., if $X$ is
defined in $\mathcal{M}$ by $\varphi (x)$, then $\mathcal{N}\models \varphi
([\mathrm{id}]_{\mathcal{U}})$. Moreover, by elementarity, since $h$ is
one-to-one on $X$ in $\mathcal{M}$, it must remain so in $\mathcal{N}$.
Hence, in light of the fact that $\mathcal{N}\models h([\mathrm{id}]_{%
\mathcal{U}})=s$, we can define $[\mathrm{id}]_{\mathcal{U}}$ within $%
\mathcal{N}$ to be the unique element of $h^{-1}(s)\cap X.$ This completes
the proof.\hfill $\square $\medskip

\noindent \textbf{3.5.~Remark.}~As Theorem 4.1 of the next section shows, in
Theorem 3.1 the model $\mathcal{N}$ cannot be required to fix $\omega ^{%
\mathcal{M}}$ and in particular $\mathcal{N}$ cannot be arranged to end
extend $\mathcal{M}$. By putting this observation together with Remark
2.6(j) we can conclude that in the proof of Corollary 3.2, $\omega ^{%
\mathcal{N}_{\kappa }}$ is $\kappa $-like (i.e., its cardinality is $\kappa $
but each proper initial segment of it is of cardinality less than $\kappa $%
), and $\omega ^{\mathcal{N}_{\mathrm{Ord}}}$ is $\mathrm{Ord}$-like (i.e.,
it is a proper class every proper initial segment of which forms a set).
\medskip

We next address the question of the extent to which the hypothesis in
Theorem 3.1 that $\mathcal{M}$\ is a model of $\exists p\left( \mathrm{V}=%
\mathrm{HOD}(p)\right) $ can be weakened. In what follows $\mathrm{AC}$
stands for the axiom of choice, and $\mathsf{DC}$ stands for the axiom of
dependent choice. \medskip

\noindent \textbf{3.6.~Proposition.}~\textit{Suppose} $\mathcal{M}$ \textit{%
and} $\mathcal{N}$ \textit{are models of }$\mathsf{ZF}$, $\mathcal{M}\prec _{%
\mathrm{cons}}\mathcal{N}$, \textit{and} $\mathcal{N}$ \textit{is taller than%
} $\mathcal{M}$\textit{.} \textit{If} \textit{every set can be linearly
ordered in} $\mathcal{N}$, \textit{then there is an }$\mathcal{M}$-\textit{%
definable global linear ordering. }\medskip

\noindent \textbf{Proof}.~Suppose $\mathcal{M}$ and $\mathcal{N}$ are as in
the assumptions of the proposition. Fix $\alpha \in \mathrm{Ord}^{\mathcal{N}%
}$ such that $\alpha $ is above all the ordinals of $\mathcal{M}$ and let $R$
denote a linear ordering of $\mathrm{V}_{\alpha }^{\mathcal{N}}$ in $%
\mathcal{N}$. By conservativity, there is a formula $r(x,y)$ (possibly with
parameters from $M$) such that $r(x,y)$ that defines $R\cap M^{2}$ in $%
\mathcal{M}$ . It is clear that $r(x,y)$ defines a global linear ordering in 
$\mathcal{M}$.\hfill $\square $\medskip

\noindent \textbf{3.7.~Corollary.}~\textit{There are models of} $\mathsf{ZFC}
$ \textit{that have no taller conservative elementary extensions to a model
of }$\mathsf{ZFC}$. \medskip

\noindent \textbf{Proof.} This follows from Proposition 3.6, and the
well-known fact that there are models of $\mathsf{ZFC}$ in which the
universe is not definably linearly ordered.\footnote{%
Easton proved (in his unpublished dissertation \cite{Easton Thesis}) that
assuming $\mathrm{Con}(\mathsf{ZF})$ there is a model $\mathcal{M}$ of $%
\mathsf{ZFC}$ that carries no $\mathcal{M}$-definable global choice function
for the class of pairs in $\mathcal{M}$; and in particular the universe
cannot be definably linearly ordered in $\mathcal{M}$. Easton's theorem was
exposited by Felgner \cite[p.231]{Felgner}; for a more recent and
streamlined account, see Hamkins' MathOverflow answer \cite{Joel-failure of
class choice for pairs}.}\hfill $\square $\medskip

\noindent \textbf{3.8.~Question}.~\textit{Can the conclusion of Proposition
3.6} \textit{be improved to }$\mathcal{M}\models \exists p\left( \mathrm{V}=%
\mathrm{HOD}(p)\right) $?\medskip

\noindent \textbf{3.9. Remark}.~If $\mathcal{M}\models \mathsf{ZF}$ and $%
\mathcal{M}\prec _{\Delta _{0},\ \mathrm{cons}}\mathcal{N}$, then each $k\in
\omega ^{\mathcal{M}}$ is fixed. To see this, suppose to the contrary that $%
r\in \mathrm{Ext}_{\mathcal{N}}(k)\backslash \mathrm{Ext}_{\mathcal{M}}(k)$.
Consider the set $S=\left\{ x\in \mathrm{Ext}_{\mathcal{M}}(k):\mathcal{N}%
\models x<r\right\} .$ It is easy to see that $S$ has no last element. On
the other hand, by the assumption $\mathcal{M}\prec _{\mathrm{cons}}\mathcal{%
N}$, the set $S=$ $\{x\in \mathrm{Ext}_{\mathcal{M}}(k):\mathcal{N}\models
x<r\}$ is a definable subset of the predecessors of $r$ and thus by the
veracity of the separation scheme in $\mathcal{M}$, $S$ is coded in $%
\mathcal{M}$ by some element, which makes it clear that $S$ must have a last
element, contradiction. This shows that in Theorem 4.1, $\mathcal{N}$ fixes
each $k\in \omega ^{M}.$ A natural question is whether a taller conservative
elementary extension of $\mathcal{M}$ can be arranged to fix $\omega ^{%
\mathcal{M}}$. As we shall see in the next section, this question has a
negative answer for models of $\mathsf{ZFC}$ (see Corollary 4.2).\medskip

\noindent \textbf{3.10. Remark}.~Theorem 4.1 should also be contrasted with
Gaifman's theorem mentioned in Remark 2.5 which bars taller elementary
self-extensions of models of (a fragment of) $\mathsf{ZF}$. Thus in
Gaifman's theorem, (a) the condition \textquotedblleft $\mathcal{N}$ is an
elementary self extension of $\mathcal{M}$\textquotedblright\ cannot be
weakened to \textquotedblleft $\mathcal{N}$ is an elementary conservative
extension of $\mathcal{M}$\textquotedblright . Also note that in the first
proof of Theorem 3.1, each $\mathcal{N}_{k}$ is a self-extension of $%
\mathcal{M}$ whose associated embedding $j$ is a $\Sigma _{k+2}$-elementary
embedding whose image is \textit{not} cofinal in $\mathcal{M}$.

\textbf{\bigskip }

\begin{center}
\textbf{4}.~\textbf{FAITHFUL EXTENSIONS}\bigskip
\end{center}

This section contains the proof of Theorem B of the abstract (see Corollary
4.2). Recall from Definition 2.4 and Remark 2.5 that the notion of a
faithful extension is a generalization of the notion of a conservative
extension. The notion of a full satisfaction class was defined in Definition
2.7(b). \medskip

\noindent \textbf{4.1.~Theorem.}~\textit{Suppose }$\mathcal{M}$\textit{\ and 
}$\mathcal{N}$\textit{\ are both models of }$\mathsf{ZFC}$\textit{\ such
that }$\mathcal{M}\prec _{\Delta _{0}^{\mathcal{P}},\mathrm{faith}}\mathcal{N%
}$, $\mathcal{N}$ \textit{fixes} $\omega ^{\mathcal{M}},$ \textit{and }$%
\mathcal{N}$\textit{\ is taller than }$\mathcal{M}$. \textit{Then:}\medskip

\begin{enumerate}
\item[\textbf{(a)}] \textit{There is some }$\gamma \in \mathrm{Ord}^{%
\mathcal{N}}$ \textit{such that} $\mathcal{M}\preceq \mathcal{N}_{\gamma }.$ 
\textit{Thus} \textit{either} $\mathcal{M=N}_{\gamma }$, \textit{in which
case} $\mathcal{N}$ \textit{is a topped rank extension of} $\mathcal{M}$, 
\textit{or }$\mathcal{M}\prec \mathcal{N}_{\gamma }$.\footnote{%
As pointed out in Remark 4.5, both cases of the dichotomy can be realized.}%
\medskip

\item[\textbf{(b)}] \textit{In particular, there is a full satisfaction
class }$S$\textit{\ for }$\mathcal{M}$\textit{\ such that }$S$\textit{\ can
be written as }$D\cap M$, \textit{where} $D$ \textit{is }$\mathcal{N}$%
\textit{-definable.}\medskip
\end{enumerate}

\noindent \textbf{Proof}\footnote{%
The proof strategy of Theorem 4.1 is a variant of the proofs of the
following results: \cite[Theorem 1.5]{Ali-TAMS}, \cite[Theorem 3.3]%
{AliJSL1987}, \cite[Theorem 2.1]{Ali+Joel}, \cite[Theorem 2.1.3]%
{Ali-Aut-and-Mahlo}, and \cite[Lemma 2.19]{Ali-ZFI}. A variant of the same
strategy is used in the proof of Theorem 4.5. The origins of the strategy
can be traced to Kaufmann's refinement \cite[Lemma 1.4]{Kaufmann} of a
Skolem hull argument due to Keisler and Silver \cite[Theorem 2.1]%
{Keisler-Silver}.}\textbf{.}~Fix some $\mathcal{N}$-ordinal $\lambda \in 
\mathrm{Ord}^{\mathcal{N}}$ that dominates each $\mathcal{M}$-ordinal, and
some $\mathcal{N}$-ordinal $\beta >\lambda $. Note that: \medskip

\noindent (1) $\mathcal{N}_{\beta }\prec _{\Delta _{0}^{\mathcal{P}}}%
\mathcal{N}.$\medskip

\noindent Thanks to the availability of $\mathrm{AC}$ in $\mathcal{N}$, we
can choose an ordering $\vartriangleleft $ in $\mathcal{N}$ such that $%
\mathcal{N}$ satisfies \textquotedblleft $\vartriangleleft $ \ is a
well-ordering of $\mathrm{V}_{\beta }$\textquotedblright . Then for each $%
m\in M$, we define the following set $K_{m}$ (again within $\mathcal{N}$)
as:\medskip

\begin{center}
$K_{m}:=\{a\in \mathrm{V}:\mathcal{N}\models a\in \mathrm{Def}(\mathrm{V}%
_{\beta },\in ,\vartriangleleft ,\lambda ,m)\},$\medskip
\end{center}

\noindent where $x\in \mathrm{Def}(\mathrm{V}_{\beta },\in ,\vartriangleleft
,\lambda ,m)$ is shorthand for the formula of set theory that
expresses:\medskip

\begin{center}
$x$ is definable in the structure $(\mathrm{V}_{\beta },\in
,\vartriangleleft ,\lambda ,m).$\medskip
\end{center}

\noindent Note that the only allowable parameters used in a definition of $x$
are $\lambda $ and $m$. Thus, intuitively speaking, within $\mathcal{N}$ the
set $K_{m}$ consists of elements $a$ of $\mathrm{V}_{\beta }$ such that $a$
is first order definable in $\left( \mathrm{V}_{\beta },\in
,\vartriangleleft ,\lambda ,m\right) $. Also note that since we are not
assuming that $\mathcal{M}$ is $\omega $-standard, an element $a\in K_{m}$
need not be definable in the structure $\left( \mathcal{N}_{\beta
},\vartriangleleft ,\lambda ,m\right) $ in the real world. Next we move
outside of $\mathcal{N}$ and define $K$ as follows:\medskip

\begin{center}
$K:=\bigcup\limits_{m\in M}K_{m}.$\medskip
\end{center}

\noindent We observe that for each $m\in M$, $K_{m}$ is coded in $\mathcal{N}
$, but there is no reason to expect that $K$ is coded in $\mathcal{N}$. Thus
every element of $K$ is definable, in the sense of $\mathcal{N}$, in $(%
\mathrm{V}_{\beta },\in ,\vartriangleleft )$ with some appropriate choice of
parameters in $\{\lambda \}\cup M.$ Let $\mathcal{K}$ be the submodel of $%
\mathcal{N}$ whose universe is $K$. Using the Tarski test for elementarity
and the fact that, as viewed from $\mathcal{N}$, $\vartriangleleft $
well-orders $\mathrm{V}_{\beta }\mathcal{\ }$we have:\medskip

\noindent (2) $\ \ \mathcal{M}\subsetneq \mathcal{K}\preceq \mathcal{N}%
_{\beta }$. \medskip

\noindent By putting (1)\ and (2) together we conclude:\medskip

\noindent (3) $\ \ \mathcal{M}\prec _{\Delta _{0}^{\mathcal{P}}}\mathcal{K}%
\preceq \mathcal{N}_{\beta }\prec _{\Delta _{0}^{\mathcal{P}}}\mathcal{N}$%
.\medskip

\noindent Let $O^{\ast }$ be the collection of `ordinals'\ of $\mathcal{K}$
that are above the `ordinals'\ of $\mathcal{M}$, i.e., \medskip

\begin{center}
$O^{\ast }=\left\{ \gamma \in \mathrm{Ord}^{\mathcal{K}}:\forall \alpha \in 
\mathrm{Ord}^{\mathcal{M}}\ \mathcal{N}\models \alpha \in \gamma \right\} .$%
\medskip
\end{center}

\noindent Clearly $O^{\ast }$ is nonempty since $\lambda \in O$ $.$ We now
consider the following two cases. As we shall see, Case I leads to the
conclusion of the theorem, and Case II\ is impossible.\medskip

\noindent Case I. $O^{\ast }$ has a least ordinal (under $\in ^{\mathcal{K}%
}).$\medskip

\noindent Case II.\ $O^{\ast }$ has no least ordinal.\medskip

\noindent Suppose Case\ I holds and let $\gamma _{0}=\min (O^{\ast })$.
Clearly $\gamma _{0}$ is a limit ordinal of $\mathcal{N}$. \ By the choice
of $\gamma _{0},$ $\mathrm{Ord}^{\mathcal{M}}$ is cofinal in $\mathrm{Ord}^{%
\mathcal{K}_{\gamma _{0}}}.$ Since by (3) and part (d) of Remark 2.6, $%
\mathrm{V}_{\alpha }^{\mathcal{M}}=\mathrm{V}_{\alpha }^{\mathcal{K}}$ for
each $\alpha \in \mathrm{Ord}^{\mathcal{M}}$, this makes it clear that:
\medskip

\noindent (4) $\ \ \mathcal{M}\preceq _{\Delta _{0}^{\mathcal{P}},\ \mathrm{%
cof}}\mathcal{K}_{\gamma _{0}}.$\medskip

\noindent Since $\mathcal{K}_{\gamma _{0}}$ is the convex hull of $\mathcal{M%
}$ in $\mathcal{K}$, by Gaifman Splitting Theorem, (4) shows:\medskip

\noindent (5) $\ \ \mathcal{M}\preceq \mathcal{K}_{\gamma _{0}}$. \medskip

\noindent On the other hand since by (3) $\mathcal{K}\preceq \mathcal{N}%
_{\beta }$, we have:\medskip

\noindent (6) $\mathcal{K}_{\gamma _{0}}=\mathrm{V}_{\gamma _{0}}^{\mathcal{K%
}}\preceq \mathrm{V}_{\gamma _{0}}^{\mathcal{N}_{\beta }}=\mathrm{V}_{\gamma
_{0}}^{\mathcal{N}}=\mathcal{N}_{\gamma _{0}}.$\medskip

\noindent By (5) and (6), $\mathcal{M}\preceq \mathcal{N}_{\gamma _{0}}$.
Thus the proof of the theorem will be complete once we show that Case II\
leads to a contradiction. Within $\mathcal{N}$ let $s_{\beta }$ be the
elementary diagram of $\left( \mathrm{V}_{\beta },\in \right) $, i.e., $%
\mathcal{N}$ views $s_{\beta }$ is the Tarskian satisfaction class for $%
\left( \mathrm{V}_{\beta },\in \right) $, and let $\Phi
:=\bigcup\limits_{m\in M}\Phi _{m}$, where:\medskip

\begin{center}
$\Phi _{m}:=\{x\in M:\mathcal{N}\models $ $x$ is (the code of) a formula $%
\varphi (c,\dot{m})$ such that $\varphi (\dot{\lambda},\dot{m})\in s_{\beta
}\}\mathrm{.}$\medskip
\end{center}

\noindent So intuitively speaking, $\Phi $ is the result of replacing $\dot{%
\lambda}$ by $c$ (where $c$ is a fresh constant) in the sentences in the
elementary diagram of $\mathcal{N}_{\beta }$ (as computed in $\mathcal{N}$)
whose constants are in $\{\dot{\lambda}\}\cup \{\dot{m}:m\in M\}.$ Since $%
\mathcal{N}$ need not be $\omega $-standard, the elements of $\Phi $ might
be nonstandard formulae.\medskip

\noindent By the assumption that $\mathcal{N}$ is a faithful extension of $%
\mathcal{M}$, $\Phi $ is $\mathcal{M}$-amenable. Next let:\medskip

\begin{center}
$\Gamma :=\left\{ 
\begin{array}{c}
t(c,\dot{m})\in M:t(c,\dot{m})\in \Phi ,\ \mathrm{and}\ \forall \theta \in 
\mathrm{Ord\ }\left( t(c,\dot{m})>\dot{\theta}\right) \in \Phi ,\ \mathrm{and%
}\  \\ 
t\ \mathrm{is\ a\ definable\ term\ in\ the\ language}\ \mathcal{L}_{\mathrm{%
set}}\cup \{c\}\cup \{\dot{m}:m\in M\}%
\end{array}%
\right\} .$ $\ $\medskip
\end{center}

\noindent Officially speaking, $\Gamma $ consists of \textit{syntactic
objects} $\varphi (c,\dot{m},x)$ in $\mathcal{M}$ that satisfy the following
three conditions in $(\mathcal{M},\Phi )$:\medskip

\begin{enumerate}
\item[$(i)$] $\left[ \exists !x\varphi (c,\dot{m},x)\right] \in \Phi .$%
\medskip

\item[$(ii)$] $\left[ \forall x\left( \varphi (c,\dot{m},x)\rightarrow x\in 
\mathrm{Ord}\right) \right] \in \Phi .$\medskip

\item[$(iii)$] $\forall \theta \in \mathrm{Ord}\ \left[ \forall x\left(
\varphi (c,\dot{m},x)\rightarrow \dot{\theta}\in x\right) \right] \in \Phi .$%
\medskip
\end{enumerate}

\noindent Note that $\Gamma $ is definable in $\left( \mathcal{M},\Phi
\right) $. Since we are considering Case II, $\left( \mathcal{M},\Phi
\right) \models \psi $, where:\medskip

\begin{center}
$\psi :=\forall t\left[ t\in \Gamma \longrightarrow \exists t^{\prime }\in
\Gamma \ (t^{\prime }\in t)\in \Phi \right] .$\medskip
\end{center}

\noindent The veracity of the dependent choice scheme in $\mathcal{M}$ (see
part (f) of Definition 2.1), together with the facts that $\Phi $ is $%
\mathcal{M}$-amenable and $\Gamma $ is definable in $\left( \mathcal{M},\Phi
\right) $\textrm{\ }make it clear that there is a sequence $s=\left\langle
t_{n}:n\in \omega ^{\mathcal{M}}\right\rangle $ in $\mathcal{M}$ such
that:\medskip

\begin{center}
$\left( \mathcal{M},\Phi \right) \models \forall n\in \omega \ \left[
t_{n}\in \Gamma \wedge \left( \left( t_{n+1}\in t_{n}\right) \in \Phi
\right) \right] .$\medskip
\end{center}

\noindent Since $s$\textbf{\ }is a countable object in $\mathcal{M}$, and $%
\mathcal{N}$ fixes $\omega ^{M}$ by assumption, $s$ is fixed in the passage%
\textbf{\ }from\textbf{\ }$\mathcal{M}$ to $\mathcal{N}$ by Remark\textbf{\ }%
2.6(c). On the other hand, since $\mathcal{N}$ has a satisfaction predicate $%
s_{\beta }$ for $\mathcal{N}_{\beta }$, this leads to a contradiction
because we have: \medskip

\begin{center}
$\mathcal{N}\models \left\langle t_{n}^{(\mathrm{V}_{\beta },\in )}(\dot{%
\lambda}):n\in \omega \right\rangle $ is an infinite decreasing sequence of
ordinals. \medskip
\end{center}

\noindent In the above $t_{n}(\dot{\lambda})$ is the term obtained by
replacing $c$ with $\dot{\lambda}$ in $t_{n},$ and $t_{n}^{(\mathrm{V}%
_{\beta },\in )}(\dot{\lambda})$ is the interpretation of $t_{n}(\dot{\lambda%
})$ in $(\mathrm{V},\in )$ using the Tarskian satisfaction class for $(%
\mathrm{V}_{\beta },\in )$.\hfill $\square $\medskip

\noindent \textbf{4.2.~Corollary.}~\textit{If }$\mathcal{M}\models \mathsf{%
ZFC}$, $\mathcal{M}\prec _{\Delta _{0}^{\mathcal{P}},\ \mathrm{cons}}%
\mathcal{N}$, \textit{and} $\mathcal{M}$ \textit{fixes} $\omega ^{\mathcal{M}%
},$ \textit{then} $\mathcal{M}$ \textit{is cofinal in} $\mathcal{N}$. 
\textit{Thus a conservative elementary extension of }$\mathcal{M}$ \textit{%
that fixes} $\omega ^{\mathcal{M}}$ \textit{is a cofinal extension}. \medskip

\noindent \textbf{Proof.}~If $\mathcal{M}\models \mathsf{ZFC}$, $\mathcal{M}%
\prec _{\mathrm{cons}}\mathcal{N}$, and $\mathcal{M}$ is not cofinal in $%
\mathcal{N}$, then $\mathcal{N}$ is taller than $\mathcal{M}$. But then by
part (b)\ of Theorem 4.1, there is an $\mathcal{M}$-definable full
satisfaction class over $\mathcal{M}$, which contradicts Tarski's
Undefinability of Truth Theorem.\hfill $\square $\medskip

In Theorem 4.4 below we show that the conclusion of Theorem 4.1 holds if $%
\mathsf{ZFC}$ is weakened to $\mathsf{ZF}\mathrm{,}$ but the assumption that 
$\mathcal{N}$\ fixes $\omega ^{\mathcal{M}}$ is strengthened to $\mathcal{M}%
\subsetneq _{\mathrm{end}}\mathcal{N}$. We first establish a lemma that will
come handy in the proof of Theorem 4.4.\medskip

\noindent \textbf{4.3.~Lemma}~\textit{Suppose} $\mathcal{M}$ \textit{and} $%
\mathcal{N}$ \textit{are models of} $\mathsf{ZF}$, \textit{and }$\mathcal{M}%
\subsetneq _{\mathrm{end,faithful}}\mathcal{N}$. \textit{Then} $\mathcal{M}%
\subsetneq _{\mathrm{rank}}\mathcal{N}$.\medskip

\noindent \textbf{Proof.}~Suppose $\mathcal{M}$ and $\mathcal{N}$ are as in
the assumptions of the lemma. Since $\mathcal{N}$ end extends $\mathcal{M}$,
in order to show that $\mathcal{N}$ rank extends $\mathcal{M}$ it is
sufficient that if $\alpha \in \mathrm{Ord}^{\mathcal{M}}$ and $m\in M$ with 
$\mathcal{M}\models m=\mathrm{V}_{\alpha }$, then $\mathcal{N}\models m=%
\mathrm{V}_{\alpha }.$ For this purpose, it is sufficient to verify that $%
\mathcal{N}$ is a powerset-preserving end extension of $\mathcal{M}$, since
the formula expressing $x=\mathrm{V}_{\alpha }$ is $\Sigma _{1}^{\mathcal{P}%
} $ and is therefore absolute for powerset-preserving end extensions (as
noted in Remark 2.6(e)). Towards this goal, let $a\in M$ and suppose $b\in N$
such that:\medskip

\noindent (1) $\ \ \mathcal{N}\models b\subseteq a.$ \medskip

\noindent Let $X_{b}=\mathrm{Ext}_{\mathcal{N}}(b)\cap \mathrm{Ext}_{%
\mathcal{M}}(a)=\mathrm{Ext}_{\mathcal{N}}(b).$ Thanks to the assumption
that $\mathcal{N}$ is a faithful extension of $\mathcal{M}$, the expansion $(%
\mathcal{M},X_{b})$ satisfies Separation in the extended language and
therefore there is an element $b^{\prime }\in M$ that codes $X_{b}.$ Thus
\medskip

\noindent (2) $\ \ \mathrm{Ext}_{\mathcal{M}}(b^{\prime })=\mathrm{Ext}_{%
\mathcal{N}}(b)\cap \mathrm{Ext}_{\mathcal{M}}(a),$ and \medskip

\noindent Since $\mathcal{N}$ end extends $\mathcal{M}$, we have:\medskip

\noindent (3) \ \ $\mathrm{Ext}_{\mathcal{M}}(a)=\mathrm{Ext}_{\mathcal{N}%
}(a)$ and $\mathrm{Ext}_{\mathcal{M}}(b^{\prime })=\mathrm{Ext}_{\mathcal{N}%
}(b^{\prime }).$\medskip

\noindent So by (1) and (2),\medskip

\noindent (4) $\mathrm{Ext}_{\mathcal{N}}(b^{\prime })=\mathrm{Ext}_{%
\mathcal{N}}(b)\cap \mathrm{Ext}_{\mathcal{N}}(a).$\medskip

\noindent By (1), $\mathrm{Ext}_{\mathcal{N}}(b)\subseteq \mathrm{Ext}_{%
\mathcal{N}}(a),$ so combined with (4) this yields:\medskip

\noindent (5) $\ \ \mathrm{Ext}_{\mathcal{N}}(b^{\prime })=\mathrm{Ext}_{%
\mathcal{N}}(b)$, \medskip

\noindent which by Extensionality makes it clear that $b^{\prime }=b.$\hfill 
$\square $\medskip

\begin{itemize}
\item The proof of Theorem 4.4 below is a variant of the proof of Theorem
4.1, and will strike the reader who has worked through the proof of Theorem
4.1 as repetitious, but the proof is a bit more involved here since we do
not have access to the axiom of choice in $\mathcal{M}$ or in $\mathcal{N}$.
Note that in Theorem 4.4 $\mathcal{N}$ is assumed to end extend $\mathcal{M}$%
, whereas in Theorem 4.1 $\mathcal{N}$ is only assumed to fix $\omega ^{%
\mathcal{M}}.$ Of course Theorem 4.4 follows from Theorem 4.1 if $\mathcal{M}
$ and $\mathcal{N}$ are models of $\mathsf{ZFC}$.\medskip
\end{itemize}

\noindent \textbf{4.4.~Theorem.}~\textit{Suppose} $\mathcal{M}$ \textit{and} 
$\mathcal{N}$ \textit{are models of} $\mathsf{ZF}$, \textit{and }$\mathcal{M}%
\subsetneq _{\mathrm{end,faithful}}\mathcal{N}$. \textit{Then:}\medskip

\begin{enumerate}
\item[\textbf{(a)}] \textit{There is some} $\gamma \in \mathrm{Ord}^{%
\mathcal{N}}\backslash \mathrm{Ord}^{\mathcal{M}}$ \textit{such that} $%
\mathcal{M}\preceq \mathcal{N}_{\gamma }.$ \textit{Thus} \textit{either} $%
\mathcal{M}=\mathcal{N}_{\gamma }$, \textit{in which case} $\mathcal{N}$ 
\textit{is a topped rank extension of} $\mathcal{M}$, \textit{or }$\mathcal{M%
}\prec \mathcal{N}_{\gamma }$.\medskip

\item[\textbf{(b)}] \textit{In particular, there is a full satisfaction
class }$S$\textit{\ for }$\mathcal{M}$\textit{\ such that }$S$\textit{\ can
be written as }$D\cap M$, \textit{where} $D$ \textit{is }$\mathcal{N}$%
\textit{-definable.}\medskip
\end{enumerate}

\noindent \textbf{Proof.} Assume $\mathcal{M}$ and $\mathcal{N}$ are as in
the assumptions of the theorem. By Lemma 4.3 we can assume that $\mathcal{N}$
rank extends $\mathcal{M}$. Fix some $\lambda \in \mathrm{Ord}^{\mathcal{N}%
}\backslash \mathrm{Ord}^{\mathcal{M}}$ and some ordinal $\beta >\lambda $.
For each $m\in M$, we can define the following set within $\mathcal{N}$%
:\smallskip

\begin{center}
$O_{m}:=\left\{ \gamma \in \beta :\mathcal{N}\models \gamma \in \mathrm{Def}(%
\mathrm{V}_{\beta },\in ,\lambda ,m)\right\} .$\smallskip
\end{center}

\noindent where $x\in \mathrm{Def}(\mathrm{V}_{\beta },\in ,\lambda ,m)$ is
shorthand for the formula of set theory that expresses:\medskip

\begin{center}
$x$ is definable in $(\mathrm{V}_{\beta },\in )$ using at most the
parameters $\lambda $ and $m.$\medskip
\end{center}

\noindent Clearly $\lambda \in O_{m}$ and $O_{m}$ is coded in $\mathcal{N}$.
Next let:\smallskip

\begin{center}
$O:=\bigcup\limits_{m\in M}O_{m}.$\smallskip
\end{center}

\noindent Let $O^{\ast }=O\backslash \mathrm{Ord}^{\mathcal{M}}$. Clearly $%
O^{\ast }$ is nonempty since $\lambda \in O^{\ast }.$ We now consider two
cases:\medskip

\begin{enumerate}
\item[(I)] $O^{\ast }$ has a least element (under $\in ^{\mathcal{N}}).$%
\medskip

\item[(II)] $O^{\ast }$ has no least element.\medskip
\end{enumerate}

\noindent Case I. Let $\gamma _{0}$ be the least element of $O^{\ast }.$ It
is easy to see that $\gamma _{0}$ is a limit ordinal in $\mathcal{N}$. We
claim that: \smallskip

\begin{center}
$\mathcal{M}\preceq \mathcal{N}_{\gamma _{0}}.$\smallskip
\end{center}

\noindent \medskip To verify the claim, by the Tarski criterion of
elementarity, it suffices to show that if $\mathcal{N}_{\gamma _{0}}\models
\exists y\varphi (y,m)$ for some $m\in M$, then there is some $m_{0}\in M$
such that $\mathcal{N}_{\gamma _{0}}\models \varphi (m_{0},m).$ Let $\gamma
_{m}$ be the first ordinal below $\gamma _{0}$ such that:\smallskip

\begin{center}
$\mathcal{N}_{\gamma _{0}}\models \exists y\in \mathrm{V}_{\gamma _{m}}\
\varphi (y,m).$\smallskip
\end{center}

\noindent So $\gamma _{m}$ is definable in $\mathcal{N}_{\beta }$ with
parameters $\gamma _{0}$ and $m$, and since $\gamma _{0}\in O$ by
assumption, this shows that $\gamma _{m}$ is definable in $\mathcal{N}%
_{\beta }$ with parameter $m.$ So $\gamma _{m}\in O.$ By the fact that $%
\gamma _{0}$ was chosen to be the least element of $O^{\ast }$, this shows
that $\gamma _{m}\in \mathrm{Ord}^{\mathcal{M}},$ which makes it clear that
there is some $m_{0}\in M$ such that:\smallskip

\begin{center}
$\mathcal{N}_{\gamma _{0}}\models \varphi (m_{0},m).$\smallskip
\end{center}

\noindent This concludes the verification of $\mathcal{M}\preceq \mathcal{N}%
_{\gamma _{0}}$. To finish the proof of Theorem 4.4 it suffices to show that
Case II\ is impossible. Suppose to the contrary that $O^{\ast }$ has no
least element. Within $\mathcal{N}$ let $s_{\beta }$ be the Tarskian
satisfaction class for $\left( \mathrm{V}_{\beta },\in \right) $. For each $%
m\in M$ we define: \smallskip

\begin{center}
$\Phi _{m}:=\{\varphi (c,\dot{m})\in M:\mathcal{N}\models $ $\varphi (\dot{%
\lambda},\dot{m})\in s_{\beta }\}.$\smallskip
\end{center}

\noindent Note the constant $c$ is interpreted as $\lambda $ in the
right-hand side of the above definition of $\Phi _{m}.$ Next we
define:\smallskip

\begin{center}
$\Phi :=\bigcup\limits_{m\in M}\Phi _{m}$.\smallskip
\end{center}

\noindent Observe that $\Phi \subseteq M$ and $\Phi $ is $\mathcal{M}$%
-amenable since $\mathcal{N}$ is a faithful rank extension of $\mathcal{M}$.
So intuitively speaking, $\Phi $ is the subset of the elementary diagram of $%
\mathcal{N}_{\beta }$ (as computed in $\mathcal{N}$) that consists of
sentences whose constants are in $\{c\}\cup \{\dot{m}:m\in M\}.$ Since $%
\mathcal{N}$ need not be $\omega $-standard, the elements of $\Phi $ need
not be standard formulae. Also note that the constant $c$ is interpreted as $%
\lambda $ in the right-hand side of the above definition of $\Phi _{m}.$ Now
let:\smallskip

\begin{center}
$\Gamma :=\left\{ 
\begin{array}{c}
t(c,\dot{m})\in M:t(c,\dot{m})\in \Phi ,\ \mathrm{and}\ \forall \theta \in 
\mathrm{Ord\ }\left( t(c,\dot{m})>\dot{\theta}\right) \in \Phi ,\ \mathrm{and%
}\  \\ 
t\ \mathrm{is\ a\ definable\ term\ in\ the\ language}\ \mathcal{L}_{\mathrm{%
set}}\cup \{c\}\cup \{\dot{m}:m\in M\}%
\end{array}%
\right\} ,$
\end{center}

\noindent So, officially speaking, $\Gamma $ consists of \textit{syntactic
objects} $\varphi (c,\dot{m},x)$ in $\mathcal{M}$ that satisfy the following
three conditions in $(\mathcal{M},\Phi )$:\medskip

\begin{enumerate}
\item[$(i)$] $\left[ \exists !x\varphi (c,\dot{m},x)\right] \in \Phi .$%
\medskip

\item[$(ii)$] $\left[ \forall x\left( \varphi (c,\dot{m},x)\rightarrow x\in 
\mathrm{Ord}\right) \right] \in \Phi .$\medskip

\item[$(iii)$] $\forall \theta \in \mathrm{Ord}\ \left[ \forall x\left(
\varphi (c,\dot{m},x)\rightarrow \dot{\theta}\in x\right) \right] \in \Phi .$%
\medskip
\end{enumerate}

\noindent Clearly $\Gamma $ is definable in $\left( \mathcal{M},\Phi \right) 
$. Since we are in Case II, $\left( \mathcal{M},\Phi \right) \models \psi $,
where $\psi $ is the sentence that expresses:\medskip

\begin{center}
$\Gamma \neq \varnothing \wedge \forall t\left[ t\in \Gamma \longrightarrow
\exists t^{\prime }\in \Gamma \ (t^{\prime }\in t)\in \Phi \right] .$\medskip
\end{center}

\noindent In contrast to the proof of Theorem 4.1 we cannot at this point
conclude that there is a countable descending chain $\left\langle t_{n}:n\in
\omega \right\rangle $ in $\mathcal{M}$ since $\mathsf{DC}$ need not hold in 
$\mathcal{M}$. Instead, we will use the following argument that takes
advantage of the Reflection Theorem. Choose $k\in \omega $ such that $\psi $
is a $\Sigma _{k}(\Phi )$-statement$,$ and use the Reflection Theorem in $%
\left( \mathcal{M},\Phi \right) $ to pick $\mu \in \mathrm{Ord}^{\mathcal{M}%
} $ such that:\smallskip

\begin{center}
$\left( \mathcal{M}_{\mu },\Phi \cap M_{\mu }\right) \prec _{\Sigma _{k}(%
\mathrm{X})}\left( \mathcal{M},\Phi \right) .$\footnote{%
Here $\mathrm{X}$ is interpreted by the left-hand side structure as $\Phi
\cap M_{\mu }$ and by the right-hand side structure as $\Phi $.}\smallskip
\end{center}

\noindent Then $\psi $ holds in $\left( \mathcal{M}_{\mu },\Phi \cap M_{\mu
}\right) $. Within $\mathcal{M}$, let \smallskip

\begin{center}
$w=\{v\in \mathrm{V}_{\mu }:v\in \Gamma \}.$ \smallskip
\end{center}

\noindent Thus $\left( \mathcal{M},\Phi \right) \models $ $w\neq \varnothing
\wedge \exists t^{\prime }\in \Gamma \ (t^{\prime }\in t)\in \Phi .$ Observe
that since $\mathcal{N}$ has access to the Tarskian satisfaction class for $%
\mathcal{N}_{\beta }$, $\mathcal{N}$ can evaluate each term $t(c)$ in $w$ as
an ordinal $\delta <\beta $, where $\mathcal{N}_{\beta }\models t(\dot{%
\lambda})=\delta $, where $t(\dot{\lambda})$ is the term obtained by
replacing all occurrences of the constant $c$ with $\dot{\lambda}$ in $t(c).$
So we can consider $s\in N$, where:\smallskip

\begin{center}
$\mathcal{N}\models s=\left\{ \delta :\exists t(c)\in w\ (\mathrm{V}_{\beta
},\in )\models \delta =t(\dot{\lambda})\right\} .$\smallskip
\end{center}

\noindent Within $\mathcal{N}$, $s$ is a nonempty set of ordinals that has
no least element, which of course is a contradiction. This shows that Case
II is impossible, thus concluding the proof. \hfill $\square $\medskip

\noindent \textbf{4.5.~Remark.}~Arguing in $\mathsf{ZFC}$, suppose $\kappa $
is the first strongly inaccessible cardinal, $\lambda $ is the second
strongly inaccessible cardinal, $\mathcal{M}=(\mathrm{V}_{\kappa },\in )$
and $\mathcal{N}=(\mathrm{V}_{\lambda },\in )$. Then $\mathcal{N}$ is a
proper faithful rank extension of $\mathcal{M}$ in which the first clause of
the dichotomy of the conclusion of part (a) of Theorem 4.1 (and 4.4) holds ,
but not the second one. On the other hand, if $\kappa $ is a weakly compact
cardinal, then as noted by Kaufmann \cite[Proposition 2.3]{Kaufmann}, $(%
\mathrm{V}_{\kappa },\in )$ has a faithful topless elementary end extension $%
\mathcal{N}$; in this scenario the second clause of the dichotomy of the
conclusion of part (a) of Theorem 4.1 (and 4.4) holds, but not the first
one.\medskip

\noindent \textbf{4.6.~Remark.}~If the assumptions of Theorem 4.4 are
strengthened by adding the assumption that $\mathcal{M}\prec \mathcal{N}$,
then as in the proof of Theorem 3.3 of \cite{AliJSL1987} we can use
`Kaufmann's trick' to conclude that there is some $\mathcal{K}$ such that $%
\mathcal{M}\prec \mathcal{K}\preceq \mathcal{N}$ and $\mathcal{K}$ is a
topped rank extension of\textit{\ }$\mathcal{M}$. This yields a
strengthening of Theorem 3.3 of \cite{AliJSL1987} by eliminating the
assumption that the axiom of choice holds in $\mathcal{M}$ and $\mathcal{N}$%
.\medskip

\noindent \textbf{4.7.~Remark.}~An inspection of the proof of Theorem 4.1
makes it clear that the assumption that $\mathcal{M}\models \mathsf{ZFC}$
can be reduced to $\mathcal{M}\models \mathsf{ZF}+\mathsf{DC}$, and the
assumption that $\mathcal{N}\models \mathsf{ZFC}$ can be reduced to $%
\mathcal{N}\models \mathsf{KPR}$ + \textquotedblleft every set can be
well-ordered\textquotedblright . Here $\mathsf{KPR}$ is as in part (d) of
Definition 2.1. Similarly, an inspection of the proof of Theorem 4.4 shows
that the assumption that $\mathcal{N}\models \mathsf{ZF}$ can be reduced to $%
\mathcal{N}\models \mathsf{KPR}$. \medskip

\noindent \textbf{4.8.~Question.}~\textit{Can Theorem} 4.1 \textit{and} 
\textit{Corollary} 4.2 \textit{be strengthened by assuming that} $\mathcal{M}
$ \textit{and} $\mathcal{N}$ \textit{are models of} $\mathsf{ZF}$?\bigskip

\pagebreak

\begin{center}
\textbf{5}.~\textbf{DEAD-END MODELS}\bigskip
\end{center}

In this section we will establish Theorem C of the abstract, the proof uses
many ingredients, including the following one that refines a result obtained
independently by Kaufmann and the author who demonstrated Theorem 5.1 for
models of $\mathsf{ZF}$ in which the Axiom of Choice holds (see
\textquotedblleft Added in Proof\textquotedblright\ of \cite{Kaufmann} and
Remark 1.6 of \cite{Ali-TAMS}). \medskip

\noindent \textbf{5.1.~Theorem.}~\textit{No model of }$\mathsf{ZF}$ \textit{%
has a conservative proper end extension satisfying} $\mathsf{ZF}$. \medskip

\noindent \textbf{Proof.}~This follows from putting part (b) of Theorem 4.4
together with Tarski's Undefinability of Truth Theorem (as in Theorem
2.8).\hfill $\square .$\medskip

The rest of the section is devoted to presenting results that in conjunction
with Theorem 5.1 will allow us to establish Theorem C of the abstract (as
Theorem 5.18).\medskip

\noindent \textbf{5.2.~Definition.}~Suppose $\mathcal{M}$\textbf{\ }is an%
\textbf{\ }$\mathcal{L}_{\mathrm{set}}$-structure. \medskip

\begin{enumerate}
\item[\textbf{(a)}] $X\subseteq M$ is \textit{a class\footnote{%
Classes of $\mathcal{M}$ are sometimes referred to as \textit{piecewise coded%
} subsets of $\mathcal{M}$.} }in $\mathcal{M}$ if $\forall a\in M\ \exists
b\in M\ X\cap \mathrm{Ext}_{\mathcal{M}}(a)=\mathrm{Ext}_{\mathcal{M}}(b).$%
\medskip

\item[\textbf{(b)}] $\mathcal{M}$ is \textit{rather classless} if every
class of $M$ is $\mathcal{M}$-definable.\medskip

\item[\textbf{(c)}] $\mathcal{M}$ is $\aleph _{1}$-like if $\left\vert
M\right\vert =\aleph _{1}$ but $\left\vert \mathrm{Ext}_{\mathcal{M}%
}(a)\right\vert \leq \aleph _{0}$ for each $a\in M.$ \medskip
\end{enumerate}

\noindent \textbf{5.3.~Theorem.}~(Keisler-Kunen \cite{Keislertree}, Shelah 
\cite{Shelah}). \textit{Every countable model of} $\mathsf{ZF}$ \textit{has
an elementary end extension to an} $\aleph _{1}$-\textit{like rather
classless model}. \medskip

\noindent \textbf{5.4.~Proposition. }\textit{No rather classless model of} $%
\mathsf{ZF}$\ \textit{has a proper rank extension to a model of} $\mathsf{ZF}
$.\medskip

\noindent \textbf{Proof.} This follows from putting Theorem 4.4 together
with the observation that a rank extension of a rather classless model of $%
\mathsf{ZF}$\ is a conservative extension, and therefore a faithful
extension.\hfill $\square $\medskip

\noindent \textbf{5.5.~Remark.}~As pointed out in Remark 1.6 of \cite%
{AliJSL1987} it is possible for a rather classless model to have a proper
end extension satisfying\ $\mathsf{ZF}$, since $\aleph _{1}$-like rather
classless models exist by Theorem 5.3, and one can use the Boolean-valued
approach to forcing to construct set generic extensions of such
models.\medskip

\noindent \textbf{5.6.~Definition.}~A \textit{ranked tree} $\tau $ is a
two-sorted structure $\tau =(\mathbb{T},\ \leq _{\mathbb{T}},\ \mathbb{L},\
\leq _{\mathbb{L}},\ \rho )$ satisfying the following three
properties:\medskip

\begin{enumerate}
\item[(1)] $(\mathbb{T},\ \leq _{\mathbb{T}})$ is a tree, i.e., a partial
order such that any two elements below a given element are
comparable.\medskip

\item[(2)] \textbf{\ }$(\mathbb{L},\ \leq _{\mathbb{L}})$ is a linear order
with no last element.\medskip

\item[(3)] $\rho $ is an order preserving map from $(\mathbb{T},\ \leq _{%
\mathbb{T}})$ onto $(\mathbb{L},\ \leq _{\mathbb{L}})$ with the property
that for each $t\in \mathbb{T},$ $\rho $ maps the set of predecessors of $t$
onto the initial segment of $(\mathbb{L},\ \leq _{\mathbb{L}})$ consisting
of elements of $L$ less than $\rho (t).$\medskip
\end{enumerate}

\noindent \textbf{5.7.~Definition.}~Suppose $\tau =(\mathbb{T},\ \leq _{%
\mathbb{T}},\ \mathbb{L},\ \leq _{\mathbb{L}},\ \rho )$ is a ranked tree. A
linearly ordered subset $B$ of $\mathbb{T}$ is said to be a \textit{branch}
of $\tau $ if the image of $B$ under $\rho $ is $\mathbb{L}$. The \textit{%
cofinality} of $\tau $ is the cofinality of $(\mathbb{L},\ \leq _{\mathbb{L}%
}).$\medskip

\noindent \textbf{5.8.\ Definition.}~Given a structure $\mathcal{M}$ in a
language $\mathcal{L}$, we say that a ranked tree $\tau $ is $\mathcal{M}$%
-definable if $\tau =\mathbf{t}^{\mathcal{M}}$, where $\mathbf{t}$ is an
appropriate sequence of $\mathcal{L}$-formulae whose components define the
corresponding components of $\tau $ in $\mathcal{M}.$ $\mathcal{M}$ is 
\textit{rather branchless} if for each $\mathcal{M}$-definable ranked tree $%
\tau $, all branches of $\tau $ (if any) are $\mathcal{M}$-definable.
\medskip

\noindent \textbf{5.9.~Theorem}. \textit{Suppose} $\mathcal{M}$\textit{\ is
a countable structure in a countable language}.\medskip

\begin{enumerate}
\item[\textbf{(a)}] (Keisler-Kunen \cite{Keislertree}, essentially). \textit{%
It is a theorem of} $\mathsf{ZFC}+\Diamond _{\omega _{1}}$ \textit{that }$%
\mathcal{M}$ \textit{can be elementarily extended to a rather branchless
model}.\medskip

\item[\textbf{(b)}] (Shelah \cite{Shelah}). \textit{It is a theorem of} $%
\mathsf{ZFC}$ \textit{that} $\mathcal{M}$ \textit{can be elementarily
extended to a rather branchless model}.\medskip
\end{enumerate}

\noindent \textbf{5.10.~Definition.}~Suppose $(\mathbb{P},\leq _{\mathbb{P}%
}) $ is a poset (partially ordered set).\medskip

\begin{enumerate}
\item[\textbf{(a)}] $(\mathbb{P},\leq _{\mathbb{P}})$ is \textit{directed}
if any given pair of elements of $\mathbb{P}$ has an upper bound.\medskip

\item[\textbf{(b)}] A subset $F$ of $\mathbb{P}$ is a \textit{filter} over $(%
\mathbb{P},\leq _{\mathbb{P}})$ if the sub-poset $(F,\leq _{\mathbb{P}})$ is
directed.\medskip

\item[\textbf{(c)}] A filter over $(\mathbb{P},\leq _{\mathbb{P}})$ is 
\textit{maximal} if it cannot be properly extended to a filter over $(%
\mathbb{P},\leq _{\mathbb{P}})$.\footnote{%
Thus a maximal filter can be described as a \textit{maximally compatible}
subset of $(\mathbb{P},\leq _{\mathbb{P}}).$}\medskip

\item[\textbf{(d)}] A subset $C$ of $\mathbb{P}$\ is \textit{cofinal} in $(%
\mathbb{P},\leq _{\mathbb{P}})$ if $\forall x\in \mathbb{P}\ \exists y\in C\
x\leq _{\mathbb{P}}y.$\medskip
\end{enumerate}

\noindent \textbf{5.11.~Definition.}~Suppose $s$ is an infinite set.\medskip

\begin{enumerate}
\item[\textbf{(a)}] $[s]^{<\omega }$ is the directed poset of finite subsets
of $s,$ ordered by containment. Note that $[s]^{<\omega }$ is a directed set
with no maximum element.\medskip

\item[\textbf{(b)}] $\mathrm{Fin}(s,2)$ is the poset of finite functions
from $s$ into $\{0,1\}$, ordered by containment (where a function is viewed
as a set of ordered pairs). \medskip
\end{enumerate}

\noindent \textbf{5.12.~Example.}~Given an infinite set $a$, and $s\subseteq
a,$ let $\chi _{s}:s\rightarrow 2$ be the characteristic function of $s,$
i.e., $\chi _{s}(x)=1$ iff $x\in s.$ Clearly $[\chi _{s}]^{<\omega }$ is a
maximal filter of $\mathrm{Fin}(a,2).$ More generally, if $\mathcal{N}%
\models \mathsf{ZF}$, and $\mathcal{N}\models $ \textquotedblleft $%
s\subseteq a$ and $a$ is infinite\textquotedblright , consider $(\mathbb{P}%
,\leq _{\mathbb{P}})$, where:\smallskip\ 

\begin{center}
$\mathbb{P}:=\mathrm{Ext}_{\mathcal{N}}(\mathrm{Fin}^{\mathcal{N}}(a,2))$,
\smallskip
\end{center}

\noindent and $\leq _{\mathbb{P}}$ is set-inclusion (among members of $%
\mathbb{P}$) in the sense of $\mathcal{N}$. Then $F_{s}$ is a maximal filter
over $(\mathbb{P},\leq _{\mathbb{P}})$, where: \smallskip

\begin{center}
$F_{s}:=$ $\mathrm{Ext}_{\mathcal{N}}\left( [\chi _{s}]^{<\omega }\right) ^{%
\mathcal{N}}$.\smallskip
\end{center}

\noindent \textbf{5.13.~Definition.}~A structure $\mathcal{M}$ is a \textit{%
Rubin} model if it has the following two properties:\medskip

\begin{enumerate}
\item[\textbf{(a)}] Every~$\mathcal{M}$-definable directed poset $(\mathbb{D}%
,\leq _{\mathbb{D}})$ with no maximum element\textit{\ }has a cofinal chain
of length $\omega _{1}.$\medskip

\item[\textbf{(b)}] Given any $\mathcal{M}$-definable poset $\mathbb{P}$,
and any maximal filter $F$ of $\mathbb{P}$, if $F$ has a cofinal chain of
length $\omega _{1}$, then $F$ is coded in $\mathcal{M}$.\medskip
\end{enumerate}

\noindent \textbf{5.14.~Remark.}~If $\tau =(\mathbb{T},\ \leq _{\mathbb{T}%
},\ \mathbb{L},\ \leq _{\mathbb{L}},\ \rho )$ is a ranked tree, then each
branch of $\tau $ is a maximal filter over $(\mathbb{T},\ \leq _{\mathbb{T}%
}).$ This makes it clear that every Rubin~model is rather branchless. Also,
a rather branchless model of $\mathsf{ZF}$ is rather classless. To see this
consider the ranked tree $\tau _{\mathrm{class}}$ defined within a model $%
\mathcal{M}$ of $\mathsf{ZF}$ as follows: The nodes of $\tau _{\mathrm{class}%
}$ are ordered pairs $(s,\alpha )$, where $s\subseteq \mathrm{V}_{\alpha }$,
the rank of $(s,\alpha )$ is $\alpha $ and $(s,\alpha )<(t,\beta )$ if $%
\alpha \in \beta $ and $s=t\cap \mathrm{V}_{\alpha }$. It is easy to see
that $\mathcal{M}$ is rather classless iff every branch of $\tau _{\mathrm{%
class}}^{\mathcal{M}}$ is $\mathcal{M}$-definable.\footnote{%
The ranked tree $\tau _{\mathrm{class}}$ was first introduced by Keisler 
\cite[Example 2.1]{Keislertree}, who noted that $\mathcal{M}$ is rather
classless iff every branch of $\tau _{\mathrm{class}}^{\mathcal{M}}$ is $%
\mathcal{M}$-definable, which is the key property that we need here. As
noted by the referee, this key property is also satisfied by the subtree $%
\sigma _{\mathrm{class}}$ of $\tau _{\mathrm{class}}$ whose tree-nodes are
of the form $\left( \alpha ,s\right) $ where the ordinal rank of $s$ is $%
\alpha .$ Note that the branches of $\tau _{\mathrm{class}}^{\mathcal{M}}$
correspond to all classes of $\mathcal{M}$ (including improper ones, i.e.,
those that form a set), whereas the branches of $\sigma _{\mathrm{class}}^{%
\mathcal{M}}$ correspond to all proper classes of $\mathcal{M}$.} Hence we
have the following chain of implications:\medskip

\begin{center}
Rubin $\Rightarrow $ rather branchless $\Rightarrow $ rather
classless.\medskip
\end{center}

\noindent \textbf{5.15.~Theorem.}~(Rubin\footnote{%
The attribution to Rubin is informed by fact that the credit for \cite[Lemma
2.3]{Rubin-Shelah} is explicitly given to Rubin (at the end of the
introduction of the paper).} \cite[Corollary 2.4]{Rubin-Shelah}). \textit{%
It\ is a theorem of} $\mathsf{ZFC}+\Diamond _{\omega _{1}}$ \textit{that if }%
$\mathcal{M}$ \textit{is a countable structure in a countable language, then 
}$\mathcal{M}$ \textit{has an elementary extension of cardinality }$\aleph
_{1}$ \textit{that is a Rubin model}.\footnote{%
The proof of Theorem 5.15 in \cite[Lemma 2.3]{Rubin-Shelah} has a small gap
(it is assumed that $T_{0}\vdash \exists \vec{x}\psi (\vec{x})$ instead of
assuming the consistency of $T_{0}+\exists \vec{x}\psi (\vec{x})$). The
improved proof was presented in \cite[Theorem 2.1.3]{Ali-thesis} and is
reproduced here in the Appendix.} \medskip

\noindent \textbf{5.16.~Definition.}~Suppose $\mathcal{M}$ is a model of $%
\mathsf{ZF}$. $\mathcal{M}$ is \textit{weakly Rubin\ }if (a) and (b) below
hold:\medskip

\begin{enumerate}
\item[\textbf{(a)}] $\mathcal{M}$ is rather classless (this is the
asymptotic case of $(ii)$ below, if one could use $a=\mathrm{V}$).\medskip

\item[\textbf{(b)}] For every element $a$ of $\mathcal{M}$ that is infinite
in the sense of $\mathcal{M}$ the following two statements hold:\medskip

\begin{enumerate}
\item[$(1)$] $\left( [a]^{<\omega }\right) ^{\mathcal{M}}$ has a cofinal
chain of length $\omega _{1}.$\medskip

\item[$(2)$] If $F$ is a maximal filter of $\mathrm{Fin}^{\mathcal{M}}(a,2)$
and $F$ has a cofinal chain of length $\omega _{1}$, then $F$ is coded in $%
\mathcal{M}$.\medskip
\end{enumerate}
\end{enumerate}

\noindent As demonstrated in the Appendix, Schmerl's strategy of $\Diamond
_{\omega _{1}}$-elimination in \cite{schmerldiamond} (which is based on an
absoluteness argument first presented by Shelah \cite{Shelah}) can be
employed to build weakly Rubin models within $\mathsf{ZFC}$. Thus we
have:\medskip

\noindent \textbf{5.17.~Theorem}~(Rubin-Shelah-Schmerl)\textit{\ It is a
theorem of }$\mathsf{ZFC}$\textit{\ that every countable model of }$\mathsf{%
ZF}$ \textit{has an elementary extension to a weakly Rubin model of
cardinality }$\aleph _{1}$\textit{.}\medskip

\noindent We are finally ready to establish the main result of this section
on the existence of `dead-end' models. Note that since every generic
extension is an end extension (even for ill-founded models), the model $%
\mathcal{M}$ in Theorem 5.18 has no proper generic extension.\medskip

\noindent \textbf{5.18.~Theorem.}~\textit{Every countable model }$\mathcal{M}%
_{0}\models \mathsf{ZF}$ \textit{has an elementary extension }$\mathcal{M}$ 
\textit{of power} $\aleph _{1}$ \textit{that has no proper end extension to
a model }$\mathcal{N}\models \mathsf{ZF}$. \textit{Thus every consistent
extension of }$\mathsf{ZF}$ \textit{has a model of power} $\aleph _{1}$ 
\textit{that has no proper end extension to a model of }$\mathsf{ZF}$%
.\medskip

\noindent \textbf{Proof}\footnote{%
The proof of Fact $\left( \nabla \right) $ presented here is a reformulation
of a proof that was suggested by the referee. Our original proof admittedly
used too much machinery.}.~We begin with a basic fact that will be called
upon in the proof.\medskip

\noindent \textbf{Fact} $\left( \nabla \right) .$ \textit{Suppose} $\mathcal{%
M}\models \mathsf{ZF}$ \textit{and} $\mathcal{N}\models \mathsf{ZF}$\textit{%
\ with} $\mathcal{M}\subseteq _{\mathrm{end}}\mathcal{N}$, \textit{and} $%
a\in M$. \textit{If} $s\in N$, $s$ \textit{is finite as viewed in }$\mathcal{%
N}$, \textit{and} $\mathcal{N}\models s\subseteq a$, \textit{then} $s\in M.$ 
\textit{Thus for all} $a\in M$, \textit{we have}:\smallskip

\begin{center}
$\left( \left[ a\right] ^{<\omega }\right) ^{\mathcal{M}}=\left( \left[ a%
\right] ^{<\omega }\right) ^{\mathcal{N}}.$\smallskip
\end{center}

\noindent \textbf{Proof. }The assumptions on $\mathcal{M}$ and $\mathcal{N}$
readily imply:\medskip

\noindent (1) $\omega ^{\mathcal{M}}=\omega ^{\mathcal{N}}$, and \medskip

\noindent (2) $\mathcal{M}\preceq _{\Delta _{1}}\mathcal{N}$. \medskip

\noindent Let $x=\mathrm{Fn}(y,a)$ be the usual formula expressing
\textquotedblleft $x$ is the set of all functions $f:y\rightarrow a$%
\textquotedblright $.$ It is easy to see that: \medskip

\noindent (3) The predicate $\left( x=\mathrm{Fn}(y,a)\wedge y\in \omega
\right) $ is $\Delta _{1}$ within $\mathsf{ZF.}$ \medskip

\noindent Observe that $x=\mathrm{Fn}(y,a)$ is clearly expressible by a $\Pi
_{1}$-formula within $\mathsf{ZF}$ (with no restriction on $y$). With the
added condition that $y\in \omega $, we can take advantage of recursion to
express $x=\mathrm{Fn}(y,a)$ by the following $\Sigma _{1}$-formula:\medskip

\begin{center}
$\exists \left\langle s_{0},\cdot \cdot \cdot ,s_{y}\right\rangle $ $\left[
s_{y}=x\wedge s_{0}=\{\varnothing \}\wedge \forall z<y\ (s_{z+1}=\left\{
f\cup \{(z,v)\}:f\in s_{z}\wedge v\in a\right\} )\right] .$\medskip
\end{center}

\noindent Now assume $s\in N$, $\mathcal{N}\models s\subseteq a$, and for
some $k\in N,$ $\mathcal{N\models }$ $k\in \omega \wedge \left\vert
s\right\vert =k$. Thanks to (1) and the assumption $\mathcal{M}\subseteq _{%
\mathrm{end}}\mathcal{N}$, $k\in M.$ Therefore by (2) and (3), if $b\in M$
such that $\mathcal{M}\models b=\mathrm{Fn}(k,a)$, then $\mathcal{N}\models
b=\mathrm{Fn}(k,a)$. Thus if $f\in N$ such that $f:k\rightarrow a$ is an
injective function such that $\mathrm{range}(f)=s$, then $f\in b$, and since 
$\mathcal{N}$ is an end extension of $\mathcal{M}$, $f\in M.$ Hence $\mathrm{%
range}(f)=s\in M.$\hfill $\square $ Fact $\left( \nabla \right) $) \medskip

Given a countable model $\mathcal{M}_{0}$ of $\mathsf{ZF}$, by Theorem 5.17
there is a weakly Rubin model $\mathcal{M}$ that elementarily extends $%
\mathcal{M}_{0}$. By Theorem 5.1, to prove Theorem 5.18 it is sufficient to
verify that every end extension $\mathcal{N}$ of $\mathcal{M}$ that
satisfies $\mathsf{ZF}$ is a conservative extension. Towards this goal,
suppose $\mathcal{M}\subsetneq _{\mathrm{end}}\mathcal{N}\models \mathsf{ZF}$%
. In light of Remark 2.6(e), it suffices to show that $\mathcal{M}\subsetneq
_{\mathrm{end}}^{\mathcal{P}}\mathcal{N}$. \medskip 

To show $\mathcal{M}\subsetneq _{\mathrm{end}}^{\mathcal{P}}\mathcal{N}$,
suppose that for $a\in M$ and $s\in N$, $\mathcal{N}\models s\subseteq a.$
We will show that $s\in M.$ By Fact $(\nabla )$ we may assume that $a$ is
infinite as viewed from $\mathcal{M}$. Note that this implies that $\mathcal{%
M}$ views $[a]^{<\omega }$ as a directed set with no maximum element. Fact $%
(\nabla )$ assures us that: \medskip

\noindent $(\ast )$ $\ \forall m\in M$ $\left( [m]^{<\omega }\right) ^{%
\mathcal{M}}=\left( [m]^{<\omega }\right) ^{\mathcal{N}},$ and \medskip

\noindent Since $\mathrm{Fin}(m,2)\subseteq \lbrack m\times
\{0,1\}]^{<\omega }$, $(\ast )$ implies:\medskip

\noindent $(\ast \ast )$ $\ \forall m\in M$ $\mathrm{Fin}^{\mathcal{N}}(m,2)=%
\mathrm{Fin}^{\mathcal{M}}(m,2)$.\medskip

\noindent Let $\chi _{s}\in N$ be the characteristic function of $s$ in $%
\mathcal{N}$, i.e., as viewed from $\mathcal{N}$, $\chi _{s}:a\rightarrow
\{0,1\}$ and $\forall x\in a(x\in s\leftrightarrow \chi _{s}(x)=1)$. Let
\medskip

\begin{center}
$F_{s}:=$ $\mathrm{Ext}_{\mathcal{N}}\left( [\chi _{s}]^{<\omega }\right) ^{%
\mathcal{N}}.$\medskip
\end{center}

\noindent As noted in Example 5.12, $F_{s}$ is a maximal filter over $%
\mathrm{Ext}_{\mathcal{N}}(\mathrm{Fin}^{\mathcal{N}}(a,2))$, so by $(\ast
\ast )$ and the assumption that $\mathcal{N}$ end extends $\mathcal{M}$, $%
F_{s}$ is a maximal filter over $\mathrm{Ext}_{\mathcal{M}}(\mathrm{Fin}^{%
\mathcal{M}}(a,2))$. Note that directed set $\mathrm{Ext}_{\mathcal{M}%
}\left( [a]^{<\omega }\right) ^{\mathcal{M}}$ has a cofinal chain $%
\left\langle p_{\alpha }:\alpha \in \omega _{1}\right\rangle $ thanks to the
assumption that $\mathcal{M}$ is weakly Rubin. Together with $(\ast )$, this
readily implies that $F_{s}$ has a cofinal chain $\left\langle q_{\alpha
}:\alpha \in \omega _{1}\right\rangle ,$ where:\smallskip

\begin{center}
$q_{\alpha }:=(\chi _{s}\upharpoonright p_{\alpha })^{\mathcal{N}}=(\chi
_{s}\upharpoonright p_{\alpha })^{\mathcal{M}}.$ \smallskip
\end{center}

\noindent Therefore, by the assumption that $\mathcal{M}$ is weakly Rubin, $%
F_{s}$ is $\mathcal{M}$-definable and thus coded in $\mathcal{M}$ by some $%
m\in M$. This makes it clear that $\mathcal{M}\models \chi _{s}=\cup m$, and
thus $s\in $ $M$ since $s$ is $\Delta _{0}$-definable from $\chi _{s}$. This
concludes the verification that $\mathcal{M}\subsetneq _{\mathrm{end}}^{%
\mathcal{P}}\mathcal{N}$, which as explained earlier, is sufficient to
establish Theorem 5.18. \hfill $\square $\medskip

An examination of the proof of Theorem 5.1, together with Remark 4.7, makes
it clear that Theorem 5.1 can be refined as follows: \medskip

\noindent \textbf{5.19.~Theorem}.~\textit{No model of} $\mathsf{KPR+\Pi }%
_{\infty }^{1}$-$\mathsf{DC}$ \textit{has a conservative proper end
extension to a model of} $\mathsf{KPR}$. \medskip

The proof of Theorem 5.18, together with Theorem 5.19 allows us to refine
Theorem 5.18 as follows:\medskip

\noindent \textbf{5.20.~Theorem}.~\textit{Every countable model }$\mathcal{M}%
_{0}\models \mathsf{KPR}+\Pi _{\infty }^{1}$-$\mathsf{DC}$ \textit{has an
elementary extension }$\mathcal{M}$ \textit{of power} $\aleph _{1}$ \textit{%
that has no proper end extension to a model }$\mathcal{N}\models \mathsf{KPR}%
\mathrm{.}$ \textit{Thus every consistent extension of }$\mathsf{KPR}+\Pi
_{\infty }^{1}$-$\mathsf{DC}$ \textit{has a model of power} $\aleph _{1}$ 
\textit{that has no proper end extension to model of }$\mathsf{KPR}\mathrm{.}
$\medskip

\noindent \textbf{5.21.~Question}.~\textit{Is there an }$\omega $-\textit{%
standard} \textit{model of }$\mathsf{ZF(C)}$ \textit{that has no proper end
extension to a model of} $\mathsf{ZF(C)}$?\medskip

\noindent \textbf{5.22.~Remark}.~Weakly Rubin models are never $\omega $%
-standard (even though $\omega $-standard rather classless models exist in
abundance). However, a slight variation of the existence proof of Rubin
models shows that a Rubin model $\mathcal{M}$ can be arranged to have an $%
\aleph _{1}$-like $\omega ^{\mathcal{M}}.$ Note, however, that the answer to
Question 5.21 is known to be in the negative if \textquotedblleft end
extension\textquotedblright\ is strengthened to \textquotedblleft rank
extension\textquotedblright ; this follows from of \cite[Theorem 1.5(b)]%
{Ali-TAMS}, (which asserts that no $\aleph _{1}$-like model of $\mathsf{ZFC}$
that elementarily end extends a model all of whose ordinals are pointwise
definable has a proper rank extension to a model of $\mathsf{ZFC}$. The
proof of \cite[Theorem 1.5(b)]{Ali-TAMS} together with part (b) of Theorem
4.4 makes it clear that the same goes for $\mathsf{ZF}$, i.e., there are $%
\omega $-standard models of $\mathsf{ZF}$ that have no proper rank extension
to a model of $\mathsf{ZF}$.\textbf{\bigskip }

\begin{center}
\textbf{6.~ADDENDUM AND CORRIGENDUM TO \cite{AliJSL1987}}\bigskip
\end{center}

Theorem 2.2 of \cite{AliJSL1987}\textbf{\ }incorrectly stated (in the
notation of the same paper, where Gothic letters are used for denoting
structures) that if $\mathfrak{A}\models \mathsf{ZF}$, $\mathfrak{A}\prec _{%
\mathrm{cons}}\mathfrak{B}$ and $\mathfrak{B}$ fixes each $a\in \omega ^{%
\mathfrak{A}},$ then $\mathfrak{A}$ is cofinal in $\mathfrak{B}$. As shown
in Theorem 3.1 of this paper, it is possible to arrange taller conservative
elementary extensions of models of $\mathsf{ZFC}$, thus Theorem 2.2 of \cite%
{AliJSL1987} is false as stated. However, by assuming that `slightly'
strengthening the conditions:\smallskip

\begin{center}
$\mathfrak{A}\models \mathsf{ZF}$, and $\mathfrak{B}$ fixes each $a\in
\omega ^{\mathfrak{A}},$\ \smallskip
\end{center}

\noindent to:

\begin{center}
$\mathfrak{A}\models \mathsf{ZFC}$, and $\mathfrak{B}$ fixes $\omega ^{%
\mathfrak{A}},$\ \smallskip 
\end{center}

\noindent one obtains a true statement, as indicated in Corollary 4.2 of
this paper.\footnote{%
To my knowledge, there is only one place in the published literature in
which Theorem 2.2 of \cite{AliJSL1987} has been used, namely in a recent
paper of Goldberg and Steel \cite{Goldberg+Steel}. In the discussion after
Lemma 3.5 of the aforementioned paper, it is noted that Theorem 2.2 of \cite%
{AliJSL1987} implies that an elementary embedding $\pi :\mathcal{M}%
\rightarrow \mathcal{N}$ (where $\mathcal{M}$ and $\mathcal{N}$ are inner
models of $\mathsf{ZFC}$) has the property that $\pi ^{-1}[S]\in M$ for each 
$S\in N$ (paraphrased in \cite{Goldberg+Steel} as `$\pi $ is a close
embedding of $\mathcal{M}$ into $\mathcal{N}$') iff for each parametrically
definable subset $D$ of $\mathcal{N}$, $\pi ^{-1}[D]$ is parametrically
definable in $\mathcal{M}$. In the context of \cite{Goldberg+Steel}, the
relevant models $\mathcal{M}$ and $\mathcal{N}$ are both $\omega $-models
(they are inner models, thus well-founded), so Corollary 4.2 of this paper
applies.} As pointed out in Remark 3.9, $\mathfrak{B}$ fixes each $a\in
\omega ^{\mathfrak{A}}$ whenever $\mathfrak{A}\prec _{\mathrm{cons}}%
\mathfrak{B}\models \mathsf{ZF}$.\medskip

Several corollaries were drawn from Theorem 2.2 of \cite{AliJSL1987}; what
follows are the modifications in them, necessitated by the above correction.
We assume that the reader has \cite{AliJSL1987} handy for ready reference.

\begin{itemize}
\item In Corollary 2.3(a) of \cite{AliJSL1987}, the assumption that $%
\mathfrak{B}$ fixes every integer of $\mathfrak{A}$\ should be strengthened
to: $\mathfrak{B}$ fixes $\omega ^{\mathfrak{A}}$.

\item Part (b) of Corollary 2.3 of \cite{AliJSL1987} is correct as stated in 
\cite{AliJSL1987}: if $\mathfrak{A}$ is a Rubin model of $\mathsf{ZFC}$ and $%
\mathfrak{A}\prec \mathfrak{B}$ that fixes every $a\in \omega ^{\mathfrak{A}%
} $, then $\mathfrak{B}$ is a conservative extension of $\mathfrak{A}$. This
is really what the proof of part (a) of Corollary 2.3 of \cite{AliJSL1987}%
\textbf{\ }shows. However, the following correction is needed: towards the
end of the proof of part (a) of Corollary 2.3 of \cite{AliJSL1987}, after
defining the set $X$, it must simply be stated as follows: But $X=Y\cap A,$
where $Y$ is the collection of members $a\in A$ satisfying $[\varphi
(b,a)\wedge a\in R_{\alpha }]$. Since $\mathfrak{B}$ is an end extension of $%
\mathfrak{C}$, $Y\in C,$ so by claim $(\ast )$, $X$\ is definable in $%
\mathfrak{A}$.

\item Corollary 2.4 of \cite{AliJSL1987} is correct as stated since it does
not rely on the full force of Theorem 2.2 of \cite{AliJSL1987} and only
relies on the fact that no elementary conservative extension of a model of $%
\mathsf{ZF}$ is an end extension. The latter follows from Corollary 4.2 of
this paper.

\item The condition stipulating that $\mathfrak{B}$\ fixes $\omega ^{%
\mathfrak{A}}$\ must be added to condition $(i)$ of Corollary 2.5 of \cite%
{AliJSL1987}\textbf{.}
\end{itemize}

\bigskip

\begin{center}
\textbf{APPENDIX: PROOFS\ OF\ THEOREMS 5.15 \& 5.17}

\bigskip
\end{center}

\noindent \textbf{A.1.} \textbf{Theorem}~(Rubin-Shelah-Schmerl)\textit{\ It
is a theorem of }$\mathsf{ZFC}$\textit{\ that every countable model of }$%
\mathsf{ZF}$ \textit{has an elementary extension to a weakly Rubin model of
cardinality }$\aleph _{1}$\textit{.}\medskip

\noindent The proof of Theorem A.1 has two distinct stages: in the first
stage we prove within $\mathsf{ZFC}+\Diamond _{\omega _{1}}$ that every
countably infinite structure in a countable language has an elementary
extension to a \textit{Rubin} model of cardinality\textit{\ }$\aleph _{1}$
(Theorem 5.15), and then in the second stage we use the result of the first
stage together with a forcing-and-absoluteness argument to establish the
theorem.\footnote{%
It appears to be unknown whether Theorem 5.15 (existence of Rubin elementary
extensions) can be proved in $\mathsf{ZFC}$ alone. However, as shown here, $%
\mathsf{ZFC}$ can prove that every countable model of $\mathsf{ZF}$ has an
elementary extension to a \textit{weakly Rubin model}. See also Remark A.7.}
\medskip

\begin{center}
\textbf{Stage 1 of the Proof of Theorem A.1}\medskip
\end{center}

\noindent Fix a $\Diamond _{\omega _{1}}$ sequence $\left\langle S_{\alpha
}:\alpha <\omega _{1}\right\rangle $. \textit{\ }Given a countable language $%
\mathcal{L}$, and countably infinite $\mathcal{L}$-structure model $\mathcal{%
M}_{0}$\textit{, }assume without loss of generality that $M_{0}=\omega .$ We
plan to inductively build two sequences $\left\langle \mathcal{M}_{\alpha
}:\alpha <\omega _{1}\right\rangle $, and $\left\langle \mathcal{O}_{\alpha
}:\alpha <\omega _{1}\right\rangle .$ The first is a sequence of
approximations to our final model $\mathcal{M}.$ The second sequence, on the
other hand, keeps track of the increasing list of `obligations'\ we need to
abide by throughout the construction of the first sequence. More
specifically, each $\mathcal{O}_{\alpha }$ will be of the form $%
\{\{V_{n}^{\alpha },\ W_{n}^{\alpha }\}:n\in \mathbb{\omega }\}$, where $%
\{V_{n}^{\alpha },\ W_{n}^{\alpha }\}$ is pair of disjoint subsets of $%
M_{\alpha }$ that are inseparable in $\mathcal{M}_{\alpha }$ and should be
kept inseparable in each $\mathcal{M}_{\beta }$, for all $\beta >\alpha .$
Here we say that two disjoint subsets $V$ and $W$ of a model $\mathcal{N}$
are \textit{inseparable} in $\mathcal{N}$ if there is no $\mathcal{N}$%
-definable $X\subseteq N$ such that $V\subseteq X,$ and $W\cap X=\varnothing
.\mathfrak{\medskip }$

We only need to describe the construction of these two sequences for stages $%
\alpha +1$ for limit ordinals $\alpha $ since:

\begin{itemize}
\item $\mathcal{O}_{0}:=\varnothing $.

\item For limit $\alpha $, $\mathcal{M}_{\alpha }:=\bigcup\limits_{\beta
<\alpha }\mathcal{M}_{\beta }$ and $\mathcal{O}_{\alpha
}:=\bigcup\limits_{\beta <\alpha }\mathcal{O}_{\beta }$.

\item For nonlimit $\alpha $ , $\mathcal{M}_{\alpha +1}:=\mathcal{M}_{\alpha
}$ and $\mathcal{O}_{\alpha +1}:=\mathcal{O}_{\alpha }$.
\end{itemize}

\noindent At stage $\alpha +1,$ where $\alpha $ is a limit ordinal, we have
access to a model $\mathcal{M}_{\alpha }$ (where $M_{\alpha }\in \omega _{1}$%
), and a collection $\mathcal{O}_{\alpha }$ of inseparable pairs of subsets
of $M_{\alpha }$. We now look at $S_{\alpha }$, and consider two cases:
either $S_{\alpha }$ is parametrically undefinable in $\mathcal{M}_{\alpha
}, $ or not$.$ In the latter case we `do nothing', and define $\mathcal{M}%
_{\alpha +1}:=\mathcal{M}_{\alpha }$ and $\mathcal{O}_{\alpha +1}:=\mathcal{O%
}_{\alpha }.$ But if the former is true, we augment our list of obligations
via:\smallskip

\begin{center}
$\mathcal{O}_{\alpha +1}:=\mathcal{O}_{\alpha }\cup \{\left\{ S_{\alpha },\
M_{\alpha }\backslash S_{\alpha }\right\} \}.$\smallskip
\end{center}

\noindent Notice that if $S_{\alpha }$ is parametrically undefinable in $%
\mathcal{M}_{\alpha },$ then $\{S_{\alpha },\ M_{\alpha }\backslash
S_{\alpha }\}$ is inseparable in $\mathcal{M}_{\alpha }.$ Then we use Lemma
A.2 below to build an elementary extension $\mathcal{M}_{\alpha +1}$ of $%
\mathcal{M}_{\alpha }$ that satisfies the following three conditions:\medskip

\noindent (1)\ For each $\{V,\ W\}\in \mathcal{O}_{\alpha +1}$, $V$\ and $W$%
\ are inseparable in $\mathcal{M}_{\alpha +1}$.\medskip

\noindent (2)\ For every\textit{\ }$\mathcal{M}_{\alpha }$-definable
directed set $\mathbb{D}$ with no last element, there is some\textit{\ }$%
d^{\ast }\in \mathbb{D}^{\mathcal{M}_{\alpha +1}}$ such that for each $d\in 
\mathbb{D}^{\mathcal{M}_{\alpha }},$ $\mathcal{M}_{\alpha +1}\models \left(
d<_{\mathbb{D}}d^{\ast }\right) .$\medskip

\noindent (3)\ $M_{\alpha +1}=M_{\alpha }+\omega $ (ordinal
addition).\medskip

\noindent This concludes the description of the sequences $\left\langle 
\mathcal{M}_{\alpha +1}:\alpha <\omega _{1}\right\rangle $, and $%
\left\langle \mathcal{O}_{\alpha }:\alpha <\omega _{1}\right\rangle .$
\medskip

\noindent Let $\mathcal{M}:=\bigcup\limits_{\alpha <\omega _{1}}\mathcal{M}%
_{\alpha },$ and $\mathcal{O}:=\bigcup\limits_{\alpha <\omega _{1}}\mathcal{O%
}_{\alpha }$; note that:$\mathfrak{\medskip }$

\noindent $(4)$ For each $\{V,\ W\}\in \mathcal{O}$, $V$\ and $W$\ are
inseparable in $\mathcal{M}$.\medskip

\noindent We now verify that $\mathcal{M}$ is a Rubin model. Suppose, on the
contrary, that for some $\mathcal{M}$-definable partial order $\mathbb{P}$,
there is a maximal filter $F\subseteq \mathbb{P}$ that is not $\mathcal{M}$%
-definable, and $F$ contains a cofinal $\omega _{1}$-chain $E$. By a
standard L\"{o}wenheim-Skolem argument there is some limit $\alpha <\omega
_{1}$ such that: \medskip

\noindent $(5)$ $\ \ (\mathcal{M}_{\alpha },S_{\alpha })\prec (\mathcal{M}%
,F),$ where $S_{\alpha }=F\cap \alpha .$\medskip

\noindent In particular,\medskip

\noindent (6) $S_{\alpha }$ is an undefinable subset of $\mathcal{M}_{\alpha
}.$\medskip

\noindent Since $S_{\alpha }$ is a countable subset of some $\mathcal{M}$%
-definable partial order $\mathbb{P}$, and $E$ is chain of length $\omega
_{1}$ that is cofinal in $\mathbb{P}$, there is some $e\in \mathbb{P}$ such
that $p<_{\mathbb{P}}e$ for each $p\in S_{\alpha }.$ Fix some $\beta <\omega
_{1}$ such that $e\in M_{\beta }.$ We claim that the formula $\varphi
(x)=(x\in \mathbb{P})\rightarrow (x<_{\mathbb{P}}e)$ separates $S_{\alpha }$
and $\mathbb{P}^{M_{\alpha }}\backslash S_{\alpha }$ in $\mathcal{M}.$ This
is easy to see since if $q\in \mathbb{P}^{M_{\alpha }}\backslash S_{\alpha }$%
, then by (5), $S_{\alpha }$ is a maximal filter on $\mathbb{P}^{M_{\alpha
}} $, and therefore there is some $p\in S_{\alpha }$ such that $\mathcal{M}%
_{\alpha }\models $ \textquotedblleft $p$ and $q$ have no common upper
bound\textquotedblright , which by (5) implies the same to be true in $%
\mathcal{M}$, in turn implying that $\lnot (q<_{\mathbb{P}}e),$ as desired.
We have arrived at a contradiction since on one hand, based on (2) and (6),
the formula $\varphi (x)$ witnesses the separability of $S_{\alpha }$ and $%
M_{\alpha }\backslash S_{\alpha }$ within $\mathcal{M},$ and on the other
hand $\{S_{\alpha },\ M_{\alpha }\backslash S_{\alpha }\}\in \mathcal{O}$ by
(5), and therefore (4) dictates that $S_{\alpha }$ and $M_{\alpha
}\backslash S_{\alpha }$ are inseparable in $\mathcal{M}.$ Thus the proof of
Theorem A.1 will be complete once we establish Lemma A.2 below.\medskip

\noindent \textbf{A.2}. \textbf{Lemma. }\textit{Suppose }$\mathcal{M}$ 
\textit{is a countable} $\mathcal{L}$-\textit{structure, where} $\mathcal{L}$
\textit{is countable.} \textit{Let }$\{\{V_{n},W_{n}\}:n\in \mathbb{\omega }%
\}$ \textit{is a countable list of }$\mathcal{M}$\textit{-inseparable pairs
of subsets of }$M$, \textit{and} $\left\{ \left( \mathbb{D}_{n},\leq _{%
\mathbb{D}_{n}}\right) :n\in \mathbb{\omega }\right\} $ \textit{be an
enumeration of} $\mathcal{M}$-\textit{definable directed sets with no last
element. Then there exists an elementary extension }$\mathcal{N}$\textit{\
of }$\mathcal{M}$\textit{\ that satisfies the following two
properties:\medskip }

\begin{enumerate}
\item[\textbf{(a)}] $V_{n}$ \textit{and} $W_{n}$\textit{\ remain inseparable
in }$\mathcal{N}$ \textit{for all }$n\in \mathbb{\omega }$\textit{.\medskip }

\item[\textbf{(b)}] \textit{For each }$n\in \mathbb{\omega }$\textit{\ there
is some }$d_{n}\in N$ \textit{such that} $\mathcal{N}\models d_{n}\in 
\mathbb{D}_{n}$ (\textit{i.e.,} $d_{n}$ \textit{satisfies the formula that
defines} $\mathbb{D}_{n}$ \textit{in} $\mathcal{M}$), \textit{and for each} $%
m\in \mathbb{D}_{n},$ $\mathcal{N}\models \left( m<_{\mathbb{D}%
_{n}}d_{n}\right) .$\textit{\medskip }
\end{enumerate}

Let $\mathcal{L}^{+}$\ be the language $\mathcal{L}$\ augmented with
constants from $C_{1}=\{\dot{m}:m\in M\}$ and $C_{2}=\left\{ d_{n}:n\in 
\mathbb{\omega }\right\} $ (where $C_{1}$ and $C_{2}$ are disjoint), and
consider the $\mathcal{L}^{+}$-theory $T$ below: 
\begin{equation*}
T:=\mathrm{Th}(\mathcal{M},m)_{m\in M}+\left\{ d_{n}\in \mathbb{D}_{n}:n\in 
\mathbb{\omega }\right\} \cup \left\{ \dot{m}<_{\mathbb{D}_{n}}d_{n}:%
\mathcal{N}\models \left( m\in \mathbb{D}_{n}\right) ,n\in \mathbb{\omega }%
\right\} .
\end{equation*}

\noindent $T$ is readily seen to be consistent since it\ is finitely
satisfiable in $\mathcal{M}\mathfrak{.}$ Moreover, it is clear that if $%
\mathcal{N}\models T$, then $\mathcal{M}\prec \mathcal{N}$ and $\mathcal{N}$
satisfies condition (b) of the theorem. To arrange a model of $T$ in which
conditions (a) also holds requires a delicate omitting types argument.
First, we need a pair of preliminary lemmas. The proofs of Lemma A.3 and A.4
are routine and left to the reader. Note that they are each other's
contrapositive, so it is sufficient to only verify one of them.\textit{%
\medskip }

\noindent \textbf{A.3.} \textbf{Lemma.} \textit{The following two conditions
are equivalent for a sentence }$\varphi (d_{k_{1}},\cdot \cdot \cdot
,d_{k_{p}})$\textit{\ of }$\mathcal{L}^{+}$.\medskip

\begin{enumerate}
\item[$(i)$] $T\vdash \varphi (d_{k_{1}},\cdot \cdot \cdot ,d_{k_{p}})$%
\textit{.}\medskip

\item[$(ii)$] $\mathcal{M}\models \exists r_{1}\in \mathbb{D}_{k_{1}}\forall
s_{1}\in \mathbb{D}_{k_{1}}\cdot \cdot \cdot \exists r_{p}\in \mathbb{D}%
_{k_{p}}\forall s_{p}\in \mathbb{D}_{k_{p}}\left[ \left(
\bigwedge\limits_{1\leq i\leq p}r_{i}<_{\mathbb{D}_{k_{i}}}s_{i}\right)
\rightarrow \varphi (s_{1},\cdot \cdot \cdot ,s_{p})\right] .$\medskip
\end{enumerate}

\noindent \textbf{A.4.} \textbf{Lemma.} \textit{The following two conditions
are equivalent for a sentence }$\varphi (d_{k_{1}},\cdot \cdot \cdot
,d_{k_{p}})$\textit{\ of }$\mathcal{L}^{+}$.\medskip

\begin{enumerate}
\item[$(i)$] $T+\varphi (d_{k_{1}},\cdot \cdot \cdot ,d_{k_{p}})$\textit{\
is consistent.}\medskip

\item[$(ii)$] $\mathcal{M}\models \forall r_{1}\in \mathbb{D}_{k_{1}}\exists
s_{1}\in \mathbb{D}_{k_{1}}\cdot \cdot \cdot \forall r_{p}\in \mathbb{D}%
_{k_{p}}\exists s_{p}\in \mathbb{D}_{k_{p}}\left[ \left(
\bigwedge\limits_{1\leq i\leq p}r_{i}<_{\mathbb{D}_{k_{i}}}s_{i}\right)
\wedge \varphi (s_{1},\cdot \cdot \cdot ,s_{p})\right] .$\medskip
\end{enumerate}

\noindent We are now ready to carry out an omitting types argument to
complete the proof of Lemma A.2. For each formula $\psi (y,\vec{x})$ of $%
\mathcal{L}^{+}$, and each $n\in \omega $, consider the following $q$-type
formulated in the language $\mathcal{L}^{+}$, where $\vec{x}=\left\langle
x_{1},\cdot \cdot \cdot ,x_{q}\right\rangle $ is a $q$-tuple of
variables:\smallskip

\begin{center}
$\Sigma _{n}^{\psi }(\vec{x}):=\{\psi (\dot{a},\vec{x}):a\in V_{n}\}\cup
\{\lnot \psi (\dot{b},\vec{x}):b\in W_{n}\}.$\smallskip
\end{center}

\noindent Note that $\Sigma _{n}^{\psi }$ expresses that $\psi (y,\vec{x})$
separates $V_{n}$ and $W_{n}.$\medskip

\noindent \textbf{A.5. Lemma.} $\Sigma _{n}^{\psi }$\textit{\ is locally
omitted by }$T$\textit{\ for each formula }$\psi (y,\vec{x})$\textit{, and
each }$n\in \mathbb{\omega }$.\medskip

\noindent \textbf{Proof. }Suppose to the contrary that there is a formula $%
\theta (\vec{x},d_{k_{1}},\cdot \cdot \cdot ,d_{k_{p}})$ of $\mathcal{L}^{+}$
and some $n\in \omega $ such that (1) through (3) below hold:\medskip

\begin{enumerate}
\item[$(1)$] $T+\exists \vec{x}\ \theta (\vec{x},d_{k_{1}},\cdot \cdot \cdot
,d_{k_{p}})$\textit{\ }is consistent.\medskip

\item[$(2)$] For all $a\in V_{n},$ $T\vdash \forall \vec{x}\left[ \theta (%
\vec{x},d_{k_{1}},\cdot \cdot \cdot ,d_{k_{p}})\rightarrow \psi (\dot{a},%
\vec{x})\right] .$\medskip

\item[$(3)$] For all $b\in W_{n},$ $T\vdash \forall \vec{x}\left[ \theta (%
\vec{x},d_{k_{1}},\cdot \cdot \cdot ,d_{k_{p}})\rightarrow \lnot \psi (\dot{b%
},\vec{x})\right] .$\medskip
\end{enumerate}

\noindent Invoking Lemmas A.3 and A.4, (1) through (3) translate to $%
(1^{\prime })$ through $(3^{\prime })$ below: \medskip

\begin{enumerate}
\item[$(1^{\prime })$] $\mathcal{M}\models \forall r_{1}\in \mathbb{D}%
_{k_{1}}\exists s_{1}\in \mathbb{D}_{k_{1}}\cdot \cdot \cdot \forall
r_{p}\in \mathbb{D}_{k_{p}}\exists s_{p}\in \mathbb{D}_{k_{p}}\left[ \left(
\bigwedge\limits_{1\leq i\leq p}r_{i}<_{\mathbb{D}_{k_{i}}}s_{i}\right)
\wedge \exists \vec{x}\ \theta (\vec{x},\vec{s})\right] ,$ where $\vec{s}%
=\left\langle s_{1},\cdot \cdot \cdot ,s_{p}\right\rangle .$\medskip

\item[$(2^{\prime })$] For all $a\in V_{n}$, $\mathcal{M}\models \lambda
(a), $ where:%
\begin{equation*}
\lambda (\mathbf{v}):=\forall \vec{x}\ \exists r_{1}\in \mathbb{D}%
_{k_{1}}\forall s_{1}\in \mathbb{D}_{k_{1}}\cdot \cdot \cdot \exists
r_{p}\in \mathbb{D}_{k_{p}}\forall s_{p}\in \mathbb{D}_{k_{p}}\left[ \left(
\bigwedge\limits_{1\leq i\leq p}r_{i}<_{\mathbb{D}_{k_{i}}}s_{i}\right)
\rightarrow \left( \theta (\vec{x},\vec{s})\rightarrow \psi (\mathbf{v},\vec{%
x}\ )\right) \right] .
\end{equation*}

\item[$(3^{\prime })$] For all $b\in W_{n},\ \mathcal{M}\models \gamma (b)$,
where:%
\begin{equation*}
\gamma (\mathbf{w}):=\forall \vec{x}\ \exists r_{1}\in \mathbb{D}%
_{k_{1}}\forall s_{1}\in \mathbb{D}_{k_{1}}\cdot \cdot \cdot \exists
r_{n}\in \mathbb{D}_{k_{p}}\forall s_{p}\in \mathbb{D}_{k_{p}}\left[ \left(
\bigwedge\limits_{1\leq i\leq p}r_{i}<_{\mathbb{D}_{k_{i}}}s_{i}\right)
\rightarrow \left( \theta (\vec{x},\vec{s})\rightarrow \lnot \psi (\mathbf{w}%
,\vec{x})\right) \right] .
\end{equation*}
\end{enumerate}

\noindent Let $\Lambda :=\{m\in M:\mathcal{M}\models \lambda (m)\}$, and $%
\Gamma :=\{m\in M:\mathcal{M}\models \gamma (m)\}$, and observe that $%
V_{n}\subseteq \Lambda $ by $(2^{\prime })$ and $W_{n}\subseteq \Gamma $ by $%
(3^{\prime })$. To arrive at a contradiction, we will show that $\Lambda
\cap \Gamma =\emptyset $, which implies that $V_{n}$ and $W_{n}$ are
separable in $\mathcal{M}$. To this end, suppose to the contrary that for
some $m\in M,$ \medskip

\begin{enumerate}
\item[$(4)$] $\mathcal{M}\models \lambda (m)\wedge \gamma (m).$ \medskip
\end{enumerate}

\noindent Since each $\left( \mathbb{D}_{n},\leq _{\mathbb{D}_{n}}\right) $
is a directed set, (4) implies:\smallskip

\begin{enumerate}
\item[$(5)$] $\mathcal{M}\models \forall \vec{x}\ \exists r_{1}\in \mathbb{D}%
_{k_{1}}\forall s_{1}\in \mathbb{D}_{k_{1}}\cdot \cdot \cdot \exists
r_{p}\in \mathbb{D}_{k_{p}}\forall s_{p}\in \mathbb{D}_{k_{p}}\left[ \left(
\bigwedge\limits_{1\leq i\leq p}r_{i}<_{\mathbb{D}_{k_{i}}}s_{i}\right)
\rightarrow \left( \psi (m,\vec{x})\wedge \lnot \psi (m,\vec{x})\right) %
\right] .$\smallskip
\end{enumerate}

\noindent Recall that no $\left( \mathbb{D}_{n},\leq _{\mathbb{D}%
_{n}}\right) $ has a last element. This shows that (5) yields a
contradiction, so the proof is complete.

\hfill $\square $ (Lemma A.5)

$\medskip $

\noindent \textbf{Proof of Lemma A.2. }Put Lemma A.5 together with the
Henkin-Orey omitting types theorem \cite[Theorem 2.2.9]{Chang-Keisler} to
conclude that there exists a model $\mathcal{N}$ of $T$ that satisfies
properties (a) and (b). \hfill $\square $\medskip

\begin{center}
\textbf{Stage 2 of the Proof of Theorem A.1}\medskip
\end{center}

\noindent In this stage we employ the result of the first stage together
with some set-theoretical considerations to establish the
Rubin-Shelah-Schmerl Theorem. The main idea here is a variant of the one
used by Schmerl, \cite{schmerldiamond}, which itself based on a method of $%
\diamondsuit _{\omega _{1}}$-elimination introduced by Shelah \cite{Shelah}.
This stage has three steps.\medskip

\noindent \textbf{Step 1 of Stage 2.} It is well-known that there is an $%
\omega _{1}$-closed notion $\mathbb{Q}_{1}$ in the universe $\mathrm{V}$ of $%
\mathsf{ZFC}$ such that the $\mathbb{Q}_{1}$-generic extension $\mathrm{V}^{%
\mathbb{Q}_{1}}$ of $\mathrm{V}$ satisfies $\mathsf{ZFC}+\diamondsuit _{%
\mathbb{\omega }_{1}}$ \cite[Theorem 8.3]{Kunen-book} (note that $\aleph
_{1} $ is absolute between $\mathrm{V}$ and $\mathrm{V}^{\mathbb{Q}_{1}}$
since $\mathbb{Q}_{1}$ is $\omega _{1}$-closed). Thus by Step 1, given a
model $\mathcal{M}$ in the universe $\mathrm{V}$ of $\mathsf{ZF}$, the
forcing extension $\mathrm{V}^{\mathbb{Q}_{1}}$ satisfies \textquotedblleft
there is an elementary extension $\mathcal{N}$ of $\mathcal{M}$ such that $%
\mathcal{N} $ is a Rubin model\textquotedblright . In $\mathrm{V}^{\mathbb{Q}%
_{1}}$, since $\mathcal{N}$ is a Rubin model, for each $s\in N$ that is
infinite in the sense of $\mathcal{N}$, there is a subset $C_{s}=\left\{
c_{\alpha }^{s}:\alpha \in \omega _{1}\right\} $ of $\mathcal{N}$ such that $%
\left\langle c_{\alpha }^{s}:\alpha \in \omega _{1}\right\rangle $ is an $%
\omega _{1}$-chain that is cofinal in the directed set $\left( [s]^{<\omega
}\right) ^{\mathcal{M}}.$ \medskip

\begin{itemize}
\item Let $\mathbb{P}_{s}$ be the suborder of $\left( \mathrm{Fin}%
(s,2)\right) ^{\mathcal{M}}$ consisting of partial $\mathcal{M}$-finite
functions $f$ from $s$ into $2$ such that $\mathrm{dom}(f)\in C_{s}$. Note
that $\mathbb{P}_{s}$ is a tree-order (the predecessors of any given element
are comparable).

\item $\mathbb{P}_{s}$ can be turned into a ranked-tree $\tau _{s}$ by
defining $\rho (f)=\alpha $ if $\mathrm{dom}(f)=c_{\alpha }$. The ranked
tree $\tau _{s}$ has the key property that each maximal filter of $\left( 
\mathrm{Fin}(s,2)\right) ^{\mathcal{M}}$ is uniquely determined by a branch
of $\tau _{s}.$

\item Let $C$ be the binary predicate on $\mathcal{N}$ defined by $C(s,x)$
iff $x\in C_{s}$, and let \smallskip
\end{itemize}

\begin{center}
$\mathcal{N}^{+}=\left( \mathcal{N},C\right) $. \smallskip
\end{center}

\noindent Note that the $\mathcal{N}^{+}$-definable ranked trees include
trees of the form $\tau _{s}$, as well as the ranked tree $\tau _{\mathrm{%
class}}^{\mathcal{N}}$ (whose branches uniquely determine classes of $%
\mathcal{N}$, as noted in Remark 5.14).

\begin{itemize}
\item This gives us $\aleph _{1}$-many ranked trees `of interest' in $%
\mathrm{V}^{\mathbb{Q}_{1}},$ namely $\tau _{\mathrm{class}}^{\mathcal{N}}$
together with $\tau _{s}$ for each infinite set $s$ of $\mathcal{N}$.\medskip
\end{itemize}

\noindent \textbf{Step 2 of Stage 2. }The concept of a weakly specializing
function plays a key role in this step. Given a ranked tree $(\mathbb{T},\
\leq _{\mathbb{T}},\ \mathbb{L},\ \leq _{\mathbb{L}},\ \rho )$, we say that $%
f:\mathbb{T\rightarrow \omega }$ \textit{weakly specializes} $(\mathbb{T},\
\leq _{\mathbb{T}})$ if $f$ has the following property:\smallskip

\begin{center}
If $x\leq _{\mathbb{T}}y$ and $x\leq _{\mathbb{T}}z$ and $f(x)=f(y)=f(z),$
then $y$ and $z$ are $\leq _{\mathbb{T}}$-comparable.\smallskip
\end{center}

\noindent The following lemma captures the essence of a weakly specializing
function for our purposes.\medskip

\noindent \textbf{A.6. Lemma. }\textit{Suppose} $\tau =(\mathbb{T},\ \leq _{%
\mathbb{T}},\ \mathbb{L},\ \leq _{\mathbb{L}},\ \rho )$ \textit{is a ranked
tree, where} $(\mathbb{L},\ \leq _{\mathbb{L}})$ \textit{has cofinality} $%
\omega _{1}$, \textit{and }$f$ \textit{weakly specializes} $(\mathbb{T},\
\leq _{\mathbb{T}})$. \textit{Then every branch of} $\tau $ \textit{is
parametrically definable in the expansion} $\left( \tau ,f\right) $ \textit{%
of} $\tau .$\medskip

\noindent \textbf{Proof.} Suppose $B$ is a branch of $\tau $. Then since ($%
\mathbb{L},\ \leq _{\mathbb{L}})$ has cofinality $\omega _{1},$ there is a
subset $S$ of $B$ of order type $\omega _{1}$ that is unbounded in $B$.
Therefore there is some subset $B_{0}\ $of $B$ that is unbounded in $B$ such
that $f$ is constant on $B_{0}.$ Fix some element $b_{0}\in B_{0}$ and
consider subset $X$ of $\mathbb{T}$ that is defined in $\left( \tau
,f\right) $ by the formula:\smallskip

\begin{center}
$\varphi (x,b_{0}):=[x\leq _{\mathbb{T}}b_{0}\vee (b_{0}\leq _{\mathbb{T}%
}x\wedge f(b_{0})=f(x))].$\smallskip
\end{center}

\noindent By the defining property of weakly specializing functions, $X$ is
linearly ordered and is unbounded in $B$. This makes it clear that $B$ is
definable in $\left( \tau ,f\right) $ as the downward closure of $X$, i.e.,
by the formula $\psi (x,b_{0}):=\exists y\left[ \varphi (y,b_{0})\wedge
x\leq _{\mathbb{T}}y\right] .$\hfill $\square $\medskip

\noindent By a remarkable theorem, due independently to Baumgartner \cite%
{Baumgartnet PFA} and Shelah \cite{Shelah}, given any ranked tree $\tau $ of
cardinality $\aleph _{1}$ and of cofinality $\omega _{1}$, there is a c.c.c.
notion of forcing $\mathbb{Q}_{\tau }$ such that $\mathrm{V}^{\mathbb{Q}%
_{1}\ast \mathbb{Q}_{\tau }}:=\left( \mathrm{V}^{\mathbb{Q}_{1}}\right) ^{%
\mathbb{Q}_{\tau }}$ contains a function $f_{\tau }$ such that weakly
specializes $\tau .$\footnote{%
This result generalizes the fact, established by Baumgartner, Malitz, and
Reinhardt \cite{Baugartnet et al} that every $\omega _{1}$-Aronszajn tree
can be turned into a special Aronszajn tree in a c.c.c. forcing extension of
the universe.} Here $\mathbb{Q}_{\tau }$ consists of all finite attempts for
building a function that weakly specializes $\tau $ (ordered by inclusion).
Moreover, given any collection $\{\tau _{\alpha }:\alpha \in \kappa \}$ of
ranked trees of cofinality $\omega _{1}$, the \textit{finite support product}
$\mathbb{Q}$ of $\left\{ \mathbb{Q}_{\tau _{\alpha }}:\alpha \in \kappa
\right\} $ is also known to have the c.c.c. property\footnote{%
It is well-known that the finite support product of a given collection of
forcing notions $\left\{ \mathbb{P}_{\alpha }:\alpha \in \kappa \right\} $
has the c.c.c. property iff every finite sub-product of $\left\{ \mathbb{P}%
_{\alpha }:\alpha \in \kappa \right\} $ has the c.c.c. property.}, as
indicated in \cite[Corollary 3.3]{Chodunsky and Zapletal}\footnote{\cite[%
Corollary 3.3]{Chodunsky and Zapletal} is overtly about specializing $\omega
_{1}$-Aronszajn trees, but as shown by Baumgartner \cite[Section 7]%
{Baumgartnet PFA}, given a tree $\tau $ of height $\omega _{1}$ with at most 
$\aleph _{1}$ many branches, there is an Aronszajn subtree $\tau _{0}$ of $%
\tau $ of height $\omega _{1}$, and moreover, any function $f_{0}$ that
specializes $\tau _{0}$ can be canonically extended to a function $f$ that
specializes $\tau .$ For more detail, see the proof of \cite[Lemma 3.5]%
{Switzer-Kaufmann}.}. Note that for each $\alpha \in \kappa $, $\mathbb{Q}$
forces the existence of a function $f_{\tau _{\alpha }}$ that weakly
specializes $\tau _{\alpha }.$

\begin{itemize}
\item Thus by choosing $\{\tau _{\alpha }:\alpha \in \omega _{1}\}$ to be
the collection of ranked trees that are parametrically definable in $%
\mathcal{N}^{+}$, there is a function $g(s,x)$ in $\mathrm{V}^{\mathbb{Q}%
_{1}\ast \mathbb{Q}}:=\left( \mathrm{V}^{\mathbb{Q}_{1}}\right) ^{\mathbb{Q}%
} $ such that for each infinite set $s$ of $\mathcal{N}$, $g(s,x)$ weakly
specializes $\tau _{s}$ ($\tau _{s}$ was defined in Step 1\ of Stage 2), and
a function $f$ that weakly specializes $\tau _{\mathrm{class}}^{\mathcal{N}%
}. $
\end{itemize}

The next lemma shows that forcing with $\mathbb{Q}$ does not add new
branches to ranked trees in the ground model that have cofinality $\omega
_{1}$. Note that when the lemma below is applied to the case when \textrm{V}
is replaced by $\mathrm{V}^{\mathbb{Q}_{1}}$\medskip , it shows that the
ample supply of rather branchless models in $\mathrm{V}^{\mathbb{Q}_{1}}$
remain rather branchless after forcing with $\mathbb{Q}$, i.e., in $\mathrm{V%
}^{\mathbb{Q}_{1}\ast \mathbb{Q}}.$\medskip

\noindent \textbf{A.7. Lemma. }\textit{Suppose} $\tau =(\mathbb{T},\ \leq _{%
\mathbb{T}},\ \mathbb{L},\ \leq _{\mathbb{L}},\ \rho )$ \textit{is a ranked
tree, where} $(\mathbb{L},\ \leq _{\mathbb{L}})$ \textit{has cofinality} $%
\omega _{1}.$ \textit{If} $B$ \textit{is a branch of} $\tau $ \textit{in} $%
\mathrm{V}^{\mathbb{Q}}$, \textit{then} $B\in \mathrm{V}$. \medskip

\noindent \textbf{Proof.} Suppose $B$ is a branch of $\tau $ \textit{in} $%
\mathrm{V}^{\mathbb{Q}}$, and let $f$ be the $\mathbb{Q}$-generic function
in $\mathrm{V}^{\mathbb{Q}}$ that specializes $\tau .$ Since $\mathbb{Q}$
does not collapse $\omega _{1}$ and the rank order of $\tau $ is assumed to
have cofinality $\omega _{1}$ in $\mathrm{V}$, there is for some $k\in
\omega $ such that: \medskip

\noindent (1) $\ \ \mathrm{V}^{\mathbb{Q}}\models $ \textquotedblleft $%
f^{-1}\{k\}\cap B$ is unbounded in $B$\textquotedblright . \medskip

\noindent In $\mathrm{V}^{\mathbb{Q}},$ let $x_{0}\in B$ with $f(x_{0})=k.$
Coupled with (1) this shows that there is some $q_{0}\in \mathbb{Q}$ such
that:\medskip\ 

\noindent (2) $\ \ \ q_{0}\Vdash $ $\left[ f(x_{0})=k\wedge \exists y>_{%
\mathbb{T}}x_{0}(y\in B\wedge f(y)=k)\right] $.\medskip

\noindent Let $B_{0}$ be the subset of $\mathbb{T}$ defined in $\mathrm{V}$
as:

\begin{center}
$B_{0}=\left\{ y\in \mathbb{T}:\exists q\in \mathbb{Q\ }(q\supseteq
q_{0}\wedge q\Vdash \left[ y>_{\mathbb{T}}x_{0}\wedge y\in B\wedge f(y)=k%
\right] \right\} .$\smallskip
\end{center}

\noindent We wish to show:\medskip

\noindent $(\ast )\ \ \mathrm{V}^{\mathbb{Q}}\models $ $B_{0}$ is a cofinal
subset of $B$. \medskip

\noindent To verify $(\ast )$, suppose $b_{0}\in B$. By (1), for some $y\in 
\mathbb{T}$, $\mathrm{V}^{\mathbb{Q}}$ satisfies \textquotedblleft $y>_{%
\mathbb{T}}b_{0}$, $y\in B$, and $f(y)=k$\textquotedblright $.$ Therefore
there is some $\mathbb{Q}$-condition $q\supseteq q_{0}$ such that:$\medskip $

$q\Vdash \left[ y>_{\mathbb{T}}x_{0}\wedge y\in B\wedge f(y)=k\right]
.\medskip $

\noindent Since $B_{0}\in \mathrm{V}$, $(\ast )$ makes it clear that $B\in 
\mathrm{V}$ as well, since $B$ can be defined in $\mathrm{V}$ as the
downward closure of $B_{0}.$\hfill $\square $\medskip

\begin{itemize}
\item Let $\Omega $ be the collection consisting of the ranked tree $\tau _{%
\mathrm{class}}^{\mathcal{N}}$ together with the family of ranked trees $%
\left\{ \tau _{s}^{\mathcal{N}}:\mathcal{N}\models \left\vert s\right\vert
\geq \aleph _{0}\right\} .$ Since $\mathcal{N}$ is a Rubin model in $\mathrm{%
V}$, each of the trees in $\Omega $ are rather branchless in $\mathrm{V}$,
and Lemma A.7 assures us the ranked trees in $\Omega $ remain rather
branchless in $\mathrm{V}^{\mathbb{Q}}$. Lemma A.6, on the other hand,
assures us that each branch of the ranked trees in $\Omega $ has a simple
definition in $\mathrm{V}^{\mathbb{Q}}.$

\item Let $\mathcal{L}_{\mathrm{set}}^{+}$ be the finite language obtained
by augmenting the usual language $\mathcal{L}_{\mathrm{set}}$ of set theory
with extra symbols for denoting the binary relation $C$, the binary function 
$g,$ and the unary function $f$. Then let $\mathcal{L}$ be the countable
language that is the result of augmenting $\mathcal{L}_{\mathrm{set}}^{+}$
with constants for each element of $\mathcal{M}$. The content of the above
bullet item allows us to write a single sentence $\psi $ in $\mathcal{L}_{%
\mathsf{\omega }_{\mathsf{1}}\mathsf{,\omega }}(\mathsf{Q})$ (where $\mathsf{%
Q}$ is the quantifier \textquotedblleft there exist uncountably
many\textquotedblright ) that captures the following salient features of the
structure $\left( \mathcal{N}^{+},C,g,f\right) $ that hold in $\mathrm{V}^{%
\mathbb{Q}}$:\medskip
\end{itemize}

\begin{enumerate}
\item[$(1)$] The $\mathcal{L}_{\mathrm{set}}$-reduct of $\left( \mathcal{N}%
^{+},C,g,f\right) $ elementarily extends $\mathcal{M}$.$\medskip $

\item[$(2)$] For each infinite $s$, $\tau _{s}^{\mathcal{N}}$ has
uncountable cofinality, is rather branchless, and $g(s,x)$ weakly
specializes $\tau _{s}^{\mathcal{N}}$.$\medskip $

\item[$(3)$] $\tau _{\mathrm{class}}^{\mathcal{N}}$ has uncountable
cofinality, is rather branchless, and $f$ weakly specializes $\tau _{\mathrm{%
class}}^{\mathcal{N}}$.\medskip
\end{enumerate}

\noindent More specifically, $\psi $ is the conjunction of the following
sentences:\medskip

\begin{itemize}
\item $\psi _{1}=$ the conjunction of the countably many sentences in the
elementary diagram of $\mathcal{M}.$

\item $\psi _{2}=$ the sentence expressing (2) above.

\item $\psi _{3}=$ the sentence expressing (3) above.
\end{itemize}

\noindent We elaborate on how to write $\psi _{2}$, a similar idea is used
to write $\psi _{3}.$ The quantifier \textsf{Q} allows us to express that $%
\tau _{s}$ has uncountable cofinality, and the range of the specializing
function $g(s,x)$ is countable. To expresses the rather branchless feature
of $\tau _{s}$, let $\psi (x,y)$ be as in the proof of Lemma A.6. Note that
in the notation of Lemma A.6, for any $b\in \mathbb{T}$ $\{x\in \mathbb{T}%
_{s}:\psi (x,b)\}$ is linearly ordered by $\leq _{\mathbb{T}}.$ By Lemma
A.6, in $\mathrm{V}^{\mathbb{Q}},$ if $s$ is infinite and $B$ is a branch of 
$\tau _{s}$, then there is some $b_{0}\in B$ such that $\psi (x,b_{0})$
defines $B$. \ Together with the fact that $\tau _{s}$ is rather branchless
in $\mathrm{V}^{\mathbb{Q}}$ (thanks to Lemma A.7) the $\mathcal{L}_{\mathsf{%
\omega }_{\mathsf{1}}\mathsf{,\omega }}$-formula below (no need for the
quantifier $\mathsf{Q}$ here) expresses the rather branchless feature of $%
\tau _{s}$. In the formula below $\mathbb{T}_{s}$ is the set of the nodes of
the ranked tree $\tau _{s}$, and \textquotedblleft unbounded
rank\textquotedblright\ is short for \textquotedblleft unbounded $\tau _{s}$%
-rank\textquotedblright .\smallskip

\begin{center}
$\forall b\in \mathbb{T}_{s}\left[ \{x\in \mathbb{T}_{s}:\psi (x,b)\}\ 
\mathrm{has\ unbounded\ rank\ }\rightarrow \bigvee\limits_{\theta (x,y)\in 
\mathcal{L}_{\mathrm{set}}^{+}}\exists y\ \forall x\ \left( \psi
(x,b)\leftrightarrow \theta (x,y)\right) \right] .$\smallskip
\end{center}

\noindent \textbf{Step 3 of Stage 2.} Let $\psi $ be the sentence described
above. Since $\psi $ has a model in $\mathrm{V}^{\mathbb{Q}_{1}\ast \mathbb{Q%
}}$, and by Keisler's completeness theorem for\textit{\ }$\mathcal{L}%
_{\omega _{1},\omega }(\mathsf{Q})$\textit{\ }\cite{Jerry-L(Q)}, the
consistency of sentences in $\mathcal{L}_{\omega _{1},\omega }(\mathsf{Q})$
is absolute for extensions of $\mathrm{V}$ that do not collapse $\aleph
_{1}, $ we can conclude that $\psi $ has a model in $\mathrm{V}$ since $%
\aleph _{1} $ is absolute between $\mathrm{V}$ and $\mathrm{V}^{\mathbb{Q}%
_{1}\ast \mathbb{Q}}$ (recall that $\mathbb{Q}_{1}$ is $\omega _{1}$-closed
and $\mathbb{Q}$ is c.c.c.$).$ Thanks to Lemmas A.6 and A.7, the $\mathcal{L}%
_{\mathrm{set}}$-reduct of any model of $\psi $ is a weakly Rubin model that
elementarily extends $\mathcal{M}$ (and the cardinality of the model can be
arranged to be $\aleph _{1}$ by a L\"{o}wenheim-Skolem argument). This
concludes the proof of Theorem A.1.\hfill $\square $ \medskip

\noindent \textbf{A.7. Remark.} Schmerl \cite{schmerldiamond} used a similar
argument as the one used in the proof of Theorem A.1 to show that every
countable model $\mathcal{M}$ of $\mathsf{ZFC}$ has an elementary extension $%
\mathcal{N}$ of cardinality $\aleph _{1}$ such that every $\omega ^{\mathcal{%
N}}$-complete ultrafilter over $\mathcal{N}$ is $\mathcal{N}$-definable (and
thus coded in $\mathcal{N}$). The proof of Theorem A.1 can be dovetailed
with Schmerl's proof so as to show that every countable model $\mathcal{M}$
of $\mathsf{ZFC}$ has an elementary extension $\mathcal{N}$ of cardinality $%
\aleph _{1}$ that is weakly Rubin and also has the additional property that
every $\omega ^{\mathcal{N}}$-complete ultrafilter over $\mathcal{N}$ is $%
\mathcal{N}$-definable.

\textit{\bigskip }

\noindent \textsc{Ali Enayat, Department of Philosophy, Linguistics, and
Theory of Science, University of Gothenburg}, \textsc{\ Sweden.}

\noindent \texttt{ali.enayat@gu.se}\medskip

\end{document}